\documentclass[a4paper,10pt]{article}
\usepackage{times}
\usepackage[utf8x]{inputenc} 
\usepackage{geometry} 
\usepackage{graphicx} 
\usepackage{amsmath}
\usepackage{amssymb}
\usepackage{float}
\usepackage{fancyhdr} 
\usepackage{sectsty}
\allsectionsfont{\sffamily\mdseries\upshape} 
\usepackage[english]{babel}
\usepackage{subfigure}
\newtheorem{theorem}{Theorem}[section]
\newtheorem{lemma}[theorem]{Lemma}

\newcommand{\qed}{\nobreak \ifvmode \relax \else
      \ifdim\lastskip<1.5em \hskip-\lastskip
      \hskip1.5em plus0em minus0.5em \fi \nobreak
      \vrule height0.75em width0.5em depth0.25em\fi}
\newcommand{\lx}{\lambda^x}
\newcommand{\ly}{\lambda^y}
\newcommand{\ipmh}{i\pm \frac{1}{2}}
\newcommand{\jpmh}{j\pm \frac{1}{2}}
\newcommand{\imh}{i-\frac{1}{2}}
\newcommand{\ipt}{i+\frac{3}{2}}
\newcommand{\Z}{{\mathbb Z}} 
\newcommand{\be} {\begin{equation}}
\newcommand{\ee} {\end{equation}}
\newcommand{\sstar}{s_*}
\newcommand{\lw}{\lambda_w}
\newcommand{\lo}{\lambda_o}
\newcommand{\rw}{\rho_w}
\newcommand{\ro}{\rho_o}

\newcommand{\re}{\mathbb{R}}
\newcommand{\half}{\frac{1}{2}}
\newcommand{\jph}{j+\half}
\newcommand{\iph}{i+\half}
\newcommand{\jmh}{j-\half}

\newcommand{\dt}{\Delta t}
\newcommand{\dx}{\Delta x}
\newcommand{\dy}{\Delta y}
\newcommand{\ud}{\textrm{d}}
\newcommand{\df}[2]{\frac{\partial #1}{\partial #2}}
\newcommand{\dd}[2]{\frac{\ud #1}{\ud #2}}
\newcommand{\norm}[1]{\left\| #1 \right\|}
\newcommand{\bmu}{\bar{\mu}}
\title{Multicomponent polymer flooding in two dimensional oil reservoir simulation }
\author{Sudarshan Kumar K \thanks{TIFR Centre for Applicable Mathematics Bangalore, sudarshan@math.tifrbng.res.in}, Praveen C \thanks{TIFR Centre for Applicable Mathematics Bangalore, praveen@math.tifrbng.res.in} and G. D Veerappa Gowda \thanks{TIFR Centre for Applicable Mathematics Bangalore,
 gowda@math.tifrbng.res.in}}
\date{} 
\begin{document}
\maketitle
\begin{abstract}
 We propose a high resolution finite volume scheme for a $(m+1)\times(m+1) $ system of nonstrictly  hyperbolic conservation laws which models 
 multicomponent polymer flooding in enhanced oil-recovery process in two dimensions. In the presence of gravity the flux functions need not be 
 monotone and hence the exact Riemann problem is complicated and computationally expensive.
 To overcome this difficulty, we use the idea of discontinuous flux to reduce the coupled system into uncoupled system of scalar conservation 
 laws with discontinuous coefficients. High order accurate scheme is constructed by introducing slope limiter
in space variable and a strong stability preserving Runge-Kutta scheme in the time variable. The performance of the 
 numerical scheme is  presented in various situations by choosing a heavily heterogeneous hard rock type medium. Also the significance of 
 dissolving  multiple polymers in aqueous phase is presented.
 \par
 For the updated vesion of this article please see\\
 Sudarshan Kumar, K.; Praveen, C.; Veerappa Gowda, G. D. A finite volume method for a two-phase multicomponent polymer flooding. J. Comput. Phys. 275 (2014), 667–695. See the link

 http://www.sciencedirect.com/science/article/pii/S0021999114004951\\

\end{abstract}
\section{Introduction}
Simulation of two phase flow in porous media plays a key role in many engineering areas such as 
oil-recovery \cite{aziz,branet,peaceman}, environmental remediation \cite{bear}
and water management in polymer electrolyte fuels cells \cite{lister}.
We are interested in multi dimensional simulation of two phase flow in heterogeneous porous media
arising in enhanced oil-recovery. It involves simultaneous flow of two immiscible phases (the aqueous phase and the oil phase) in a 
heterogeneous porous medium. We have assumed that $m$ chemical components are dissolved in the aqueous phase. These components 
could, for example, be different polymers that all have different influence on the flow properties. 
We propose a high order finite volume scheme for the numerical simulation of Buckley-Leverett model
with multicomponent polymer flooding by using the idea of discontinuous numerical flux developed in \cite{adidflu,dflu}.
For simplicity we let $\Omega=[0,1] \times [0,1]$ denote the two dimensional reservoir. Let $s \in [0,1]$ denote the saturation of aqueous phase
and $c=(c_1,c_2,....,c_m) \in [0, c_0]^m$ denote the concentration of the polymers dissolved in the aqueous phase, where $c_0$ is some  non negative real number. 
Then in the absence of capillary pressure the governing equations 
form a  $(m+1)\times(m+1)$ system of hyperbolic conservation laws ~\cite{thormord-3,thormord-1} given by
\begin{equation}
\begin{array}{rll}
s_t + \nabla \cdot F(s,c_1,c_2,.....c_m,x) &=& 0 \\
(sc_l+a_l(c_l))_t + \nabla \cdot (c_lF(s,c_1,c_2,.....c_m,x)) &=& 0,  \,\,\,l=1,2,....,m
\end{array}
\label{eq:system}
\end{equation}
where $ (x,t)\in \Omega\times(0,\infty) ,\,\,a_l:[0,1]\rightarrow \re $ are given smooth functions and  the flux $F:[0,1] \times [0, c_0]^m \times \Omega \rightarrow  \re^2$ is given by
$F=(F_1,F_2),$
\begin{equation}
F_1(s,c,x) = v_1(x)f(s,c), \quad f(s,c) = \frac{\lw(s,c)}{\lw(s,c) + \lo(s)}
\label{eq:f}
\end{equation}
\begin{equation}
F_2(s,c,x) = [v_2(x) - (\rw-\ro)g \lo(s,c) K(x) ] f(s,c).
\label{eq:f1}
\end{equation}\\
Here $\rw, \ro$ are the densities of water and oil, $g$ is the acceleration due to gravity. 
The quantities $\lw$ and $\lo$ are the mobilities of the water and oil phase respectively and $v=(v_1,v_2) \in \re^2$ is the total velocity given by Darcy law \cite{ewing1983}.
\begin{equation}
v = -\left((\lw + \lo) K(x) \frac{\partial p}{\partial x_1},(\lw + \lo) K(x) \frac{\partial p}{\partial x_2}+ (\lw \rho_w + \lo \rho_o)g K(x)\right)
\label{eq:pr}
\end{equation}
where $K : \Omega \to [0,\infty)$ is the permeability of the rock which can be discontinuous in $x $ and $p : \Omega \to \re$ is the pressure. If we assume incompressibility of the flow and if there are no sources, then the velocity is governed by
\begin{equation}
\nabla \cdot v = 0 \qquad \textrm{in } \quad \Omega
\label{eq:v}
\end{equation}
with some suitable boundary conditions  for  pressure on $\partial \Omega$. For instance in the inlet part of the boundary, water is pumped in at high pressure $p=p_I$ while a lower pressure $p=p_O$
is maintained on outlet, see Fig.\ref{fig:bc}. On the remaining part of the boundary, the normal velocity is set to zero, which gives a Neumann boundary condition
on pressure. Equations (\ref{eq:system}) and (\ref{eq:v}) form a system of coupled algebraic-differential equations and there is no time derivative  involved in 
equation (\ref{eq:v}). A commonly used model for the mobilities are
\begin{equation}
\lw(s,c) = \frac{s^2}{\mu_w(c)}, \qquad \lo(s) = \frac{(1-s)^2}{\mu_o}
\label{eq:mobility}
\end{equation}
where $\mu_w$, $\mu_o$ are the viscosities of water and oil and $\mu_w=\mu_w(c)$ which is increasing in each of its variable  $c_i.$
The term $a_l$ in (\ref{eq:system}) models the adsorption of the component $l$ on the porous medium.

\par
 In the absence of polymer flooding or equivalently if the flux function is independent of $c,$
then this problem (\ref{eq:system}) reduces to   scalar equation. In  \cite{karlsen-2}  by using a fast marching method and in \cite{karlsen-1} by using semi-Godunov scheme
method the problem is studied in the absence of polymer. Also in \cite{masson-3} two-phase flow problems are studied by using gradient schemes.
It is well known that in the 
heterogeneous media, that is when  the permeability $K(x)$ is discontinuous , fingering instability \cite{glimm} will develop and which results in an inefficient oil-recovery.
For example see Fig.\ref{fig1}(a). As the 
concentration $c$ increases, viscosity of  water increases and the fingering effects reduces which leads to an efficient oil-recovery
see Fig.\ref{fig1}(b). In the presence of the concentration $c$ the system (\ref{eq:system}) becomes coupled and non-strictly hyperbolic . When the 
concentration $c$ is smooth, existence and uniqueness theory is established in \cite{winther} but we deal here with the case when $c$  need not be  smooth. For this system, developing 
a Godunov type upwind schemes are difficult as it needs a solution of Riemann problems. Most often numerical methods requires the calculation of eigenvalues
and eigenvectors of the Jacobian matrix of the system.
Here by using the idea of discontinuous flux  we reduce the   system to an
uncoupled  scalar equations with discontinuous coefficients. Next we study each scalar equation by using the idea 
of discontinuous flux. This approach does not require detailed information about the eigenstructure of the full system.
Also in \cite{karlsen-1}, the idea of discontinuous flux is used to study a coupled  system arising in three-phase flows in porous media and shown its successfulness. 
Scalar conservation laws with discontinuous flux have been studied by many authors \cite{kyoto,burger,burger-1,chavent,diel1,diel2,gimse,jafrrey2,kaa,langtangen}. In particular, in \cite{dflu} a Godunov type  finite volume scheme is proposed
and convergence to a proper entropy solution is proved, provided the flux functions satisfies certain conditions like in \S \ref{1d-1}.
In one dimensional case for a $(2 \times 2)$ system this problem was studied in \cite{adidflu} and  there
proposed a finite volume scheme and named  numerical flux as DFLU. This DFLU flux works even in cases where the upstream mobility gives an entropy violating solution \cite{Mishra}.
Here we are extending DFLU to a multi dimensional case with high order accuracy.
The difficulties of developing an upwind type numerical schemes  in a highly heterogeneous media
in the presence of gravity attracts the importance of the proposed work.
\par
  The paper is organized as follows. 
From  \S \ref{1d-1} the idea of discontinuous flux for one dimensional problem is briefly explained and also  numerical experiments for high order schemes are performed 
to show their efficiency. In \S 3 two dimensional problem is introduced and 
the idea of the one dimensional discontinuous flux  is extended. Also high order accurate scheme is constructed by introducing slope limiter in
space variable and a strong stability preserving Runge-Kutta 
scheme in the time variable \cite{rk}. The resulting schemes are shown to respect a maximum principle. Also two dimensional
numerical results in various situation are shown for a quarter five-spot geometry. 
\section{System of equations in one dimension}\label{1d-1}
\begin{figure}
\centering
 \includegraphics[scale=0.5]{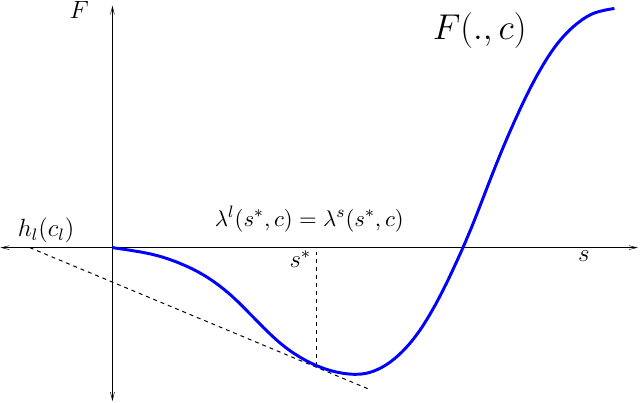}
 \caption{$\lambda^l=\lambda^s$} 
 \label{eigen}
\end{figure}

 The corresponding$(m+1)\times (m+1) $ system of equations in one-dimension in the presence of gravity is given by
\begin{equation}
\begin{array}{rll}
 s_t + \frac{\partial }{\partial x}F(s,c_1,c_2,...,c_m,x)&=& 0 \\
(sc_l+a_l(c_l))_t + \frac{\partial }{\partial x}c_lF(s,c_1,c_2,....,c_m,x) &=& 0, \quad l=1,2,...,m
\end{array}
\label{system1-D}
\end{equation}
where $t>0$ and $x\in \re, (s,c_1,c_2,....,c_m)=(s,c) \in [0,1] \times [0 ,c_0]^m $   and 
\begin{equation}
F(s,c,x) = [v - (\rw-\ro)g \lo(s,c) K(x) ] f(s,c).
\label{fluxfn}
\end{equation}
 In one dimension the solution  $v $ of the equation (\ref{eq:v}) reduces to a constant. 
We assume that  the flux function satisfies following conditions:
\begin{enumerate}
 \item $F(0,c_1,c_2,...,c_m,x)=0, \,\,F(1,c_1,c_2,....,c_m,x)=v \,\,\forall\,\, x, c_l, \quad l=1,2,...,m$
\item The function $s \rightarrow F(s,c_1,c_2,...,c_m,x) $ is of convex type i.e, has no local maximum in the interior of $  [0,1] \times [0 ,c_0]^m $ 
see Fig.\ref{eigen}
\item The adsorption term $a_l=a_l(c_l)$ satisfies
$h_l(c_l)=\frac{da_l}{dc_l}(c_l)>0, \forall\, c_l\in [0,1].$
\end{enumerate}
The case when $v=0$ and  $ F$ does not change sign  is studied in \cite{adidflu}. Here we assume  $v$ need not be zero and allow $F$
to change   sign, see Fig.~\ref{eigen}. 
In the absence of gravity, $F_s=v f_s$ is non-negative or non-positive depending on $v\ge 0$ or $v \le 0.$ Hence $F$ is increasing 
or decreasing in $s$ accordingly. In the presence of gravity $F_s$  becomes,
\begin{equation}
\nonumber
F_s = \frac{2s(1-s)}{\mu_w \mu_o (\lw+\lo)^2} [ v + (\rw-\ro)gK(x)(s \lw - (1-s)\lo)]
\end{equation}
which vanishes at $s=0,1$. Depending on the values of $v,\rw,\ro,g,K$, there can be a root $\sstar\in (0,1)$ which makes $F$ non-monotone in $s$, as shown
in Fig.\ref{eigen}. If such a root exists,
it is a root of the following cubic equation
\begin{equation}
\nonumber
r(c) s^3 - (1-s)^3 + z = 0, \quad r(c) = \frac{\mu_o}{\mu_w(c)}, \quad z = \frac{v \mu_o}{(\rw-\ro)gK}.
\end{equation}
This cubic equation has one real and two complex roots, the real root is given by
$$
s_*=\frac{1}{1+r}\left[ 1 -\frac{3\sqrt[3]{2} r}{\left(\alpha+\sqrt{\beta}\right)^{1/3}}+\frac{1}{3\sqrt[3]{2}}\left(\alpha+\sqrt{\beta}\right)^{1/3} \right]
$$
where
$$
\alpha = -27 r+27 r^2-27 z-54 r z-27 r^2 z, \qquad \beta = 2916 r^3 + \alpha^2.
$$
Since 
\begin{equation}
\nonumber
F_{ss}(\sstar,c,x) = \frac{6\sstar(1-\sstar)(\rw-\ro)gK(x)}{\mu_w \mu_o (\lw+\lo)^2} \left[ \frac{\sstar^2}{\mu_w} + \frac{(1-\sstar)^2}{\mu_o} \right]
\end{equation}
then $F$ attains the maximum(minimum) at $s=\sstar$ if $\rw > \ro$ ($\rw < \ro$). Note that the nature of the extremum depends only on the densities and is independent of the polymer 
concentrations 
$c_l$ and the permeability $K$.

\par
If $F(s,c,x)=F(s,c)$ then the system (\ref{system1-D}) can be put in the  matrix form as 
$$U_t+A(U)U_x=0, \quad  U=\begin{bmatrix} s \ \ c_1 \ \ c_2 \ \hdots \ c_m \end{bmatrix}^\top,$$
where $A(U)$ is the $(m+1)\times(m+1)$ Jacobian matrix
\[A(U)=
 \begin{pmatrix}
\frac{\partial F}{\partial s}& \frac{\partial F}{\partial c_1}&\frac{\partial F}{\partial c_2} & \cdots &\cdots& \frac{\partial F}{\partial c_n}\\
 0 &\frac{F}{s+h_1}& 0&\cdots&\cdots&0\\
 0&0&\frac{F}{s+h_2}&0&\cdots&0 \\
 \vdots&&\ddots&\ddots&&\vdots\\
 \vdots&&\ddots&\ddots&&\vdots\\
 0&\cdots&\cdots&\cdots&0& \frac{F}{s+h_m}
 \end{pmatrix}
 \]
The eigenvalues of this system are given by $$\lambda^s =\lambda(s,c)= \frac{\partial F}{\partial s}(s,c) $$ 
$$ \lambda^l =\lambda^l(s,c)= \frac{F(s,c)}{s+h_l(c_l)}, \quad l=1,2,..,m.$$ 
We can observe that for any $c=(c_1,c_2,...,c_m)\in [0,1]\times[0,c_0]^m $ and for some $l\in\{1,2,...,m\}$ 
there exist at least one point $ s^*=s^*(c)\in [0,1] $ 
such that  (see Fig.\ref{eigen}).
$$\lambda^l(s^*,c)=\lambda^s(s^*,c) .$$ 
For this couple $ (s^*,c), \lambda^l=\lambda^s, $ hence eigenvalues may coincide and the problem 
is non strictly hyperbolic.
The Rankine-Hugoniot condition corresponding to  
(\ref{system1-D}) is given by 
\begin{equation}
  \begin{array}{rll}
   F(s^R,c^R,x^+)-F(s^L,c^L,x^-)&=&\sigma(s^R-s^L)\\
   c_l^RF(s^R,c^R,x^+)-c_l^LF(s^L,c^L,x^-)&=&\sigma(s^Rc_l^R+a_l(c_l^R)-s^Lc_l^L-a_l(c_l^L)) \\
   &&\hspace{2 cm}\forall \,\, l=1,2,..,m.
  \end{array}
  \label{RHcondition1}
 \end{equation}
 For details see \cite{thormord-2,thormord-1}.
If $c^L=c^R$( i.e. $c_l^L=c_l^R \quad \forall \,\, l=1,2,...,m$)  then second equation reduces to  the first
equation of  (\ref{RHcondition1}). This corresponds to the Rankine-Hugoniot condition for single Buckely-Leverett equation (\ref{eq:system}).
Now we are interested in the case $c^L \neq c^R, $ i.e. $c_l^L\neq c_l^R$ for some $l, 1\leq l\leq m.$ If we combine the two equations 
(\ref{RHcondition1}) then we may write 
\begin{equation}
(c_l^R - c_l^L)F(s^L,c^L,x^-)=\sigma(c_l^R - c_l^L)s^L +\sigma (a_l(c_l^R)-a_l(c_l^L))
\label{rhc2}
\end{equation}
Define the functions $h^L_l$ by 
\begin{align}
 h^L_l(c_l)=\begin{cases}
  \frac{a_l(c_l)-a_l(c_l^L)}{c_l-c_l^L} & \mbox{ if }\quad c_l\neq c_l^L,\\
  h_l(c_l) & \mbox{ if } \quad c_l=c_l^L.
 \end{cases}
\end{align}
Now from (\ref{RHcondition1}) and (\ref{rhc2}), finally we get
\begin{equation}
\dfrac{F(s^R,c^R,x^+)}{s^R+\bar h}=\dfrac{F(s^L,c^L,x^-)}{s^L+\bar h }=\sigma,\\
\label{rh2}
\end{equation}
where $\bar h=h^L_l(c_l^R)\quad \left(=h_l^L(c_l^R) \,\,\forall\,\, l \right).$
Thus the Rankine-Hugoniot condition reduces to (\ref{rh2}).  This gives an idea how to obtain a weak solution of the Riemann 
problem to (\ref{system1-D}).
\subsection{Riemann problem}\label{rproblem}
\begin{figure}
\centering
\subfigure{
\includegraphics[scale=0.5]{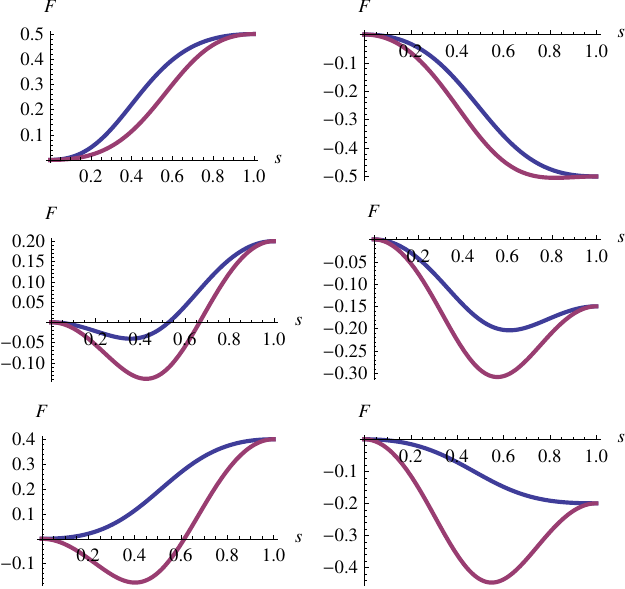}
}\\
\subfigure{
\includegraphics[scale=0.5]{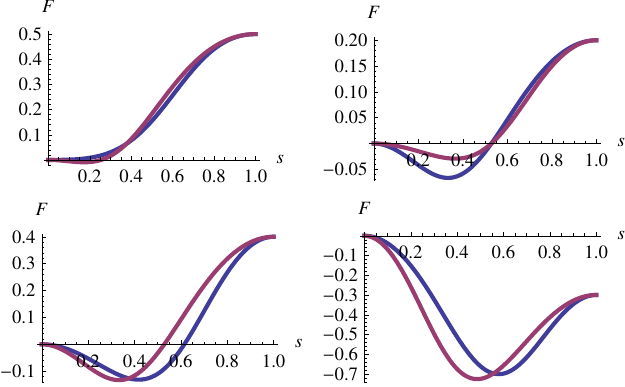}
}
\caption{Possible shapes of flux functions for different choices of $c_L$ and $c_R$.}
\label{flux}
\end{figure}

{\label{1d-3}} 
  For simplicity we restrict our study to the case when $m=2$  in equation (\ref{system1-D}), i.e $c=(c_1,c_2).$
  Also we assume that $F(s,c,x)=F(s,c).$
   Consider the Riemann problem associated to the system (\ref{system1-D}) with the initial condition 
\be s(x,0) = \left\{ \begin{array}{lll}
s_L &\mbox{if}& x<0, \\ s_R &\mbox{if}& x>0 \end{array} \right.  , \quad
c(x,0) = \left\{ \begin{array}{lll}
(c_1^L,c_2^L) &\mbox{if}& x<0, \\ (c_1^R,c_2^R) &\mbox{if}& x>0. \end{array} \right.
\label{riemannpb}
\ee
Solution to (\ref{system1-D}) and  (\ref{riemannpb}) is constructed by connecting states so that it should satisfies the Rankine-Hugoniot condition. There are two families of waves
that arise in the solution of the Riemann problem referred to as $ s $ and $ c $ waves. $ s $ waves consists of rarefaction and shocks (or contact discontinuity) across which s changes continuously and
 discontinuously respectively, but across which both $ c_1$ and $c_2$ remain constant. $ c $ waves consists solely of contact discontinuity across
 which both $s$ and $c_1, c_2$ changes such that $\frac{F}{s+\bar h}$ remains constant in the sense of (\ref{rh2}). 
For different choices of $c_L$ and $c_R$, the possible shapes of $F(s,c_L)$ and $F(s,c_R)$ are shown in
  Fig.\ref{flux}.
We restrict to the case when $c_L>c_R (i.e., \,c_l^L>c_l^R, l=1,2).$  When $c_L>c_R$  the flux functions
$s\rightarrow F(s,c^L) $ and $s\rightarrow F(s,c^R) $ are one of  the shapes  given in Fig.\ref{flux}. 
To explain the Riemann problem, for simplicity we consider the shape of the  flux functions as in Fig.\ref{rp-1}
\begin{itemize}
 \item Case 1:  $ s^L\leq s^*$\\
  Draw a line through the points $(-\bar h,0)$ and $(s^*,F(s^L,c^L_1,c^L_2)).$ This intersects the curve $F(s,c^R_1,c^R_2)$ at the point $\bar{s},$ where
  $F_s(\bar{s},c^R_1,c^R_2)\ge 0.$ We divide this in two subcases.
  \item Case 1a: $ s^R>\bar{s}$\\
  (a) Connect $(s^L,c^L_1,c^L_2)$ to $(s^*,c^L_1,c^L_2)$ by a s-rarefaction wave  (see Fig.\ref{rp-1}a).\\\\
  (b) Connect $(s^*,c^L_1,c^L_2)$ to $(\bar{s},c^R_1,c^R_2)$  by  a $c$-wave  with speed (see Fig.\ref{rp-1}a).
  $$\sigma_c=\frac{F(\bar{s},c^R_1,c^R_2)}{\bar{s}+\bar h}=F_s(s^*,c^L_1,c^L_2)$$
  (c) Connect $(\bar{s},c^R_1,c^R_2)$ to $(s^R,c^R_1,c^R_2)$ by a $s$-rarefaction wave (see Fig.\ref{rp-1}a).
   For example if $F(s,c^L_1,c^L_2)$ and $ F(s,c^R_1,c^R_2)$ are strictly convex functions then the corresponding solution of
   the Riemann problem is given by (see Fig.\ref{rp-1}b)
  \[ (s(x,t),c_1(x,t),c_2(x,t)) = \left\{ \begin{array}{lll}
(s^L,c^L_1,c^L_2) &\mbox{if}& x<\sigma_s t, \\ 
((F_s)^{-1}(\frac{x}{t},c^L_1,c^L_2),c^L_1,c^L_2) &\mbox{if}&  \sigma_s t < x < \sigma_c t,\\
(\bar{s},c^R_1,c^R_2) &\mbox{if}& \sigma_c t <x < \sigma_1 t,\\
((F_s)^{-1}(\frac{x}{t},c^R_1,c^R_2),c^R_1,c^R_2) &\mbox{if}&  \sigma_1 t < x < \sigma_2 t,\\
(s^R,c^R_1,c^R_2) &\mbox{if}&  x > \sigma_2 t.
 \end{array} \right.
\]
\begin{figure}[H]
\centering
\includegraphics[scale=0.35]{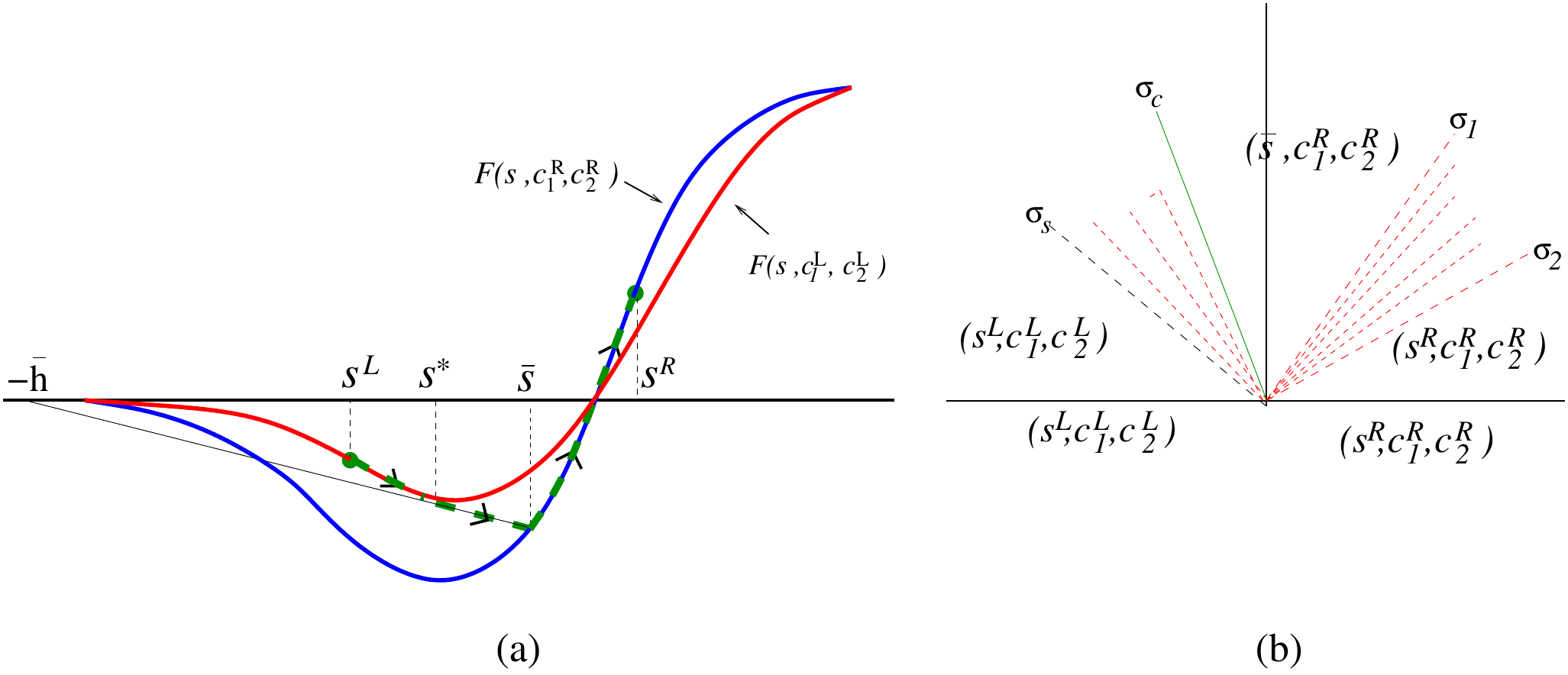}
\caption{Solution of the Riemann problem (\ref{riemannpb}) with $ s^L\leq s^*$ and  $ s^R>\bar{s}.$ }
\label{rp-1}
\end{figure} 
 \item Case 1b: $  s^R\leq \bar{s}$\\
 (a) Connect $(s^L,c^L_1,c^L_2)$ to $(s^*,c^L_1,c^L_2)$ by a s-rarefaction wave  (see Fig.\ref{rp-2}a).\\\\
  (b) Connect $(s^*,c^L_1,c^L_2)$ to $(\bar{s},c^R_1,c^R_2)$  by  a $c$-wave  with speed (see Fig.\ref{rp-2}a).
  $$\sigma_c=\frac{F(\bar{s},c^R_1,c^R_2)}{\bar{s}+\bar h}=F_s(s^*,c^L_1,c^L_2)$$
 (c) Connect $(\bar{s},c^R_1,c^R_2)$ to $(s^R,c^R_1,c^R_2)$ by a $s$ - shock wave with speed (see Fig.\ref{rp-2}a).
   $$\sigma_s= \frac{F(\bar{s},c^R_1,c^R_2)-F(s^R,c^R_1,c^R_2)}{\bar{s}-s^R}$$
  In the case of convex fluxes we can write the solution of the Riemann problem as (see Fig.\ref{rp-2}b)
 \[ (s(x,t),c_1(x,t),c_2(x,t)) = \left\{ \begin{array}{lll}
(s^L,c^L_1,c^L_2) &\mbox{if}& x<\sigma_1 t, \\ 
((F_s)^{-1}(\frac{x}{t},c^L_1,c^L_2),c^L_1,c^L_2) &\mbox{if}&  \sigma_1 t < x < \sigma_c t,\\
(\bar{s},c^R_1,c^R_2) &\mbox{if}& \sigma_c t <x < \sigma_s t,\\
(s^R,c^R_1,c^R_2) &\mbox{if}&  x > \sigma_s t.
 \end{array} \right.
\]
\begin{figure}[H]
\centering
\includegraphics[scale=0.35]{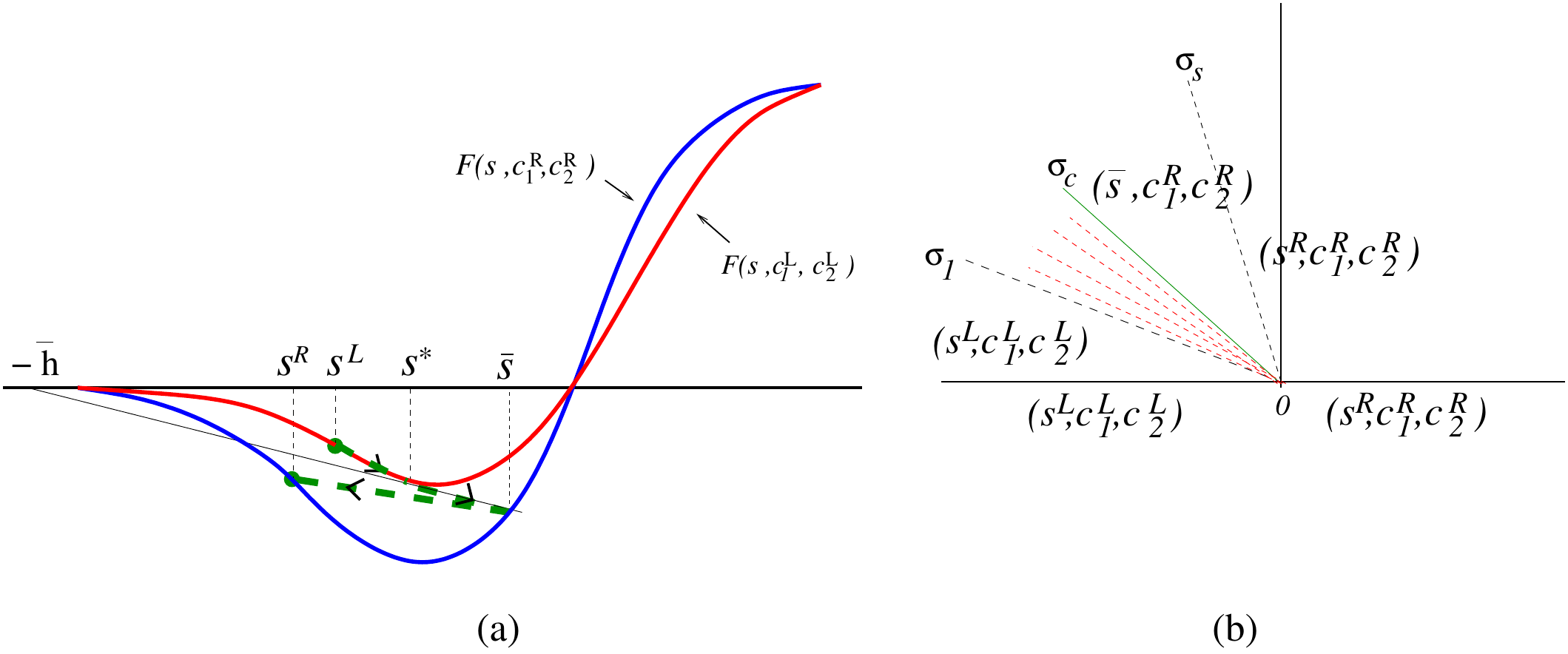}
\caption{Solution of the Riemann problem (\ref{riemannpb}) with $ s^L\leq s^*$ and  $ s^R\leq\bar{s}.$}
\label{rp-2}
\end{figure}
 \item Case 2: $ s^L>s^*$\\
  Draw a line joining the points $(-\bar h,0)$ and $(s^L,F(s^L,c^L_1,c^L_2)).$ Let $ (\bar{s},c^R_1,c^R_2)$ be the point where this line meets
  the curve $F(s,c^R_1,c^R_2),$  where $F_s(\bar{s},c^R_1,c^R_2)\geq 0.$ Consider the following subcases.
  \item Case 2a: $s^R\geq\bar{s}$\\
  (a) Connect $(s^L,c^L_1,c^L_2)$ to $(\bar{s},c^R_1,c^R_2)$ by a $c$- shock wave with speed (see Fig.\ref{rp-3}a).
    $$\sigma_c= \frac{F(s^L,c^L_1,c^L_2)-F(\bar{s},c^R_1,c^R_2)}{s^L-\bar{s}}$$
   (b) Connect $(\bar{s},c^R_1,c^R_2)$ to $(s^R,c^R_1,c^R_2)$ by a $s$- rarefaction wave (see Fig.\ref{rp-3}a). \\
In the case of convex flux the solution of the Riemann problem is given by (see Fig.\ref{rp-3}b)
\[ (s(x,t),c_1(x,t),c_2(x,t)) = \left\{ \begin{array}{lll}
(s^L,c^L_1,c^L_2) &\mbox{if}& x<\sigma_c t, \\ 
(\bar{s},c^R_1,c^R_2) &\mbox{if}& \sigma_c t <x < \sigma_1 t,\\
((F_s)^{-1}(\frac{x}{t},c^R_1,c^R_2),c^R_1,c^R_2) &\mbox{if}&  \sigma_1 t < x < \sigma_2 t,\\
(s^R,c^R_1,c^R_2) &\mbox{if}&  x > \sigma_2 t.
 \end{array} \right.
\]  
  \begin{figure}[H]
\centering\includegraphics[scale=0.35]{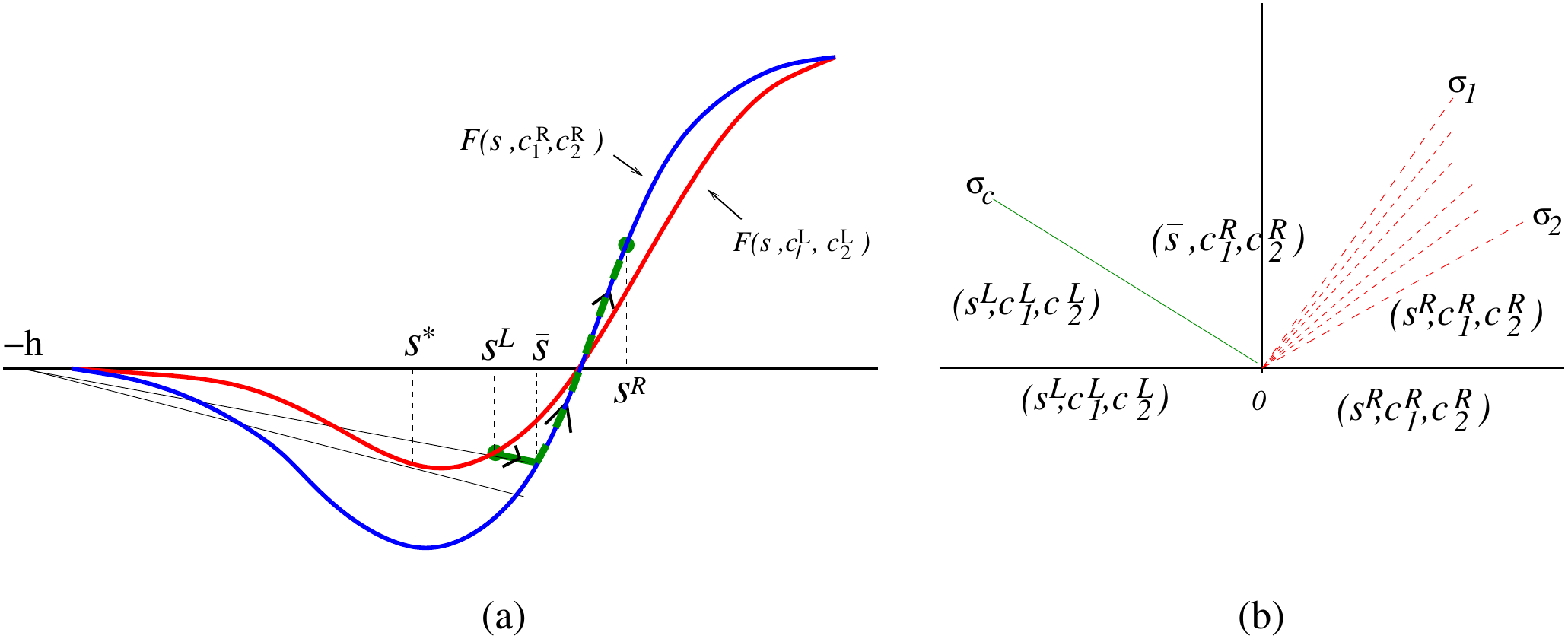}
\caption{Solution of the Riemann problem (\ref{riemannpb}) with $ s^L> s^*$ and  $ s^R\geq \bar{s}.$}
\label{rp-3}
\end{figure}
   \item Case 2b: $s^L<\bar{s}$\\
  (a) Connect $(s^L,c^L_1,c^L_2)$ to $(\bar{s},c^R_1,c^R_2)$ by a $c$- shock wave with speed (see Fig.\ref{rp-4}a).
    $$\sigma_c= \frac{F(s^L,c^L_1,c^L_2)-F(\bar{s},c^R_1,c^R_2)}{s^L-\bar{s}}$$
    (b) Connect $(\bar{s},c^R_1,c^R_2)$ to $(s^R,c^R_1,c^R_2)$ by a $s$-shock  wave with speed (see Fig.\ref{rp-4}a).\\
   $$\sigma_s= \frac{F(\bar{s},c^R_1,c^R_2)-F(s^R,c^R_1,c^R_2)}{\bar{s}-s^R}$$
In the case of convex flux the solution of the Riemann problem is given by (see Fig.\ref{rp-4}b)
\[ (s(x,t),c_1(x,t),c_2(x,t)) = \left\{ \begin{array}{lll}
(s^L,c^L_1,c^L_2) &\mbox{if}& x<\sigma_c t, \\ 
(\bar{s},c^R_1,c^R_2) &\mbox{if}& \sigma_c t <x < \sigma_st,\\
(s^R,c^R_1,c^R_2) &\mbox{if}&  x > \sigma_s t.
 \end{array} \right.
\]
\begin{figure}[H]
\centering
\includegraphics[scale=0.35]{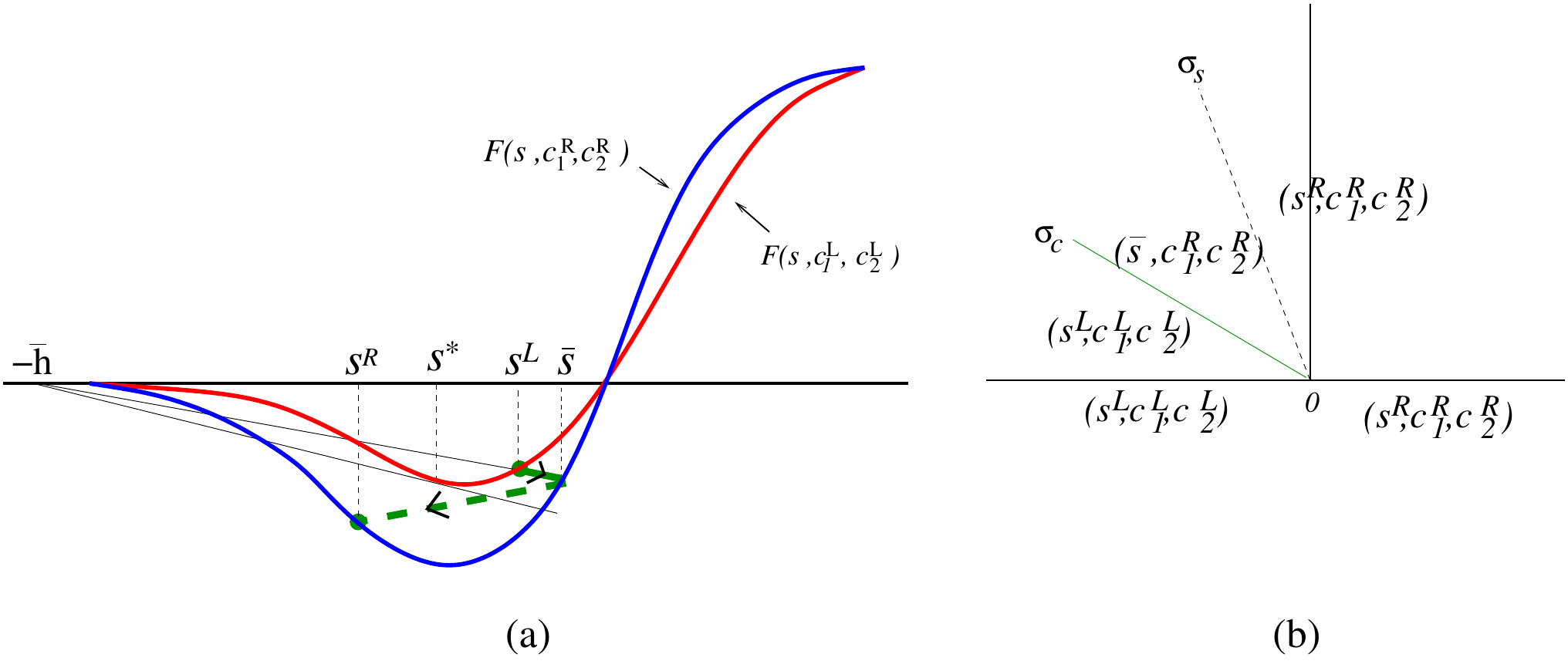}
\caption{Solution of the Riemann problem (\ref{riemannpb}) with $ s^L> s^*$ and  $ s^R<\bar{s}.$}
\label{rp-4}
\end{figure}
\end{itemize}
{\bf Remark:}
 When the flux function $F(s,c,x)$ is smooth in $s$ and $c$ and discontinuous in the $x$ variable then the construction of Riemann problem 
 is explained in the appendix of \cite{adidflu}. Here also we can construct the solution of Riemann problem in a similar way.

\subsection{Finite volume scheme }{\label{1d-4}}
We define the space grid points as $x_{\iph}=ih ,\,\, h> 0 \mbox{ and } i \in \Z $  and for
$\Delta t>0 $ define the time discretization points $t_n=n\Delta t$ for all non-negative 
integer $n,$ and $\lambda=\frac{\Delta t}{h}.$ The Finite volume scheme for the system (\ref{system1-D}) is given by 
\begin{equation}
\begin{array}{lll}
 s_i^{n+1} &= s_i^{n} - \lambda ( F^n_{\iph} -  F^n_{\imh} )\\
{c_1}_i^{n+1} s_i^{n+1}+a_1({c_1}^{n+1}_i)&= {c_1}_i^n s_i^{n}+a_1({c_1}^{n}_i)- \lambda ({G_1}^n_{\iph} -  {G_1}^n_{\imh})\\
{c_2}_i^{n+1} s_i^{n+1}+a_2({c_2}^{n+1}_i)&= {c_2}_i^n s_i^{n}+a_2({c_2}^{n}_i) - \lambda ({G_2}^n_{\iph} -  {G_2}^n_{\imh}).
\end{array}
\label{fvm1}
 \end{equation}
 where the numerical flux $F^n_{\iph},G^n_{1\iph} $ and $ G^n_{2\iph} $ are associated with the flux functions $F(s,c,x)$ 
 and ${G_l}(s,c,x)={c_l}F(s,c,x),\quad l=1,2$ and are functions of the left and right values of the saturation $s$ and the concentration $c$ at $x_{\iph}$:
 \[ F^n_{\iph} = \bar {F}(s_i^n, {c_1}_i^n,{c_2}^n_i, s_{i+1}^n, {c_1}_{i+1}^n,{c_2}^n_{i+1},x_{\iph}), \quad 
 {G_l}^n_{\iph} = \bar {G}_l(s_i^n, {c_1}_i^n,{c_2}^n_i, s_{i+1}^n, {c_1}_{i+1}^n,{c_2}^n_{i+1},x_{\iph}). \]
 The choice of the numerical flux functions $\bar F$ and $\bar {G}_l(\,l=1,2)$ determines the numerical scheme. Once we compute $s^{n+1}_i$ from the first equation of 
 (\ref{fvm1}) then we recover ${c_1}^{n+1}_i $ and ${c_2}^{n+1}_i$ from second and third  equation respectively using an iterative method, like Newton-Raphson method.

 Now we briefly explain the DFLU flux of \cite{adidflu} and Godunov flux.

 \subsection {The DFLU numerical flux}{\label{1d-5}}
 The DFLU flux is an extension of the Godunov scheme that was proposed and analyzed in
 \cite{dflu} for scalar conservations laws with a flux function discontinuous in space.
We define
 \begin{equation}
{G_l}^n_{\iph} = \begin{cases}
{c^n_l}_i F^n_{\iph} & \textrm{if } F^n_{\iph} > 0 \\
{c^n_l}_{i+1}  F^n_{\iph} & \textrm{if } F^n_{\iph} \le 0 \quad l=1,2.
\label{dflu-g}
\end{cases}
\end{equation}
Now the choice of the numerical scheme depends on the choice
of $F^n_{i+1/2}$. To do so we treat $c(x,t)$
in $F(s,c,x)$  as  a known function which may be discontinuous at the space discretization points and 
$F$ is allowed to be discontinuous in the $x$ variable at the same space  discretization points.
Therefore on each rectangle $(x_{\imh},x_{\iph}) \times (t_n,t_{n+1})$, we
 consider the conservation law:
$$ s_t+F(s,{c_1}_i^n,{c_2}_i^n,x)_x=0 $$ 
with initial condition $s(x,0)=s_i^{n}$ for $x_{\imh}< x< x_{\iph}$ (see Fig.\ref{fl}).
\begin{figure}[htbp] 
\begin{center}
\begin{picture}(350,70)(0,-5)
\thicklines
\put(50.,20){\parbox{90pt}{$$s_t + F(s, {c_1}_i^n,{c_2}_i^n,x_{i})_x =0 $$ $s(t_n) =s_i^n$}}
\put(185.,20){\parbox{100pt}{$$ s_t + F(s, {c_1}_{i+1}^n,{c_2}_{i+1}^n,x_{i+1})_x = 0$$ $s(t_n) =s_{i+1}^n$}}
\put(0.,-3.){\line(1,0){350}}
\put(0.,50.){\line(1,0){350}}
\put(175.,-2.){\line(0,1){52}}
\put(40.,-2.){\line(0,1){52}}
\put(330.,-2.){\line(0,1){52}}
\put(168.,-10){$x_{\iph}$}
\put(30.,-10){$x_{\imh}$}
\put(320.,-10){$x_{i+3/2}$}
\put(0.,5){$t=t_n$}
\put(0.,55){$t=t_{n+1}$}
\end{picture}
\end{center}
\caption{The flux functions $F(. , c_1,c_2,x)$ is discontinuous in $c_1,c_2$ and $x$ at the discretization points.}
\label{fl}
\end{figure}   
 The above problem can be considered as a conservation law with flux function discontinuous in $x$ for which DFLU flux can be used. Then the  DFLU flux is given as
 \begin{eqnarray}
\nonumber
 F^{n}_{\iph} & = & F^{\tiny{DFLU}}(s_i^n, {c_1}_i^n,{c_2}_i^n s_{i+1}^n, {c_1}_{i+1}^n,{c_2}_{i+1}^n) \\
             & = & \max\{ F(\max\{s_i^n,{\theta_i^n}\},{c_1}_i^n,{c_2}_i^n,x_{i}),
             F(\min\{s_{i+1}^n,{\theta_{i+1}^{n}}\},{c_1}_{i+1}^n,{c_2}_{i+1}^n,x_{i+1})\},\nonumber
\end{eqnarray}
where $ {\theta^n_i}=\mbox{argmin} F(.,{c_1}_i^n,{c_2}_i^n,x_{i}^{}).$         
\subsection{ The Godunov flux}
The Godunov flux at the  grid  point  $x_{\iph}$  is calculated by using the solution of  the  Riemann problem: 
\begin{eqnarray}
\nonumber
s_t + F(s,c,x)_x &=& 0 \\
\nonumber
(sc_1+a_1(c_1))_t + (c_1F(s,c,x))_x &=& 0\\
(sc_2+a_2(c_2))_t + (c_2F(s,c,x))_x &=& 0
\end{eqnarray}
in the domain 
$(x_{\imh},x_{\iph}) \times (t_n,t_{n+1}),$
with the initial condition
\be
\nonumber
(s(x,t_n),c_1(x,t_n),c_2(x,t_n))= \left\{ \begin{array}{lll}
   (s^n_i,{c_1}^n_i,{c_2}^n_i) &\mbox{if} & x<x_{\iph}   \\ (s^n_{i+1},{c_1}^n_{i+1},{c_2}^n_{i+1}) &\mbox{if} &  x>x_{\iph}.
\end{array} \right.
\ee
The numerical fluxes are given by
$$F^n_{\iph}=F(s(x_{\iph},t),c(x_{\iph},t),x_{\iph}), \quad t_n<t<t_{n+1}$$
and
\begin{equation}
\nonumber
{G_l}^n_{\iph} =c_l(x_{\iph},t)F^n_{\iph},\quad l=1,2.
\end{equation}
{\bf{Remark:}}
 In general Godunov and DFLU flux  may differ, for details see \cite{adidflu}.
\subsection{The  Upstream Mobililty flux} This flux is designed by petroleum engineers from 
physical consideration. It is an ad-hoc flux for two-phase flow in porous media which corresponds to the 
approximate solution to the Riemann problem \cite{jeffry}. To define the upstream mobility flux, assume that the absolute permeability $K(x)>0 $ and 
we redefine  the flux function in  (\ref{fluxfn})
as 
\begin{equation}
F(s,c,x) = \frac{K\lambda_w}{K\lambda_w+K\lambda_o} [v - (\rw-\ro)g K \lo(s,c)]
\label{um}
\end{equation}
Now if we take $\lambda_w=K\lambda_w$ and $\lambda_o=K\lambda_o,$ the flux function becomes 
$$F(s,c,x) = \frac{\lambda_w}{\lambda_w+\lambda_o} [v - (\rw-\ro)g  \lo(s,c)]$$
 Now the numerical fluxes are given by 
\begin{equation}
\nonumber
\begin{array}{l}
F^{n}_{i+\frac{1}{2}}(s^n_i,{c^n_1}_i,{c^n_2}_i,s^{n}_{i+1},{c^n_1}_{i+1},{c^n_2}_{i+1})= \displaystyle{ \,
         \frac{\lambda_w^*}{\lambda_w^* + \lambda_o^*}
         [ v - (\rho_w-\rho_o)g\lambda_o^* ]} ,\\
\lambda^*_\ell=
\left\{ \begin{array}{ll} \lambda_\ell(s^n_i,{c^n_1}_i,{c^n_2}_{i},K_i) & \mbox{if }
         v-(\rho_\ell-\rho_i)g\lambda_\ell >0, \;i=w,o, i\neq \ell,\\[3mm]
        \lambda_\ell(s^n_{i+1},{c^n_1}_{i+1},{c^n_2}_{i+1},K_{i+1}) & \mbox{if }
        v-(\rho_\ell-\rho_i)g\lambda_\ell \leq0, \; i=w,o, i\neq \ell
\end{array}
\right.
\end{array}
\end{equation}
and ${G^{n}_l}_{\iph} (l=1,2)$ are given as in (\ref{dflu-g}).\\
{\bf Remark:}
 The Upstream mobility flux works only for the flux function which is of the form as in (\ref{um}) where as DFLU flux 
 can be applied for any flux function which satisfies the  assumptions of \S \ref{1d-1}.

\subsection{High-order schemes}
In order to develop the second order scheme, we follow the method of lines approach in which space and time discretization are performed 
separately. In the first step, spatial discretization using piecewise linear reconstruction is made which leads to a system of ODE which can be written as
\begin{equation}
\dd{U}{t} + R(U) = 0, \quad U=\begin{bmatrix} s \\ sc_1+a_1(c_1) \\sc_2+a_2(c_2) \end{bmatrix}, \quad R(U)_{i} = \frac{1}{\dx }  \begin{bmatrix}
F_{\iph} - F_{\imh} \\
{G_1}_{\iph} - {G_1}_{\imh} \\
{G_2}_{\iph} - {G_2}_{\imh}
\end{bmatrix}
\label{eq:ode-1d}
\end{equation}
The high order accurate fluxes are given by
 $$F_{i+\half}=F(s^{L}_{i+\half},s^{R}_{i+\half},{c_1}^{L}_{i+\half},{c_2}^{L}_{i+\half},{c_1}^{R}_{i+\half},{c_2}^{R}_{i+\half})$$
\begin{equation}
{G_l}_{\iph} = \begin{cases}
{c_l}^L_{\iph} F_{\iph} & \textrm{if } F_{\iph} > 0 \\
{c_l}^R_{\iph}  F_{\iph} & \textrm{if } F_{\iph} \le 0\quad l=1,2,
\end{cases}
\end{equation}
 The quantities with superscripts $L$ and $R$ denote the reconstructed 
values of the variables to the left and right of the corresponding cell face. For any quantity $u$, we can define the reconstruction as follows:
\begin{equation}
u^L_{\iph} = u_{i} + \half \delta_{i}, \qquad
u^R_{\iph} = u_{i+1} - \half \delta _{i+1}
\end{equation}
where
\begin{equation}
\label{minimod}
\delta_{i} = \textrm{minmod}\left( \theta (u_{i} - u_{i-1}), \half (u_{i+1} - u_{i-1}), \theta (u_{i+1} - u_{i}) \right), \theta \in [1,2].
\end{equation}
Finally the time integration of the ODE (\ref{eq:ode-1d}) must be high order accurate in order for the scheme to be high order accurate. 
A third order accurate, strong stability preserving Runge-Kutta scheme due to Shu-Osher is given by
\begin{eqnarray*}
V^{(0)} &=& U^n \\
V^{(1)} &=& V^{(0)} - \dt R(V^{(0)}) \\
V^{(2)} &=& \frac{3}{4} U^n + \frac{1}{4} [V^{(1)} - \dt R(V^{(1)})] \\
V^{(3)} &=& \frac{1}{3} U^n + \frac{2}{3} [V^{(2)} - \dt R(V^{(2)})] \\
U^{n+1} &=& V^{(3)}
\end{eqnarray*}
If the explicit scheme (\ref{fvm1}) is stable in the  norm $\|.\|$, i.e., if
\begin{eqnarray}
 \dt \le \dt_c \Longrightarrow \norm{U - \dt R(U)} \le \norm{U},\mbox{ where $\dt_c$ is the CFL restricted time step,}
\end{eqnarray}
then the above Runge-Kutta scheme is also stable in the same norm under the same time-step
restriction \\(cf.\cite{rk1,rk}).
\subsection{Maximum principle on saturation}\label{stability}
Let us write $$\bar{s}^n=(\bar{s}_i^n,)^4_{i=1}=(s^{nL}_{\imh},s^{nR}_{\imh}
,s^{nL}_{\iph},s^{nR}_{\iph})$$
$$\bar{c}^n=(\bar{c}_i^n)^{8}_{i=1}=({c_1}^{nL}_{\imh},{c_1}^{nR}_{\imh},{c_1}^{nL}_{\iph},{c_1}^{nR}_{\iph},
{c_2}^{nL}_{\imh},{c_2}^{nR}_{\imh},{c_2}^{nL}_{\iph},{c_2}^{nR}_{\iph})$$
The updated value of the saturation (\ref{fvm1}) can be written as
$$s^{n+1}_{i}=H(\bar{s}^n,\bar{c}^n).$$
Where $H$ is Lipschitz continuous in saturation and concentration.
Since the slope limiter preserves the average value of the solution in each cell, we can express this as
\begin{equation}
\label{finalscheme}
s^{n+1}_{i}= \dfrac{s^{nL}_{\iph}+s^{nR}_{\imh}}{2}-\lambda (F^n_{i+\half}-F^n_{i-\half} ).
\end{equation}
If we differentiate $H$ with respect to  its variables $\bar{s}^n_i$ we  can observe that $\df{}{\bar{s}^n_i}H\geq0$ provided 
\begin{equation}
\lambda|\df{}{\bar{s}^n_i}F^n_{i\pm\half}|\leq \frac{1}{2}.
\label{cfl-1-d}
\end{equation}
Let
$$ M=\sup_{s} \{ \df{F_1}{s},\df{F_2}{s}, \frac{F_1}{s+h_l}, \frac{F_2}{s+h_l} \} , $$
then the  condition (\ref{cfl-1-d}) reduces to,  
\begin{equation}
\lambda M  \leq\frac{1}{2}.
 \label{cfl-1d}
\end{equation}
This shows that $H$ is monotone in each of its variable. Using these facts we have the following lemmas. 
\begin{lemma} 
 Let $s_0 \in [0,1]$ be the initial data and let $\{s^n_{i}\}$ be the corresponding solution calculated by the finite
volume scheme (\ref{fvm1}) using DFLU flux along with slope limiter. If the CFL given in (\ref{cfl-1d}) holds  then 
\begin{equation} 
\label{a}
 0\leq s^n_{i} \leq 1 \,\,\,\forall,\, i \mbox{ and } n.
\end{equation}
\end{lemma}
{\bf Proof:}
 From the property of slope limiter we can observe that whenever $ 0\leq s^n_{i} \leq 1 $ then 
the reconstructed values satisfies 
$$ 0\leq s^{nL}_{\ipmh},s^{nR}_{\ipmh}\leq1  \,\,\forall \,\,i \mbox{ and }n$$
Using this property and the monotonicity of the  $H$ , we get 
$$
\begin{array}{lll}
0=H({\bf{0}},\bar{c}^n )\leq H(\bar{s}^n,\bar{c}^n)=s^{n+1}_{i} \leq H({\bf{1}},\bar{c}^n)=1&&
\end{array}
$$
This proves that $$0\leq s^{n+1}_{i}\leq 1 \,\, \forall\,\, i,n.$$
\hfill \qed
\subsection{Maximum principle and TVD for concentration}
\begin{theorem}\label{tvd}
Let $\{{c_1}^n_i\}$ , $\{{c_2}^n_i\}$ be the solution calculated by the finite
volume scheme (\ref{fvm1})  using DFLU flux  with slope limiter. Under the CFL condition $\lambda M \leq \frac{1}{2},$ concentration $c=(c_1,c_2)$ satisfies
 \begin{enumerate}
  \item [(a)] $\min\{{c_l}^n_{i-1},{c_l}^n_{i},{c_l}^n_{i+1}\}\leq {c_l}^{n+1}_i \leq \max\{{c_l}^n_{i-1},{c_l}^n_{i},{c_l}^n_{i+1}\}\,\,\,\forall
  \,\, n\in\Z^+,\,\,i \in \Z \quad l=1,2.$
\item[(b)]$\displaystyle \sum_{i}^{} |{c_l}^{n+1}_i-{c_l}^{n+1}_{i-1}|\leq \sum_{i}^{} |{c_l}^{n}_i-{c_l}^{n}_{i-1}|  \,\,\,\forall \,\, n\in\Z^+, \quad l=1,2.$
 \end{enumerate}
\end{theorem}
  {\bf Proof:}
  From the finite volume scheme (\ref{fvm1}),
  \begin{equation*}
\begin{array}{lll}
 s_i^{n+1} &=& s_i^{n} - \lambda ( F^n_{\iph} -  F^n_{\imh} )\\
{c_1}_i^{n+1} s_i^{n+1}+a_1({c_1}_i^{n+1})&=& {c_1}_i^n s_i^{n}+a_1({c_1}_i^{n}) -\lambda ({G_1}^n_{\iph} -  {G_1}^n_{\imh})\\
{c_2}_i^{n+1} s_i^{n+1}+a_2({c_2}_i^{n+1})&=& {c_2}_i^n s_i^{n}+a_2({c_2}_i^{n}) -\lambda ({G_2}^n_{\iph} -  {G_2}^n_{\imh}).
\end{array}
 \end{equation*}
 We can express the numerical flux ${G_1}_{\iph},{G_2}_{\iph}$  as (here we suppress the index $n$ for fluxes )
 $${G_1}_{\iph}={c_1}^L_{\iph}F^+_{\iph}+{c_1}^R_{\iph}F^-_{\iph}$$
 $${G_2}_{\iph}={c_2}^L_{\iph}F^+_{\iph}+{c_2}^R_{\iph}F^-_{\iph},$$
 where
$$F^+_{\iph}=\max\{F_{\iph},0\} , \, \,\, F^-_{\iph}=\min\{F_{\iph},0\}$$
We write the scheme (\ref{fvm1}) as 
$$s^{n+1}_i{c_l}^{n+1}_i+a_l({c_l}_i^{n+1})-s^n_i{c_l}^n_i-a_l({c_l}_i^{n})+\lambda ({c_l}^{nL}_{\iph}F^+_{\iph}+{c_l}^{nR}_{\iph}F^-
_{\iph}-({c_l}^{nL}_{\imh}F^+_{\imh}+{c_l}^{nR}_{\imh}F^-_{\imh}))=0$$
By adding and subtracting the term $s^{n+1}_i{c_l}^n_i,$ we get
\begin{align*}
& (s^{n+1}_i+a'_l(\zeta^{n+\half}_i))({c_l}^{n+1}_i-{c_l}^n_i)+{c_l}^n_i(s^{n+1}_i-s^n_i)\notag \\
&\hspace{2.4 cm}+\lambda ({c_l}^{nL}_{\iph}F^+_{\iph}+{c_l}^{nR}_{\iph}F^-_{\iph}-({c_l}^{nL}_{\imh}F^+_{\imh}+{c_l}^{nR}_{\imh}F^-_{\imh}))=0.
\end{align*}
where $a_l({c_l}^{n+1}_i)-a_l({c_l}_i^{n})=a'_l(\zeta^{n+\half}_i)({c_l}^{n+1}_i-{c_l}^{n}_i),$ for some $\zeta^{n+\half}_i $ between 
${c_l}^{n+1}_i$ and ${c_l}^{n}_i.$
By replacing $s^{n+1}_i-s^n_i$ by $-\lambda(F_{\iph}-F_{\imh})$   and splitting   $F_{\ipmh}$  by $(F^+_{\ipmh}+F^-_{\ipmh})$  we have 
\begin{align*}
&(s^{n+1}_i+a'_l(\zeta^{n+\half}_i))({c_l}^{n+1}_i-{c_l}^n_i)-\lambda {c_l}^n_i(F^+_{\iph}+F^-_{\iph}-F^+_{\imh}-F^-_{\imh})\\
&\hspace{2.3 cm}+\lambda ({c_l}^{nL}_{\iph}F^+_{\iph}+{c_l}^{nR}_{\iph}F^-_{\iph}-({c_l}^{nL}_{\imh}F^+_{\imh}+{c_l}^{nR}_{\imh}F^-_{\imh}))=0.
\end{align*}
By rearranging the terms in the above equation we get 
\begin{align*}
&(s^{n+1}_i+a'_l(\zeta^{n+\half}_i))({c_l}^{n+1}_i-{c_l}^n_i)+\lambda F^+_{\iph} ( {c_l}^{nL}_{\iph}-{c_l}^n_i )\\
&\hspace{1.4 cm}+\lambda F^-_{\iph}({c_l}^{nR}_{\iph}-{c_l}^n_i)+\lambda F^+_{\imh}({c_l}^n_i-{c_l}^{nL}_{\imh})+\lambda F^-_{\imh}({c_l}^n_i
-{c_l}^{nR}_{\imh})=0.
\end{align*}
Note that  $${c_l}^{nL}_{\iph}={c_l}^n_i+\frac{\delta_i}{2}, \,\, \,\,{c_l}^{nR}_{\imh}={c_l}^n_i-\frac{\delta_i}{2},$$ 
and $\delta_i $ is the slope limiter given by
$$\delta_{i} = \textrm{minmod}\left( \theta({c_l}^n_{i} - {c_l}^n_{i-1}), \half ({c_l}^n_{i+1} - {c_l}^n_{i-1}), \theta({c_l}^n_{i+1} - {c_l}^n_{i}) \right).$$
After substituting the values for ${c_l}^{nL}_{\ipmh}$ and ${c_l}^{nR}_{\ipmh}$ the above equation becomes

\begin{align*}
&(s^{n+1}_i+a'_l(\zeta^{n+\half}_i))({c_l}^{n+1}_i-{c_l}^n_i)+\lambda F^+_{\iph}\frac{\delta_i}{2} +\lambda
F^-_{\iph}(1-\frac{\delta_{i+1}}{2({c_l}^n_{i+1}-{c_l}^n_{i})})({c_l}^n_{i+1}-{c_l}^n_i)\\
&\hspace{3 cm} +\lambda F^+_{\imh}(1-\frac{\delta_{i-1}}{2({c_l}^n_{i}-{c_l}^n_{i-1})})({c_l}^n_{i}-{c_l}^n_{i-1})+\lambda F^-_{\imh}\frac{\delta_i}{2}=0.
\end{align*}
Now we can write 

\begin{align}
&{c_l}^{n+1}_i={c_l}^n_i-\lambda\frac{ F^+_{\iph}}{2(s^{n+1}_i+a'_l(\zeta^{n+\half}_i)) ({c_l}^n_i-{c_l}^n_{i-1})}\delta_i({c_l}^n_i-{c_l}^n_{i-1})\\
&\hspace{1.7 cm}-\lambda\frac{F^-_{\iph}}{(s^{n+1}_i+a'_l(\zeta^{n+\half}_i))}(1-\frac{\delta_{i+1}}{2({c_l}^n_{i+1}-{c_l}^n_{i})})({c_l}^n_{i+1}-{c_l}^n_i)\notag\\
&\hspace{1.8 cm}-\lambda\frac{F^+_{\imh}}{(s^{n+1}_i+a'_l(\zeta^{n+\half}_i))}(1-\frac{\delta_{i-1}}{2({c_l}^n_{i}-{c_l}^n_{i-1})})({c_l}^n_{i}-{c_l}^n_{i-1})\\
&\hspace{1.7 cm}-\lambda \frac{F^-_{\imh}}{(s^{n+1}_i+a'_l(\zeta^{n+\half}_i))}\frac{\delta_i}{2({c_l}^n_{i+1}-{c_l}^n_i)}({c_l}^n_{i+1}-{c_l}^n_i)\notag\\
&\hspace{0.8 cm}={c_l}^n_i-\alpha^1_{\imh}({c_l}^n_i-{c_l}^n_{i-1})+\alpha^2_{\iph}({c_l}^n_{i+1}-{c_l}^n_i)
-\alpha^3_{\imh}({c_l}^n_i-{c_l}^n_{i-1})+\alpha^4_{\iph}({c_l}^n_{i+1}-{c_l}^n_i)\notag\\
&\hspace{0.8 cm}={c_l}^n_i- (\alpha^1_{\imh}+\alpha^3_{\imh})({c_l}^n_i-{c_l}^n_{i-1})+(\alpha^2_{\iph}+\alpha^4_{\iph})({c_l}^n_{i+1}-{c_l}^n_i)\notag\\
&\hspace{0.8 cm}={c_l}^n_i- C^n_{\imh}({c_l}^n_i-{c_l}^n_{i-1})+D^n_{\iph}({c_l}^n_{i+1}-{c_l}^n_i),\label{max-1d}
\end{align}
where
\begin{align*}
 C^n_{\imh}=\alpha^1_{\imh}+\alpha^3_{\imh}, \quad D^n_{\iph}=\alpha^2_{\iph}+\alpha^4_{\iph}
\end{align*}
and
$$\alpha^1_{\imh}=\lambda\frac{ F^+_{\iph}}{(s^{n+1}_i+a'_l(\zeta^{n+\half}_i)) 2({c_l}^n_i-{c_l}^n_{i-1})}\delta_i,\,\,\alpha^2_{\iph}=-\lambda\frac{
F^-_{\iph}}{(s^{n+1}_i+a'_l(\zeta^{n+\half}_i))}(1-\frac{\delta_{i+1}}{2({c_l}^n_{i+1}-{c_l}^n_{i})})$$
$$ \alpha^3_{\imh}=\lambda\frac{F^+_{\imh}}{(s^{n+1}_i+a'_l(\zeta^{n+\half}_i))}(1-\frac{\delta_{i-1}}{2({c_l}^n_{i}-{c_l}^n_{i-1})}),\,\,\,
\alpha^4_{\iph}=-\lambda \frac{F^-_{\imh}}{(s^{n+1}_i+a'_l(\zeta^{n+\half}_i))}
\frac{\delta_i}{2({c_l}^n_{i+1}-{c_l}^n_i)}.
$$
From the property of the limiter it is easy to see that
\begin{align}
 0\leq \frac{\delta_{i+1}}{2({c_l}^n_{i+1}-{c_l}^n_{i})}\leq 1
 \end{align}
which in turn implies
\begin{align}
  C^n_{\imh}, D^n_{\iph}\geq 0 \quad  \forall\,\, i.
 \label{positive-1d}
\end{align}
Now we prove the maximum principle for $c_l (l=1,2) $ by considering the following cases.\\
\par 
{\bf Case1:} Suppose that ${c_l}^n_i$ lies between $ {c_l}^n_{i-1}$ and $ {c_l}^n_{i+1}$ then 
\begin{align*}
 {c_l}^n_i=\theta {c_l}^n_{i-1}+(1-\theta){c_l}^n_{i+1} \quad \mbox{ for some } \,\,\theta \in [0,1]
\end{align*}
and
\begin{align*}
 {c_l}^n_i-{c_l}^n_{i-1}&=(1-\theta)({c_l}^n_{i+1}-{c_l}^n_{i-1})\\
 {c_l}^n_{i+1}-{c_l}^n_{i}&=\theta({c_l}^n_{i+1}-{c_l}^n_{i-1}).
\end{align*}
Now from (\ref{max-1d}) we write 
\begin{align}
 {c_l}^{n+1}_i&=(1-\theta)({c_l}^n_{i+1}-{c_l}^n_{i-1})-C^n_{\imh}(1-\theta)({c_l}^n_{i+1}-{c_l}^n_{i-1})\notag\\
 &\hspace{5 cm} + D^n_{\iph} \theta ({c_l}^n_{i+1}-{c_l}^n_{i-1})\notag\\
 &=\lambda_1 {c_l}^n_{i-1}+\lambda_2 {c_l}^n_{i+1}\label{max1-1d},
\end{align}
where 
\begin{align*}
 \lambda_1=\theta(1-D^n_{\iph})+C^n_{\imh}(1-\theta)\\
 \lambda_2=(1-\theta)(1-C^n_{\imh})+\theta D^n_{\iph}.
\end{align*}
 Note that $\lambda_1+\lambda_2=1,$ under the CFL condition $\lambda M\leq \frac{1}{2}$
  we have $C^n_{\imh},D^n_{\iph}\leq 1$ which gives  $\lambda_l, \lambda_2\geq0.$
 Hence from (\ref{max1-1d}) the maximum principle (a) follows.\\
 \par 
{\bf Case2:}
Suppose ${c_l}^n_i $ does not lies between ${c_l}^n_{i-1}$ and ${c_l}^n_{i+1},$ then we have $\delta_i=0.$
i.e., $$C^n_{\imh}=\alpha^3_{\imh} \mbox { and } D^n_{\iph}=\alpha^2_{\iph}$$
The equation (\ref{max-1d}) can be rewritten as 
\begin{align*}
 {c_l}^{n+1}_i=(1-C^n_{\imh}-D^n_{\iph}) {c_l}^n_i+C^n_{\imh}{c_l}^n_{i-1}+D^n_{\iph}{c_l}^n_{i+1}.
\end{align*}
Note that 
$$C^n_{\imh}+D^n_{\iph} \leq 1 \mbox{ under the CFL condition } \lambda M \leq \frac{1}{2}.$$
This proves  the maximum principle(a).\\
\par

To prove the TVD property, consider
\begin{align}
 C^n_{\iph}+D^n_{\iph} &=\lambda\frac{ F^+_{\ipt}}{(s^{n+1}_{i}+a'_l(\zeta^{n+\half}_i)) 2({c_l}^n_{i+1}-{c_l}^n_{i})}\delta_{i+1} \\
&\hspace{0.5 cm}+\lambda\frac{F^+_{\iph}}{(s^{n+1}_{i}+a'_l(\zeta^{n+\half}_i))}(1-\frac{\delta_{i}}{2({c_l}^n_{i+1}-{c_l}^n_{i})})\notag\\
&\hspace{0.5 cm}-\lambda\frac{F^-_{\iph}}{(s^{n+1}_i+a'_l(\zeta^{n+\half}_i))}(1-\frac{\delta_{i+1}}{2( {c_l}^n_{i+1}-{c_l}^n_{i})})\notag\\
&\hspace{0.5 cm}-\lambda \frac{F^-_{\imh}}{(s^{n+1}_i+a'_l(\zeta^{n+\half}_i))}\frac{\delta_i}{2({c_l}^n_{i+1}-{c_l}^n_i)}\notag\\
&\leq  \lambda M \left(\frac{\delta_{i+1}}{ 2({c_l}^n_{i+1}-{c_l}^n_{i}) }+1-\frac{\delta_{i}}{2({c_l}^n_{i+1}-{c_l}^n_{i})}\right.\notag\\
&\hspace{2 cm} \left.+1-\frac{\delta_{i+1}}{ 2({c_l}^n_{i+1}-{c_l}^n_{i}) }+\frac{\delta_i}{2({c_l}^n_{i+1}-{c_l}^n_i)}\right)\notag\\
& = 2 \lambda M\leq 1,\label{bound-1d}
\end{align}
under the CFL condition  $\lambda M \leq \frac{1}{2}.$ From (\ref{positive-1d}) and (\ref{bound-1d})  the TVD property(b) follows from the Harten's lemma.
\hfill  \qed 

{\bf {Remark:}}
 Note that saturation  $s$ need not  be of total variation bounded because of $f=f(s,c)$ and $c=c(x,t)$ is discontinuous (see \cite{cpam}). The singular
mapping technique as in \cite{dflu} to prove the convergence of $s^n_{i,j}$ looks very difficult to apply. However by using the method of compensated compactness,
Kalrsen, Mishra, Risebro \cite{karlsen} showed the convergence of approximated solution in the case of a triangular system. By using their results 
in the case of  a single component polymer $(m=1)$ under suitable assumptions, in \cite{adidflu}  convergence analysis of the  saturation is studied.
\subsection{ Numerical results}{\label{1d-5}}
\par 
   Here we have chosen the flux function for the above system of equations (\ref{system1-D}) with 
  $v=0.2$ , $K\equiv1$ , $\lambda_w=\frac{s^2}{0.5+c_1+c_2}$ , $\lambda_o={(1-s)^2}$, $\rho_w g=2$ and $\rho_o g=1.$ 
The adsorption term is given by $a_l(c_l)=1+0.5c_l\quad(l=1,2).$
  In the numerical experiment the initial data is chosen so that the flux function $F$
is allowed to change the sign, equivalently eigenvalues $\lambda^l (l=1,2)$ of the system (\ref{system1-D})  allowed to change the sign.
For this purpose the initial data is chosen as  
\be
\nonumber
(s(x,0),c_1(x,0),c_2(x,0))= \left\{ \begin{array}{lll}
   (0.1,1,0.6) &\mbox{if} & x<0.4   \\ (1.0,0,0) &\mbox{if} &  x>0.4.
\end{array} \right.
\ee

\begin{figure}
\begin{minipage}[b]{0.45\linewidth}
\includegraphics[scale=0.32]{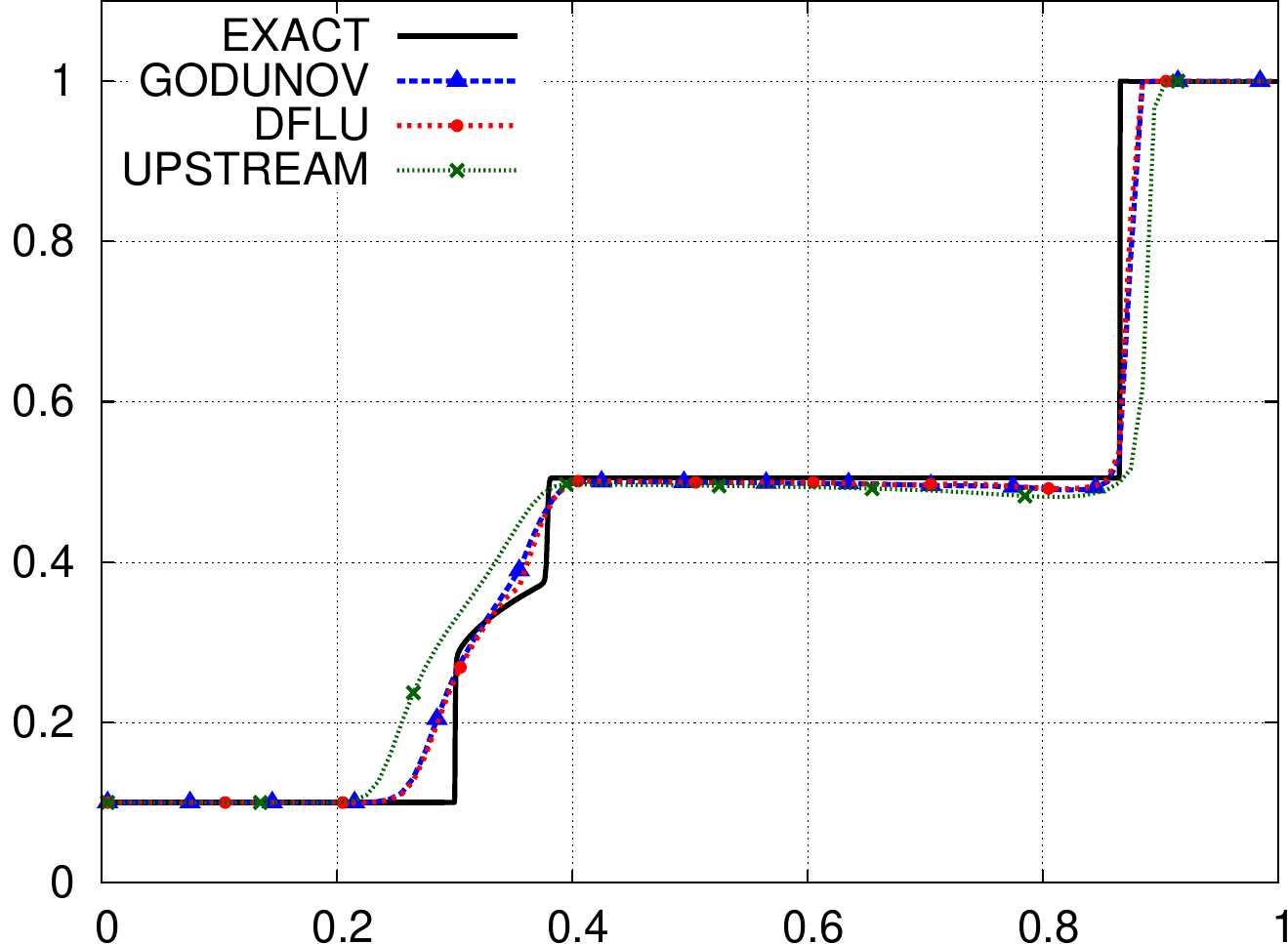}
\end{minipage}
\hspace{0.5cm}
\begin{minipage}[b]{0.45\linewidth}
\includegraphics[scale=0.32]{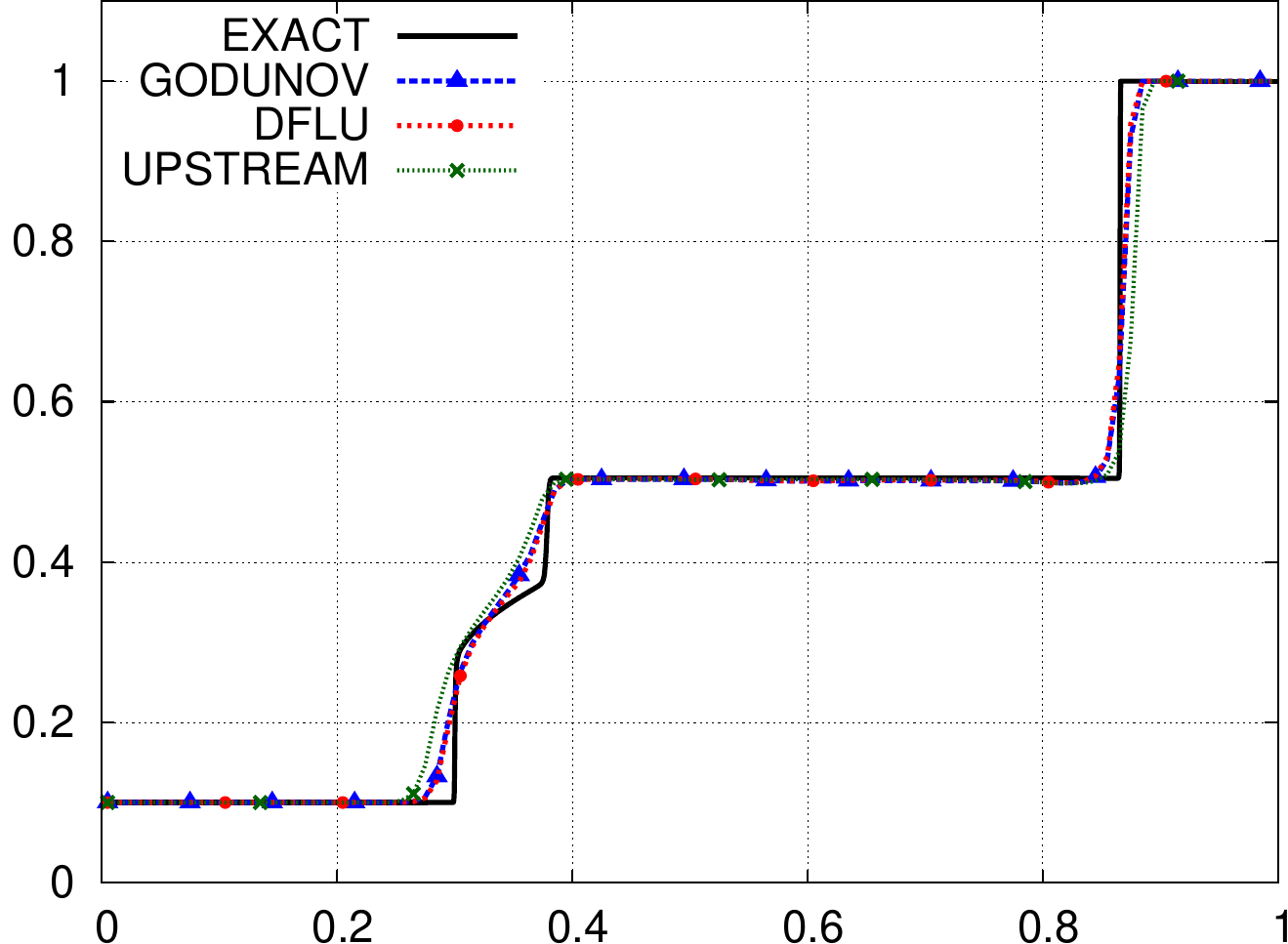}
\end{minipage}
\caption{Saturation $s$ for first order scheme (left), high order scheme (right) at time $t=1$, mesh size $h=\frac{1}{100}.$}  
\label{expt1}

\end{figure}
\begin{figure}
\begin{minipage}[b]{0.45\linewidth}
\includegraphics[scale=0.32]{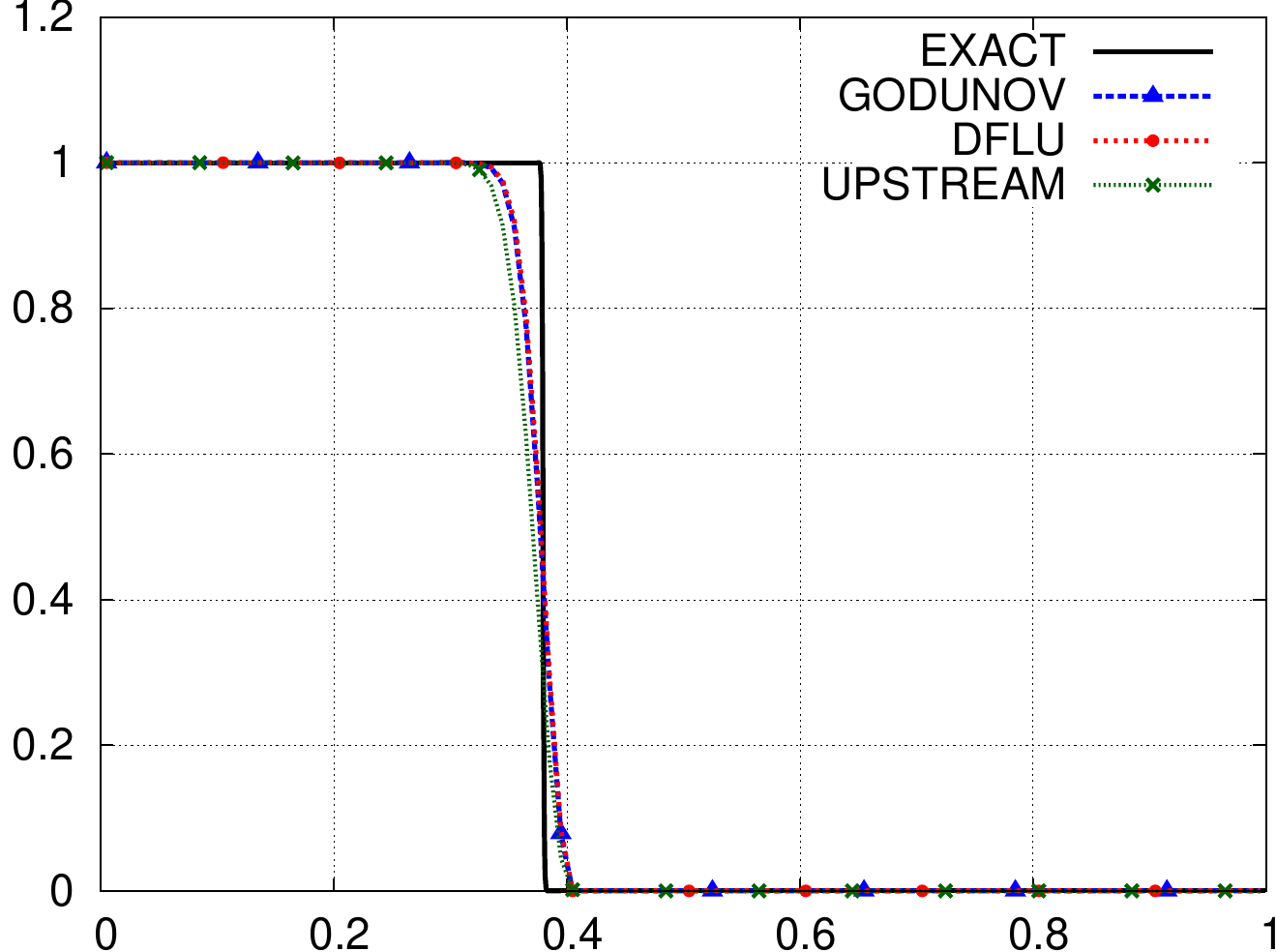}
\end{minipage}
\hspace{0.5cm}
\begin{minipage}[b]{0.45\linewidth}
\includegraphics[scale=0.32]{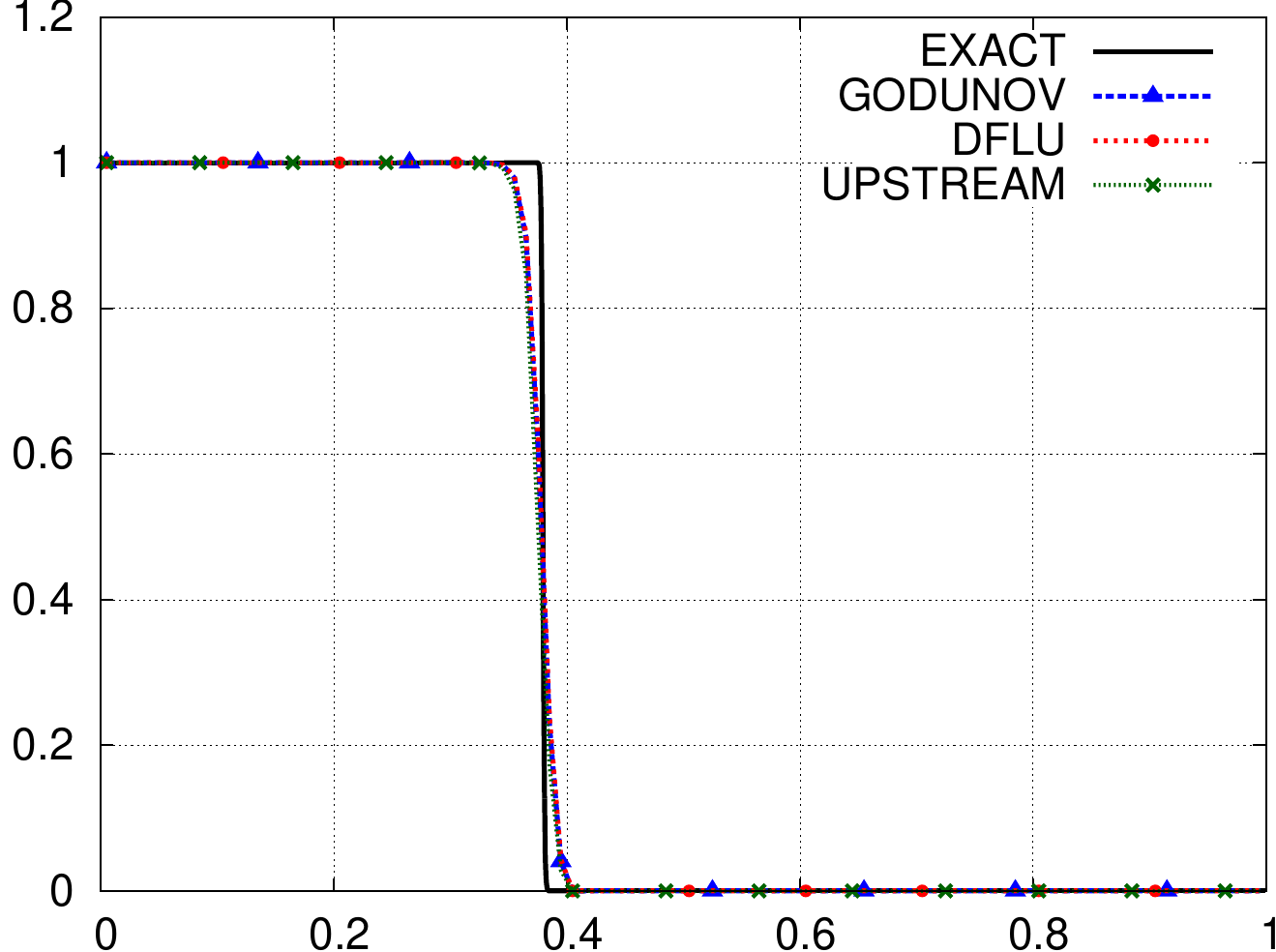}
\end{minipage}
\caption{Concentration $c_1$ for first order scheme (left), high order scheme (right) at time $t=1 ($right), mesh size $h=\frac{1}{100}.$}  
\label{expt2}
\end{figure}
\begin{figure}
\begin{minipage}[b]{0.45\linewidth}
\centering
\includegraphics[scale=0.32]{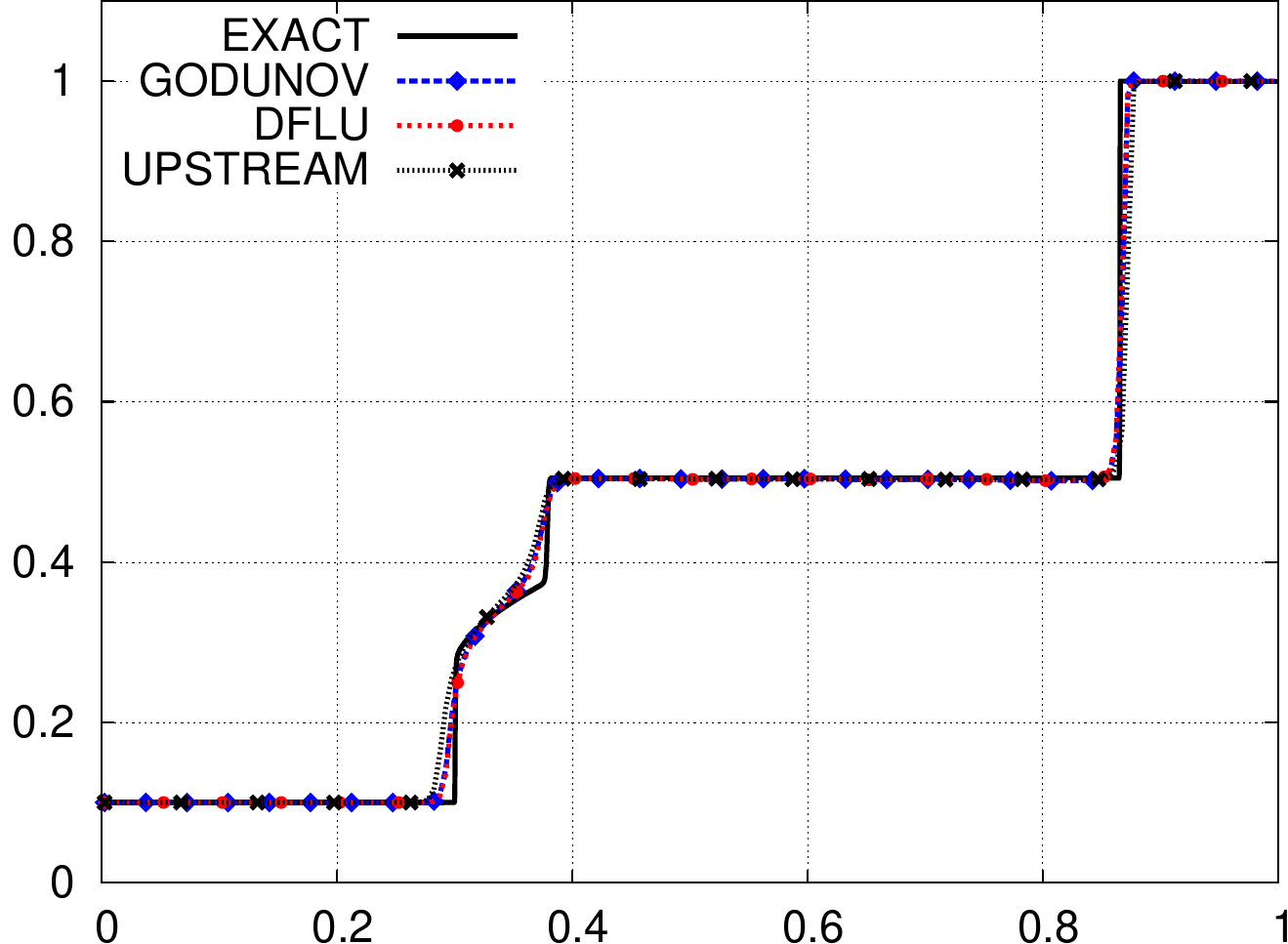}
\end{minipage}
\hspace{0.5cm}
\begin{minipage}[b]{0.45\linewidth}
\centering
\includegraphics[scale=0.32]{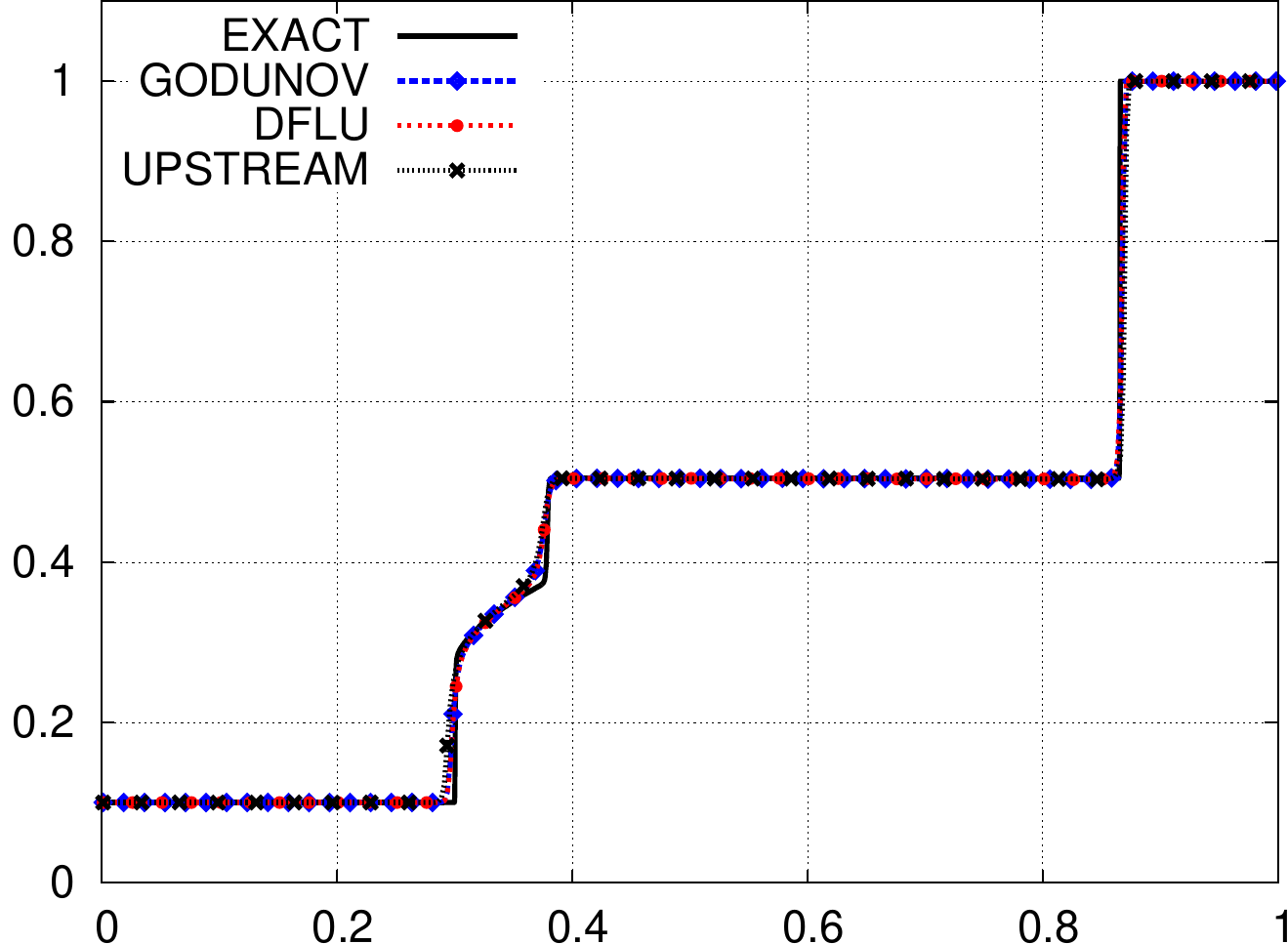}
\end{minipage}
\caption{Saturation $s$  at time $t=1$ with mesh size  $h=\frac{1}{200}$ (left) and $h=\frac{1}{400}$ (right), with high order accuracy.}
\label{expt3}
\end{figure}

\begin{figure}
\begin{minipage}[b]{0.45\linewidth}
\centering
\includegraphics[scale=0.32]{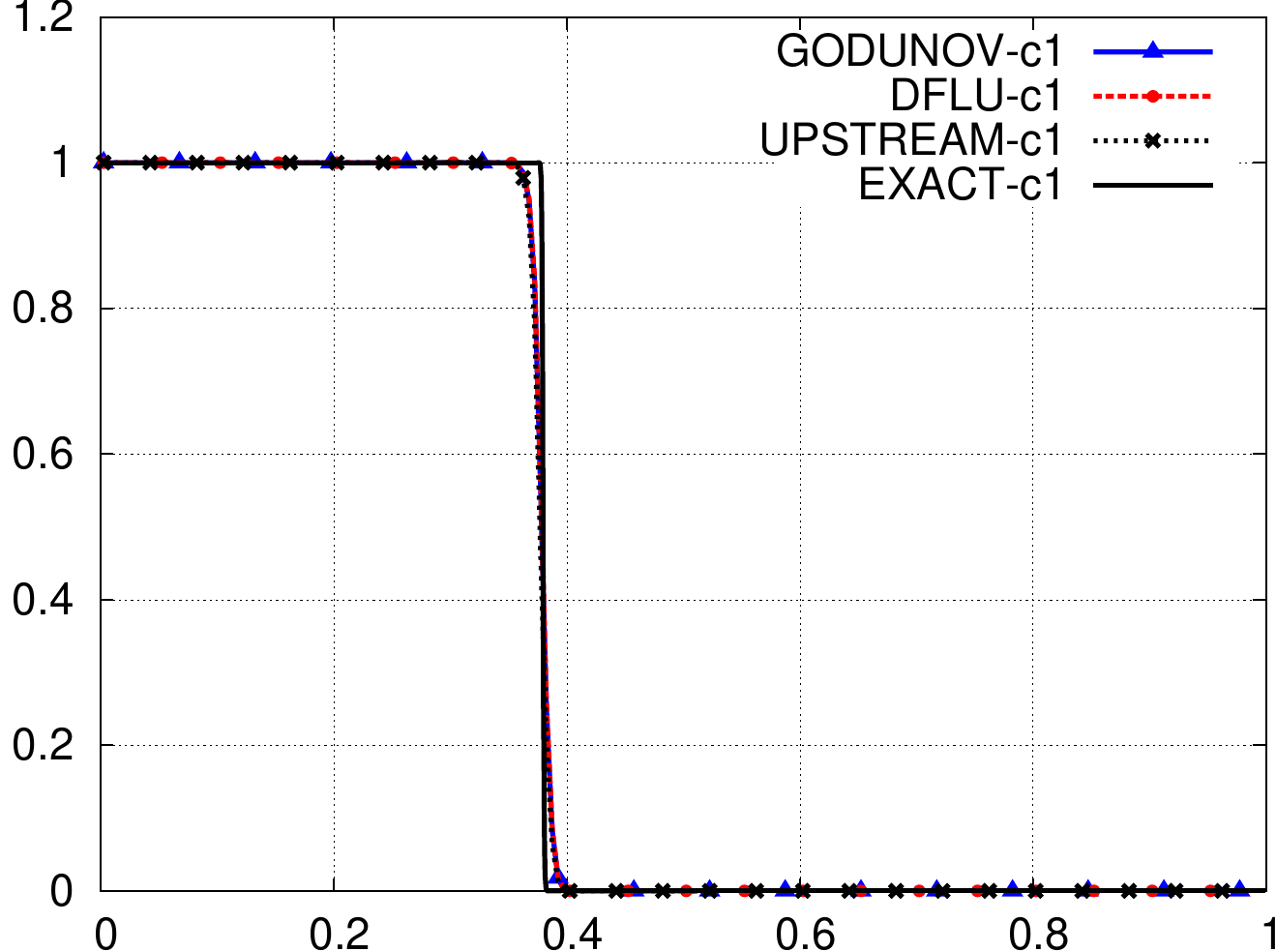}
\end{minipage}
\hspace{0.5cm}
\begin{minipage}[b]{0.45\linewidth}
\centering
\includegraphics[scale=0.32]{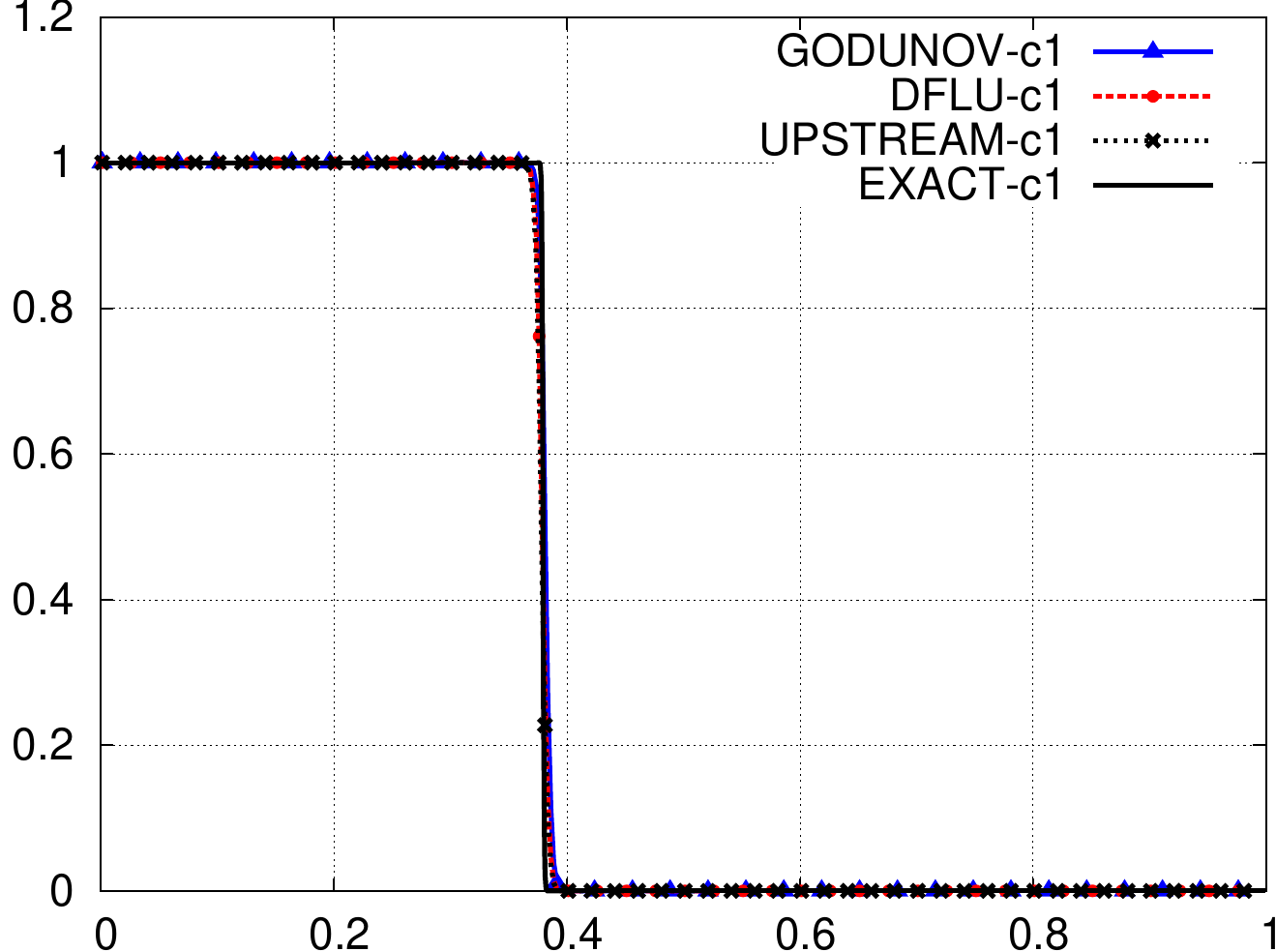}
\end{minipage}
\caption{Concentration $c_1$  at time $t=1$ with mesh size  $h=\frac{1}{200}$ (left) and $h=\frac{1}{400}$ (right), with high order accuracy.}
\label{expt4}
\end{figure}

\begin{table}
\centering
    \begin{tabular}{ | l || l | l || l | l ||l|l| }
        \hline
           & DFLU                 &  & GODUNOV &  & UPSTREAM & \\ 
           h   & $\|s-s^h\|_{L^1}$ & $\alpha$  &  $\|s-s^h\|_{L^1}$ &  $\alpha$ & $\|s-s^h\|_{L^1}$ &  $\alpha$\\ \hline
        \hline 
        1/50  & 4.2336$\times 10^{-2}$ & ~ & 4.8839$\times 10^{-2}$ & ~        &    6.3189$\times 10^{-2}$   & ~  \\
        \hline
        1/100 & 2.4366$\times 10^{-2}$ & 0.7970& 2.7735$\times 10^{-2}$ & 0.8163  & 3.6055$\times 10^{-2}$  &0.8095    \\
        \hline
        1/200 & 1.3605$\times 10^{-2}$ & 0.8407& 1.5268$\times 10^{-2}$& 0.8612& 1.9805$\times 10^{-2}$   &  0.8643 \\ 
        \hline
        1/400 & 6.2334$\times 10^{-3}$ & 1.1260 & 6.9589$\times 10^{-3}$& 1.133&9.2108$\times 10^{-3}$  & 1.1045  \\
        \hline
        1/800 & 2.2233$\times 10^{-3}$ & 1.4873 & 2.4398$\times 10^{-3}$& 1.5121&3.3674$\times 10^{-3}$  & 1.4517  \\

        \hline
    \end{tabular}
    
    \vspace{0.5 cm}
\centering
    \begin{tabular}{|l||l|l||l|l||l|l|}
        \hline
            & DFLU& & GODUNOV& &UPSTREAM &  \\ 
        h& $||c_1-{c^h_1}||_{L^1}$ & $\alpha$ & $\|c_1-{c^h_1}||_{L^1}$ &  $\alpha$ &$\|c_1-{c^h_1}||_{L^1}$ &  $\alpha$ \\ \hline
        \hline
        1/50  & 3.3257$\times 10^{-2}$& ~        & 3.9971$\times 10^{-2}$  & ~         & 5.0529$\times 10^{-2}$  &~\\ 
        \hline
        1/100 & 2.2303$\times 10^{-2}$& 0.5764  & 2.5938$\times 10^{-2}$& 0.6239& 3.4946$\times 10^{-2}$ & 0.5319   \\ 
        \hline
        1/200 & 1.2304$\times 10^{-2}$& 0.8582   & 1.4014$\times 10^{-2}$& 0.8881& 1.874$\times 10^{-2}$ &  0.899  \\ 
        \hline
        1/400 & 4.8878$\times 10^{-3}$& 1.3318 & 5.4714$\times 10^{-3}$& 1.3569&7.9071$\times 10^{-3}$  &  1.2449   \\
       \hline
        1/800 & 1.6586$\times 10^{-3}$& 1.5592 & 1.8413$\times 10^{-3}$& 1.5712&2.8197$\times 10^{-3}$  &  1.4876   \\

        \hline
    \end{tabular}
    
    \vspace{0.5 cm}
\centering
    \begin{tabular}{|l||l|l||l|l||l|l|}
        \hline
           & DFLU & & GODUNOV & & UPSTREAM & \\
        h&  $||c_2-{c^h_2}||_{L^1}$ & $\alpha$ & $\|c_2-{c^h_2}||_{L^1}$ &  $\alpha$ &  $\|c_2-{c^h_2}||_{L^1}$ &  $\alpha$\\ \hline
        \hline
        1/50  & 1.9954$\times 10^{-2}$   & ~        & 2.3983$\times 10^{-2}$   & ~   &3.0318$\times 10^{-2} $ & ~   \\ 
        \hline
        1/100 & 1.3382$\times 10^{-2}$ & 0.5764   & 1.5563$\times 10^{-2}$    &0.6239& 2.0968$\times 10^{-2} $ &0.5319  \\ 
        \hline
        1/200 & 7.3821$\times 10^{-3}$ & 0.8581   & 8.4086$\times 10^{-3}$  & 0.8881 &1.1244$\times 10^{-2} $  &0.899 \\ 
        \hline
        1/400 & 2.9327$\times 10^{-3}$ & 1.3318 & 3.2829$\times 10^{-3}$ & 1.3569 & 4.7455$\times 10^{-3} $  &1.2445\\
        \hline
        1/800 & 9.9518$\times 10^{-4}$ & 1.5592& 1.10479$\times 10^{-3}$ & 1.5712 & 1.6924$\times 10^{-3} $  &1.4875\\
         \hline
        \end{tabular}
        \caption{$L^1$ error for saturation $s$ and concentrations $c_1$ and $c_2.$}
\end{table}
\begin{figure}
\begin{minipage}[b]{0.45\linewidth}
\centering
\includegraphics[scale=0.32]{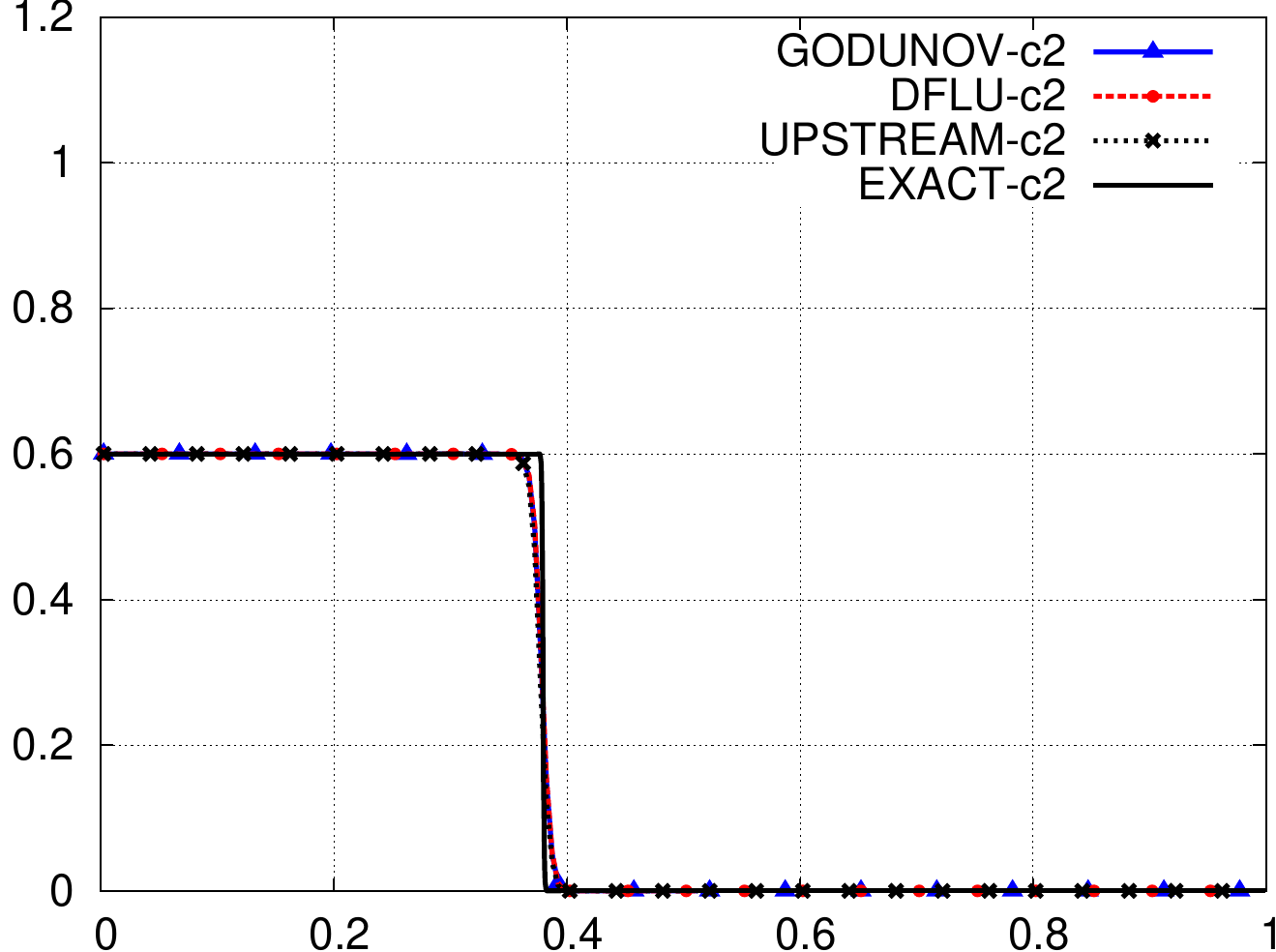}
\end{minipage}
\hspace{0.5cm}
\begin{minipage}[b]{0.45\linewidth}
\centering
\includegraphics[scale=0.32]{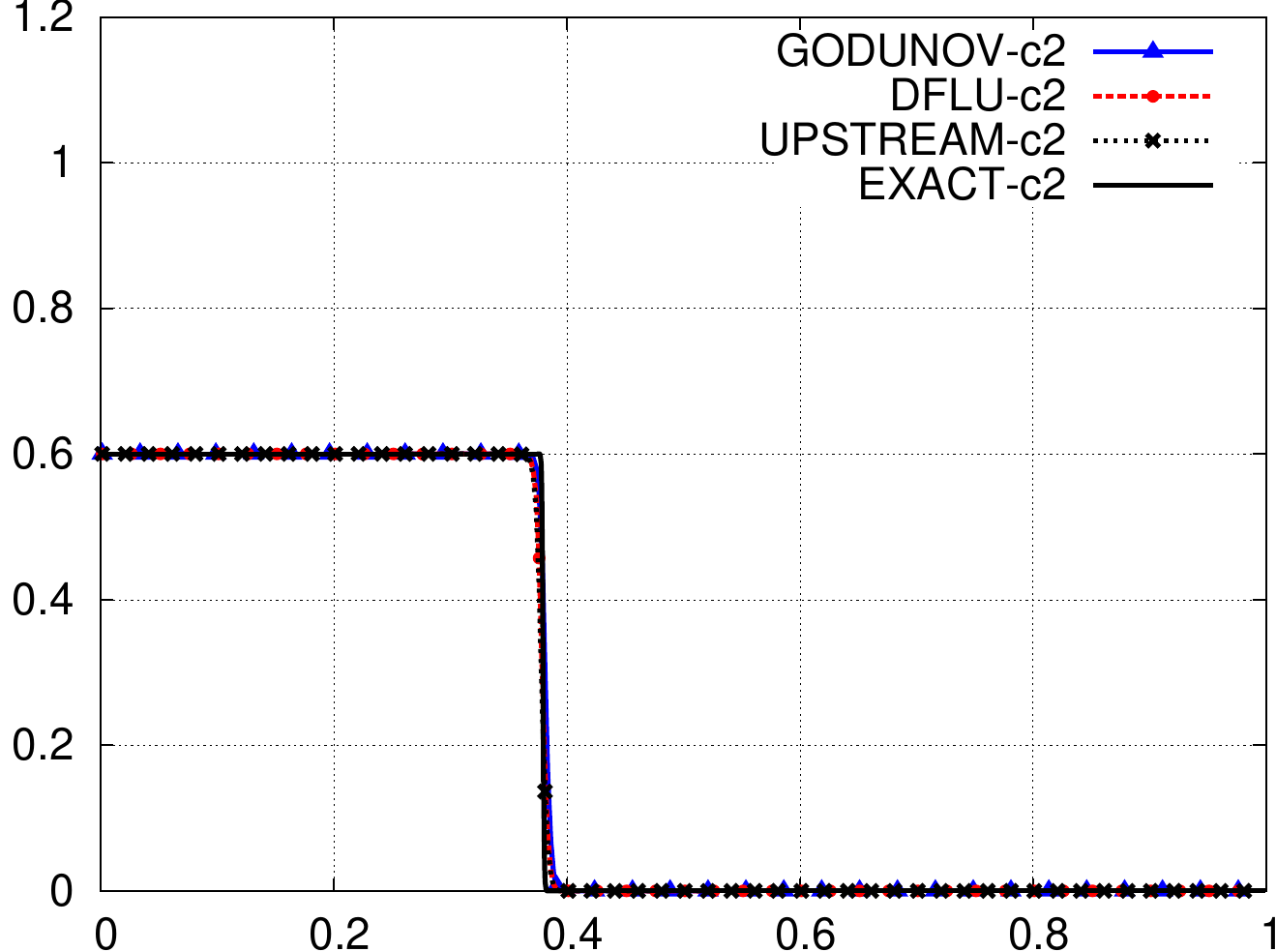}
\end{minipage}
\caption{Concentration $c_2$  at time $t=1$ with mesh size  $h=\frac{1}{200}$ (left) and $h=\frac{1}{400}$ (right), with high order accuracy.}
\label{expt5}
\end{figure}

Numerical experiments are done for DFLU flux, Upstream mobility flux and compared with Godunov flux.
In these experiments data are chosen so that DFLU flux differ from Godunov flux. 
The performance of the DFLU flux is as good as the Godunov flux.
High order accurate schemes corresponding to DFLU and Godunov  are constructed by introducing 
slope limiter in space variable and a strong stability preserving Runge-Kutta scheme
in the time variable, a comparison with first order scheme is shown in Fig.\ref{expt1},\ref{expt2}.
For first order and high order scheme it clearly shows that the DFLU flux is as good as the Godunov flux.
Note that Godunov flux requires the solution of the
Riemann problem  of a system where as DFLU flux requires the solution of the  Riemann problem of  a scalar equation.
Fig.\ref{expt3},\ref{expt4},\ref{expt5} shows that numerical solution computed by DFLU is as good as Godunov and converges faster than Upstream mobility scheme.
\par
For high order scheme the $L^1$ error , order of accuracy $\alpha$ for DFLU, Godunov and Upstream mobility schemes are given in the table \ref{}.
The order of accuracy $\alpha $ is calculated as follows:
$$e_1=\|s-s_{h_1}\|_{L^1}  \mbox{ with }h_1=h,\quad  e_2=\|s-s_{h_2}\|_{L^1} \mbox{ with } h_2=h/2,\quad \alpha=\frac{\ln(e_1/e_2)}{\ln 2}.$$
Note that as $h\rightarrow 0,$ $\alpha $ in DFLU is better than $\alpha $ in  Upstream and more close to $\alpha$ in GODUNOV.
Here the exact solutions $s,c=(c_1,c_2)$ are computed from Godunov scheme for very small values of $ h $ and $\Delta t$ with $\frac{\Delta t}{h}M=0.5.$
\section{2-D  model}\label{2d-2}
In this section we are extending the numerical schemes explained in \S~2 for one dimension to a multi dimensional space.For simplicity we explain only in two dimensions
and higher dimension can be handled in a similar way.
In dimension two the equation (\ref{eq:system}) can be rewritten as
\begin{equation}
\begin{array}{rll}
s_t + \frac{\partial F_1}{\partial x_1}(s,c,x)+\frac{\partial F_2}{\partial x_2}(s,c,x) &=& 0 \\
(sc_1+a_1(c_1))_t + \frac{\partial c_1F_1}{\partial x_1}(s,c,x)+\frac{\partial c_1F_2}{\partial x_2}(s,c,x) &=& 0\\
(sc_2+a_2(c_2))_t + \frac{\partial c_2F_1}{\partial x_1}(s,c,x)+\frac{\partial c_2F_2}{\partial x_2}(s,c,x) &=& 0\\
\end{array}
\label{eq1}
\end{equation}
where $ (x,t)\in \Omega\times(0,\infty), \,\,x=(x_1,x_2) $ and  the flux $F_1,F_2:[0,1] \times [0, c_0]^2 \times \Omega \rightarrow  \re$ are given by
\begin{equation}
F_1(s,c,x) = v_1(x)f(s,c), \quad f(s,c) = \frac{\lw(s,c)}{\lw(s,c) + \lo(s)}
\label{eq:ff}
\end{equation}
\begin{equation}
F_2(s,c,x) = [v_2(x) - (\rw-\ro)g \lo(s,c) K(x) ] f(s,c)
\label{eq:ff1}
\end{equation}\\
To compute $F_1$ and $F_2$ we need the velocity component $v=(v_1,v_2).$ This velocity (pressure) is governed by the incompressibility of the flow:
\begin{equation}
\nabla \cdot v = 0 \qquad \textrm{in } \quad \Omega
\label{eq:v1}
\end{equation}
with some suitable boundary condition for velocity (pressure) on $\partial \Omega$ as explained in \S~1.
\par 
Basic numerical approach for finite volume method is outlined in the following algorithm: 
\begin{enumerate}
\item Set time step $n=0$ and initialize $s^0, c^0=(c_1^0,c_2^0).$
\item Assume $s^{n}$ and $c^{n}=(c_1^n,c_2^n)$ are known  at $t=t_n$.
\item Solve for the pressure $p^n$  from ( \ref{eq:pr}) and (\ref{eq:v}).
\item Compute velocity $v^n$ from (\ref{eq:pr}).
\item Chose  time step $\Delta t^n$ so that 	CFL condition is satisfied see \S \ref{stability} .
\item Update saturation and concentration at $t=t_{n+1}$ level by 
\begin{eqnarray*}
s^{n+1} &=& s^n - \Delta t^n \nabla \cdot ( F(s^n,c^n,v^n)) \\
s^{n+1} {c_1}^{n+1}+a_1(c_1^{n+1}) &=& s^n {c_1}^n+a_1(c_1^{n}) - \Delta t^n \nabla \cdot ({c_1}^n F(s^n,c^n,v^n))\\
s^{n+1} {c_2}^{n+1}+a_1(c_2^{n+1}) &=& s^n {c_2}^n+a_2(c_2^{n+1}) - \Delta t^n \nabla \cdot ({c_2}^n F(s^n,c^n,v^n))
\end{eqnarray*}
\item Set $n=n+1$ and Go to step 2.
\end{enumerate}
\subsection{Discretization of the domain $\Omega=[0,1] \times [0,1]$}
\begin{figure}
\begin{center}
\includegraphics[width=0.3\textwidth]{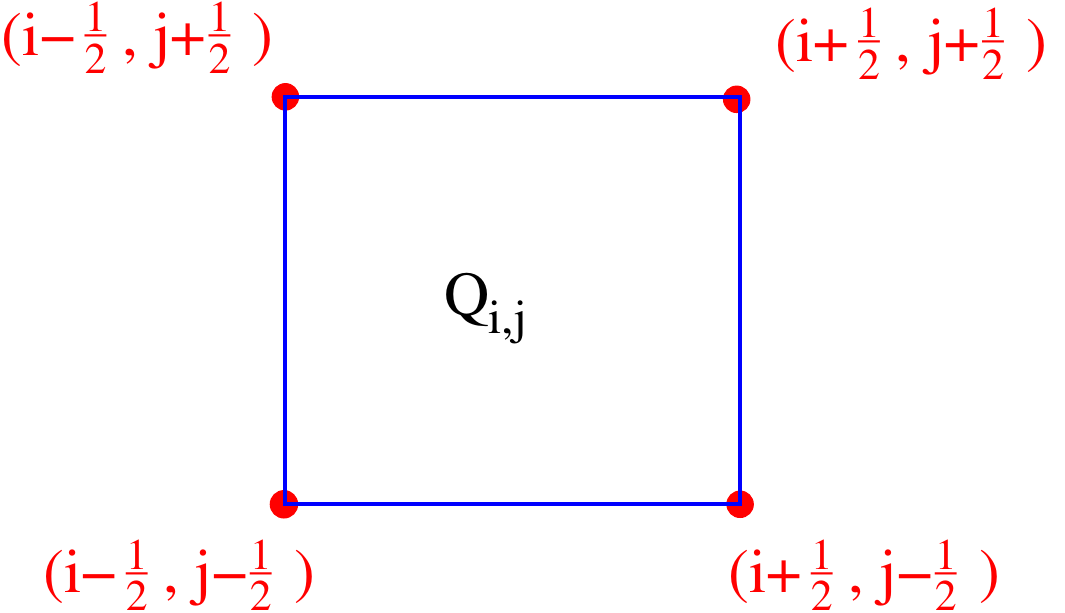}
\end{center}
\caption{Definition of cell $Q_{i,j}$ by grid points}
\label{fig:cell}
\end{figure}

Consider the Cartesian grid obtained by taking the cross product of the one-dimensional partitions 
$\{x_i, \ i=1,\ldots,n_x\}$ and $\{y_j, \ j=1,\ldots,n_y\}$ with $x_1=y_1=0$ and $x_{n_x}=y_{n_y}=1$. 
We also introduce one layer of grid points on all four sides of $\Omega$ which will be referred to as
ghost points. Thus the grid point indices range over $0 \le i \le n_x+1$ and $0 \le j \le n_y+1$. 
The grid defines the cell $Q_{i,j} = [x_{\imh},x_{\iph}] \times [y_{\jmh},y_{\jph}]$, see Fig.\ref{fig:cell},
for $0 \le i \le n_x$ and $0 \le j \le n_y$. The number of true cells where
the solution is supposed to be computed in the domain $\Omega$ is $n_c = (n_x-1) \times (n_y-1)$ (excluding the ghost cells).
\subsection{Numerical approximation for the  pressure }
Define $\mu := (\lw + \lo)K$ and  $\theta := (\lw \rw + \lo \ro) g K$. Integrating equation (\ref{eq:v}) over cell $Q_{i,j}$ and 
using the divergence theorem, we obtain the finite volume approximation
\begin{equation}
(v_{\iph,j} - v_{\imh,j}) \dy + (v_{i,\jph} - v_{i,\jmh}) \dx = 0
\label{eq:p}
\end{equation}
where the velocity at the cell face is given by
\begin{equation*}
v_{\iph,j} = \left. -\mu \df{p}{x} \right|_{\iph,j}, \quad v_{i,\jph} = \left. - ( \mu \df{p}{y}+\theta )\right|_{i,\jph}
\end{equation*}
We  approximate these  as follows:
\begin{itemize}
 \item {Along the x-direction}
\begin{eqnarray*}
p_{i+1,j} - p_{i,j} &=&  \int_{x_i}^{x_{i+1}} \mu \df{p}{x} \frac{1}{\mu} \ud x 
\approx \left. \mu \df{p}{x} \right|_{\iph,j} \int_{x_i}^{x_{i+1}} \frac{1}{\mu} \ud x \\
&\approx& \left. \mu \df{p}{x} \right|_{\iph,j} \frac{1}{2} \left( \frac{1}{\mu_{i,j}}  + \frac{1}{\mu_{i+1,j}} \right) \dx
\end{eqnarray*}
This leads to the following approximation for the velocity flux
\begin{equation}
v_{\iph,j} = - \bmu_{\iph,j} \frac{p_{i+1,j} - p_{i,j}}{\dx}, \qquad
\frac{1}{\bmu_{\iph,j}} = \frac{1}{2} \left(\frac{1}{\mu_{i,j}}  + \frac{1}{\mu_{i+1,j}} \right)
\label{eq:velx}
\end{equation}
\item {Along the y-direction}
\begin{eqnarray*}
p_{i,j+1} - p_{i,j} &=& \int_{y_j}^{y_{j+1}} \frac{1}{\mu} \left[ \mu \df{p}{y} +\theta - \theta \right] \ud y 
\approx -v_{i,\jph} \int_{y_j}^{y_{j+1}} \frac{1}{\mu} - \int_{y_j}^{y_{j+1}} \frac{\theta}{\mu} \\
&\approx& -v_{i,\jph} \frac{\dy}{2} \left( \frac{1}{\mu_{i,j}} + \frac{1}{\mu_{i,j+1}} \right) - \frac{\dy}{2} \left( \frac{\theta_{i,j}}{\mu_{i,j}} + \frac{\theta_{i,j+1}}{\mu_{i,j+1}} \right).
\end{eqnarray*}
Hence we get the approximation
\begin{equation}
v_{i,\jph} = - \bar{\mu}_{i,\jph} \frac{p_{i,j+1} - p_{i,j}}{\dy} - \bar{\theta}_{i,\jph}
\label{eq:velyg}
\end{equation}
where
\begin{equation}
\frac{1}{\bar{\mu}_{i,\jph}} := \half \left( \frac{1}{\mu_{i,j}} + \frac{1}{\mu_{i,j+1}} \right), \qquad \bar{\theta}_{i,\jph} := \frac{\bar{\mu}_{i,\jph}}{2} \left( \frac{\theta_{i,j}}{\mu_{i,j}} + \frac{\theta_{i,j+1}}{\mu_{i,j+1}} \right). 
\end{equation}
\end{itemize}
The velocity on the inlet boundary is computed as
\begin{equation}
v_{\half,j} = -\bmu_{\half,j} \frac{p_{1,j} - p_I}{\dx}, \quad \frac{1}{\bmu_{\half,j}} = \frac{1}{2} \left(\frac{1}{\mu_{0,j}}  + \frac{1}{\mu_{1,j}} \right)
\end{equation}
with similar expressions for the other inlet/outlet parts of the boundary. On the rest of the boundary, the normal velocity is zero
which is equivalent to saying that flux is zero. The system of equations (\ref{eq:p}) for the pressure can be put in the form
\begin{equation}
Ap=b
\end{equation}
where $A \in \re^{n_c \times n_c}$ and $b \in \re^{n_c}$. This matrix equation is solved using conjugate the gradient method. 
\subsection{Finite volume scheme} 
By integrating equations in (\ref{eq:system}) over the cell $Q_{i,j}$, we obtain the following finite volume approximations 
\begin{align}
s^{n+1}_{i,j}& = s^n_{i,j} - \frac{\dt}{\dx \dy} \left\{ [F^n_{\iph,j} - F^n_{\imh,j}] \dy + [F^n_{i,\jph} - F^n_{i,\jmh}] \dx \right\},\label{eq:fvm-1} \\
s^{n+1}_{i,j} {c_1}^{n+1}_{i,j}+a_1({c_1}^{n+1}_{i,j})& = s^n_{i,j} {c_1}^n_{i,j} +a_1({c_1}^n_{i,j})- \frac{\dt}{\dx \dy} \left\{ [(c_1F)^n_{\iph,j} - (c_1F)^n_{\imh,j}] \dy  \right. \notag \\
 &\hspace{5 cm}\left. +[(c_1F)^n_{i,\jph} - (c_1F)^n_{i,\jmh}] \dx \right\},\label{eq:fvm-2}\\
s^{n+1}_{i,j} {c_2}^{n+1}_{i,j}+a_2({c_2}^{n+1}_{i,j}) &= s^n_{i,j} {c_2}^n_{i,j}+a_2({c_2}^n_{i,j}) - \frac{\dt}{\dx \dy} \left\{ [(c_2F)^n_{\iph,j} -
(c_2F)^n_{\imh,j}] \right.\dy \notag\\
&\hspace{5 cm}\left.+ [(c_2F)^n_{i,\jph} - (c_1F)^n_{i,\jmh}] \dx \right\}.\label{eq:fvm-3}
\end{align}
Here we introduce the DFLU numerical flux for two dimensional finite volume scheme by using the idea explained in \S 2.
The corresponding numerical fluxes 
are given by 
\begin{eqnarray}\nonumber 
 F^{n}_{\iph,j} &=&\max\{ F_1(\max(s^n_{i,j},(\theta^n_{F_1})_{i,j}),c^n_{i,j},K_{i,j} ),{F_1}(\min(s^n_{i+1,j},(\theta^n_{F_1})_{i+1,j}),c^n_{i+1,j},K_{i+1,j} ) \}
\end{eqnarray}
\begin{eqnarray}\nonumber
 F^{n}_{i,\jph} &=&\max\{ {F_2}(\max(s^n_{i,j},(\theta^n_{F_2})_{i,j}),c^n_{i,j},K_{i,j} ),{F_2}(\min(s^n_{i,j+1},(\theta^n_{F_2})_{i,j+1}),c^n_{i,j+1},K_{i,j+1} ) \}
\end{eqnarray}
where $c^n_{i,j}=({c_1}^n_{i,j},{c_2}^n_{i,j})$ and 
$$  (\theta^n_{F_1})_{i,j}=\mbox{argmin} \,\,{F_1}(. ,c^n_{i,j},K_{i,j}) \mbox{ and } (\theta^n_{F_2})_{i,j}=\mbox{argmin} \,\,{F_2}(. ,c^n_{i,j},K_{i,j}) $$
\begin{equation}\nonumber
(c_l F)^{n}_{\iph,j} = \begin{cases}
{c_l}_{i,j} F^{n}_{\iph,j} & \textrm{if } F^{n}_{\iph,j} > 0 \\
{c_l}_{i+1,j}  F^{n}_{\iph,j} & \textrm{if } F^{n}_{\iph,j} \le 0 \quad l=1,2.
\end{cases}
\end{equation}
\subsection{High-order scheme}
In order to develop the second order scheme, we follow the method of lines approach in which space and time discretization are performed separately. In the first step, spatial discretization using piecewise linear reconstruction is made which leads to a system of ODE which can be written as
\begin{equation}
\dd{U}{t} + R(U) = 0, \quad U=\begin{bmatrix} s \\ sc_1+a_1(c_1) \\sc_2+a_2(c_2) \end{bmatrix}
\label{eq:ode}
\end{equation}
where
\begin{equation}
R(U)_{i,j} = \frac{1}{\dx \dy}  \begin{bmatrix}
[F_{\iph,j} - F_{\imh,j}] \dy + [F_{i,\jph} - F_{i,\jmh}] \dx \\
[{G_1}_{\iph,j} -{G_1} _{\imh,j}] \dy + [{G_1}_{i,\jph} - {G_1}_{i,\jmh}] \dx\\
[{G_2}_{\iph,j} - {G_2}_{\imh,j}] \dy + [{G_2}_{i,\jph} - {G_2}_{i,\jmh}] \dx
\end{bmatrix}
\end{equation}
The high order accurate fluxes are given by
 $$F_{i+\half,j}=\bar F(s^{L}_{i+\half,j},s^{R}_{i+\half,j},{c_1}^{L}_{i+\half,j},{c_2}^{L}_{i+\half,j},{c_1}^{R}_{i+\half,j},{c_2}^{R}_{i+\half,j},v_{\iph,j},K_{i,j},K_{i+1,j})$$
\begin{equation}
{G_l}_{\iph,j} = \begin{cases}
{c_l}^L_{\iph,j} F_{\iph,j} & \textrm{if } F_{\iph,j} > 0 \\
{c_l}^R_{\iph,j}  F_{\iph,j} & \textrm{if } F_{\iph,j} \le 0\quad l=1,2,
\end{cases}
\end{equation}
and similar expression for $F_{i,j+\half}$. The quantities with superscripts $L$ and $R$ denote the reconstructed 
values of the variables to the left and right of the cell face. For any quantity $u$, we can define the reconstruction in $x-$direction  as follows:
\begin{equation}
u^L_{\iph,j} = u_{i,j} + \half \delta^x_{i,j}, \qquad
u^R_{\iph,j} = u_{i+1,j} - \half \delta^x _{i+1,j}
\end{equation}
where
\begin{equation}
\label{minimod}
\delta^x_{i,j} = \textrm{minmod}\left( \theta (u_{i,j} - u_{i-1,j}), \half (u_{i+1,j} - u_{i-1,j}), \theta (u_{i+1,j} - u_{i,j}) \right), \quad \theta \in [1,2].
\end{equation} Similarly in the $y-$direction we can define $u^L_{i,\jph} $  and $u^R_{i,\jph}.$

\subsection{Stability results }\label{stability}
Let us write $$\bar{s}^n=(\bar{s}_i^n)^8_{i=1}=(s^{n,L}_{\imh,j},s^{n,R}_{\imh,j}
,s^{n,L}_{\iph,j},s^{n,R}_{\iph,j},s^{n,L}_{i,\jmh},s^{n,R}_{i,\jmh},s^{n,L}_{i,\jph},s^{n,R}_{i,\jph})$$
$$\bar{c}^n=(\bar{c}_i^n)^{16}_{i=1}=({c_l}^{n,L}_{\imh,j},{c_l}^{n,R}_{\imh,j}
,{c_l}^{n,L}_{\iph,j},{c_l}^{n,R}_{\iph,j},{c_l}^{n,L}_{i,\jmh},{c_l}^{n,R}_{i,\jmh},{c_l}^{n,L}_{i,\jph},{c_l}^{n,R}_{i,\jph})_{l=1,2}$$
The updated value of the saturation (\ref{eq:fvm-1}) can be written as
$$s^{n+1}_{i,j}=H(\bar{s}^n,\bar{c}^n,v_{i\pm\half,j},v_{i,j\pm\half},K_{i\pm1,j},K_{i,j},K_{i,j\pm1}).$$
Where $H$ is Lipschitz continuous in saturation and concentration with the property
\begin{eqnarray*}
 H({\bf 0},v_{i\pm\half,j},v_{i,j\pm\half},K_{i\pm1,j},K_{i,j},K_{i,j\pm1})=0\\
 H({\bf 1},v_{i\pm\half,j},v_{i,j\pm\half},K_{i\pm1,j},K_{i,j},K_{i,j\pm1})=1.
\end{eqnarray*}
Since the slope limiter preserves the average value of the solution in each cell, we can express this as
\begin{equation}
\label{finalscheme}
s^{n+1}_{i,j}= \frac{s^{n,L}_{\iph,j}+s^{n,R}_{\imh,j}}{4}+\frac{s^{n,L}_{i,\jph}
+s^{n,R}_{i,\jmh}}{4}-\frac{\Delta t}{\Delta x}(F^n_{i+\half,j}-F^n_{i-\half,j} )-\frac{\Delta t}{\Delta y}(F^n_{i,j+\half}-F^n_{i,j-\half} ).
\end{equation}
If we differentiate $H$ with respect to  its variables $\bar{s}^n_i$ we  can observe that $\df{}{\bar{s}^n_i}H\geq0$ provided 
\begin{equation}
\lambda^x|\df{}{\bar{s}^n_i}F^n_{i\pm\half,j}|,\lambda^y|\df{}{\bar{s}^n_i}F^n_{i,j\pm \half}|\leq \frac{1}{4}.
\label{cfl-1}
\end{equation}
Let
$$ M=\sup_{s} \{ \df{F_1}{s},\df{F_2}{s}, \frac{F_1}{s+h_l}, \frac{F_2}{s+h_l} \} , $$
then the  condition (\ref{cfl-1}) reduces to,  
\begin{equation}
\max\{\lambda^x M,\lambda^y M\}\leq \frac{1}{4}, \mbox{ where } \lambda^x=\frac{\Delta t}{\Delta x},\,\, \lambda^y=\frac{\Delta t}{\Delta y}.
 \label{cfl}
\end{equation}
This shows that $H$ is monotone in each of its variable. Using these facts we have the following lemmas. 
\begin{lemma} 
 Let $s_0 \in [0,1]$ be the initial data and let $\{s^n_{i,j}\}$ be the corresponding solution calculated by the finite
volume scheme (\ref{eq:fvm-1}) using DFLU flux along with slope limiter. If the CFL given in (\ref{cfl}) holds  then 
\begin{equation} 
\label{a}
 0\leq s^n_{i,j} \leq 1 \,\,\,\forall,\, i,j \mbox{ and } n.
\end{equation}
\end{lemma}
{\bf Proof:}
 From the property of slope limiter we can observe that whenever $ 0\leq s^n_{i,j} \leq 1 $ then 
the reconstructed values satisfies 
$$ 0\leq s^{n,L,R}_{\ipmh,j},s^{n,L,R}_{i,\jpmh}\leq1  \,\,\forall \,\,i,j \mbox{ and }n.$$
Using this property and the monotonicity of the  $H$ , we get 
$$
\begin{array}{lll}
0=H({\bf{0}},\bar{c}^n,v_{i\pm\half,j},v_{i,j\pm\half},K_{i\pm1,j},K_{i,j},K_{i,j\pm1} )&&\\
\hspace{3  cm}\leq H(\bar{s}^n,\bar{c}^n,v_{i\pm\half,j},v_{i,j\pm\half},K_{i\pm1,j},K_{i,j},K_{i,j\pm1})=s^{n+1}_{i,j} &&\\
\hspace{3 cm} \leq H({\bf{1}},\bar{c}^n,v_{i\pm\half,j},v_{i,j\pm\half},K_{i\pm1,j},K_{i,j},K_{i,j\pm1})=1&&
\end{array}
$$
This proves that $$0\leq s^{n+1}_{i,j}\leq 1 \,\, \forall\,\, i,j\mbox.$$
\hfill \qed

 Now we prove the lemma that gives the maximum principle for the concentration.
\begin{lemma}{\label{th1}}
 Let $\{{c^n_1}_{i,j}\}$,$\{{c^n_2}_{i,j}\}$ be the solution calculated by the finite
volume scheme (\ref{eq:fvm-2}) and (\ref{eq:fvm-3}) by using DFLU flux with slope limiter. Under the CFL condition (\ref{cfl})  concentration $c=(c_1,c_2)$ satisfies the followig 
maximum principle 
\begin{align*}
 &(a)\,\min\{{c^n_l}_{i,j},{c^n_l}_{i \pm 1,j},{c^n_l}_{i,j\pm 1}\}\leq {c_l}^{n+1}_{i,j}
\leq \max\{{c^n_l}_{i,j},{c^n_l}_{i \pm 1,j},{c^n_l}_{i,j\pm 1 }\},\\
&\hspace{7 cm} \forall \,\,n\in\Z^+,\,\,i\in\Z,l=1,2.\\
\end{align*}
\end{lemma}
{\bf Proof:}
We can express the high order  numerical fluxes ${G_l}_{\iph,j}, {G_l}_{i,\jph} (l=1,2)$ in the finite volume scheme (\ref{eq:fvm-2}) and (\ref{eq:fvm-3}) as
\begin{align*}
{G_l}_{\iph,j}&=c^{nL}_{l_{\iph,j}}F^+_{\iph,j}+c^{nR}_{l_{\iph,j}}F^-_{\iph,j}\\
{G_l}_{i,\jph}&=c^{nL}_{l_{i,\jph}}F^+_{i,\jph}+c^{nR}_{l_{i,\jph}}F^-_{i,\jph},
\end{align*}
where  
\begin{align*}
F^+_{\iph,j}=\max\{F_{\iph,j},0\},\quad F^-_{\iph,j}=\min\{F_{\iph,j},0\}\\
F^+_{i,\jph}=\max\{F_{i,\jph},0\},\quad F^-_{i,\jph}=\min\{F_{i,\jph},0\}.
\end{align*}
We write the scheme (\ref{eq:fvm-2}) and (\ref{eq:fvm-3})  $(l=1,2)$ as 
\begin{align*}
s^{n+1}_{i,j}c^{n+1}_{l_{i,j}}+a_l(c^{n+1}_{l_{i,j}})-s^n_{i,j}c^n_{l_{i,j}}-a_l(c^{n}_{l_{i,j}})+\lx( c^L_{l_{\iph,j}}F^+_{\iph,j}+c^R_{l_{\iph,j}}F^-_{\iph,j}
-c^L_{l_{\imh,j}}F^+_{\imh,j}-c^R_{l_{\imh,j}}F^-_{\imh,j})\\
+\ly(c^L_{l_{i,\jph}}F^+_{l_{i,\jph}}+c^R_{l_{i,\jph,j}}F^-_{i,\jph}-c^L_{l_{i,\jmh}}F^+_{i,\jmh}-c^R_{l_{i,\jmh,j}}F^-_{i,\jmh}  )=0.
\end{align*}
By adding and subtracting the terms  $ s^{n+1}_{i,j}c^n_{i,j} $ we get
\begin{align*}
 &(s^{n+1}_{i,j}+a'_l(\zeta^{n+\half}_i))(c^{n+1}_{l_{i,j}}-c^n_{l_{i,j}})+c^n_{l_{i,j}}(s^{n+1}_{i,j}-s^n_{i,j})\\
 &\hspace{2 cm}+\lx( c^{nL}_{l_{\iph,j}}F^+_{\iph,j}+c^{nR}_{l_{\iph,j}}F^-_{\iph,j}-c^{nL}_{l_{\imh,j}}F^+_{\imh,j}-c^{nR}_{l_{\imh,j}}F^-_{\imh,j})\\
&\hspace{2 cm}+\ly(c^{nL}_{l_{i,\jph}}F^+_{i,\jph}+c^{nR}_{l_{i,\jph}}F^-_{i,\jph}-c^{nL}_{l_{i,\jmh}}F^+_{i,\jmh}-c^{nR}_{l_{i,\jmh}}F^-_{i,\jmh}  )=0.
\end{align*}
where $a_l({c_l}^{n+1}_{i,j})-a_l({c_l}_{i,j}^{n})=a'_l(\zeta^{n+\half}_{i,j})({c_l}^{n+1}_{i,j}-{c_l}^{n}_{i,j}),$ for some $\zeta^{n+\half}_{i,j} $ between 
${c_l}^{n+1}_{i,j}$ and ${c_l}^{n}_{i,j}.$
Replacing $s^{n+1}_{i,j}-s^n_{i,j}$ by $-\lx (F_{\iph,j}-F_{\imh,j})-\ly(F_{i,\jph}-F_{i,\jmh})$ and splitting 
 $F_{\ipmh,j}$ by $F^+_{\ipmh,j}+F^-_{\ipmh,j}$ (similarly for $F_{i,\jpmh}$ ) and by rearranging the terms we  have
 \begin{align*}
 & (s^{n+1}_{i,j}+a'_l(\zeta^{n+\half}_{i,j}))(c^{n+1}_{l_{i,j}}-c^n_{l_{i,j}})+\lx F^+_{\iph,j}(c^{nL}_{l_{\iph,j}}-c^n_{l_{i,j}})+\lx F^-_{\iph,j}(c^{nR}_{l_{\iph,j}}-c^n_{l_{i,j}})\\
&\hspace{2.7 cm}+\lx F^+_{\imh,j}(c^n_{l_{i,j}}-c^{nL}_{l_{\imh,j}})+\lx F^-_{\imh,j}(c^n_{l_{i,j}}-c^{nR}_{l_{\imh,j}}) \\
&\hspace{2.7 cm}+\ly F^+_{i,\jph}(c^{nL}_{i,\jph}-c^n_{i,j})+\ly F^-_{i,\jph}(c^{nR}_{i,\jph}-c^n_{i,j})\\
&\hspace{2.7 cm}+\ly F^+_{i,\jmh}(c^n_{l_{i,j}}-c^{nL}_{l_{i,\jmh}})+\ly F^-_{i,\jmh}(c^n_{l_{i,j}}-c^{nR}_{l_{i,\jmh}})=0.  
 \end{align*}

 Note that 
 $$c^{nL}_{l_{\iph,j}}=c^n_{l_{i,j}}+\frac{\delta^x_{i,j}}{2},\quad c^{nR}_{l_{\imh,j}}=c^n_{l_{i,j}}-\frac{\delta^x_{i,j}}{2},\quad
 c^{nL}_{l_{i,\jph}}=c^n_{l_{i,j}}+\frac{\delta^y_{i,j}}{2},\quad c^{nR}_{l_{i,\jmh}}=c^n_{l_{i,j}}-\frac{\delta^y_{i,j}}{2}
 $$
and $\delta^x_{i,j}$ and $\delta^y_{i,j}$ are the slope limiter given by
\begin{align*}
 \delta^x_{i,j} = \textrm{minmod}\left( \theta(c^n_{l_{i,j}} - c^n_{l_{i-1,j}}), \half (c^n_{l_{i+1,j}} - c^n_{l_{i-1,j}}), \theta(c^n_{l_{i+1,j}} - c^n_{l_{i,j}}) \right),\\
  \delta^y_{i,j} = \textrm{minmod}\left( \theta(c^n_{l_{i,j}} - c^n_{l_{i,j-1}}), \half (c^n_{l_{i,j+1}} - c^n_{l_{i,j-1}}), \theta(c^n_{l_{i,j+1}} - c^n_{l_{i,j}}) \right).
\end{align*}
After substituting the values of $c^{nL}_{l_{\ipmh,j}},\, c^{nR}_{l_{\ipmh,j}},\,c^{nL}_{l_{i,\jpmh}},\, c^{nR}_{l_{i,\jpmh}}$ the above equation becomes 
\begin{align*}
&c^{n+1}_{l_{i,j}}=c^n_{l_{i,j}}-\lx \frac{F^+_{\iph,j}}{(s^{n+1}_{i,j}+a'_l(\zeta^{n+\half}_{i,j}))(c^n_{l_{i,j}}-c^n_{l_{i-1,j}})} \frac{\delta^x_{i,j}}{2}(c^n_{l_{i,j}}-c^n_{l_{i-1,j}})\\
&\hspace{1.8 cm}-\lx \frac{F^-_{\iph,j}}{(s^{n+1}_{i,j}+a'_l(\zeta^{n+\half}_{i,j}))}(1-\frac{\delta^x_{i+1,j}}{2(c^n_{l_{i+1,j}}-c^n_{l_{i,j}})})(c^n_{l_{i+1,j}}-c^n_{l_{i,j}})\\
&\hspace{1.8 cm}-\lx \frac{F^+_{\imh,j}}{(s^{n+1}_{i,j}+a'_l(\zeta^{n+\half}_{i,j}))}(1-\frac{\delta^x_{i-1,j}}{2(c^n_{l_{i,j}}-c^{n}_{l_{i-1,j}})})(c^n_{l_{i,j}}-c^{n}_{l_{i-1,j}})\\
&\hspace{1.8 cm}-\lx \frac{F^-_{\imh,j}}{(s^{n+1}_{i,j}+a'_l(\zeta^{n+\half}_{i,j}))(c^n_{l_{i+1,j}}-c^n_{l_{i,j}})}\frac{\delta^x_{i,j}}{2}(c^n_{l_{i+1,j}}-c^n_{l_{i,j}})\\   
&\hspace{1.8 cm}-\ly\frac{F^+_{i,\jph}}{(s^{n+1}_{i,j}+a'_l(\zeta^{n+\half}_{i,j}))(c^n_{l_{i,j}}-c^n_{l_{i,j-1}})} \frac{\delta^y_{i,j}}{2}(c^n_{l_{i,j}}-c^n_{l_{i,j-1}})\\
&\hspace{1.8 cm}-\ly \frac{F^-_{i,\jph}}{(s^{n+1}_{i,j}+a'_l(\zeta^{n+\half}_{i,j}))}(1-\frac{\delta^y_{i,j+1}}{2(c^n_{l_{i,j+1}}-c^n_{l_{i,j}})})(c^n_{l_{i,j+1}}-c^n_{l_{i,j}})\\
&\hspace{1.8 cm}-\ly \frac{ F^+_{i,\jmh}}{(s^{n+1}_{i,j}+a'_l(\zeta^{n+\half}_{i,j}))}(1-\frac{\delta^y_{i,j-1}}{2(c^n_{l_{i,j}}-c^{n}_{l_{i,j-1}})})(c^n_{l_{i,j}}-c^{n}_{l_{i,j-1}})\\
&\hspace{1.8 cm}-\ly \frac{ F^-_{i,\jmh}}{(s^{n+1}_{i,j}+a'_l(\zeta^{n+\half}_{i,j}))(c^n_{l_{i,j+1}}-c^n_{l_{i,j}})}\frac{\delta^y_{i,j}}{2}(c^n_{l_{i,j+1}}-c^n_{l_{i,j}}),   
\end{align*}
Now we write it as 
\begin{align}
 c^{n+1}_{l_{i,j}}&=c^n_{l_{i,j}}-\alpha^1_{\imh,j}(c^n_{l_{i,j}}-c^n_{l_{i-1,j}})
+\alpha^2_{\iph,j}(c^n_{l_{i+1,j}}-c^n_{l_{i,j}})-\alpha^3_{\imh,j}(c^n_{l_{i,j}}-c^{n}_{l_{i-1,j}})\notag\\
&\hspace{1.1 cm}+\alpha^4_{\iph,j}(c^n_{l_{i+1,j}}-c^n_{l_{i,j}})-\alpha^1_{i,\jmh}(c^n_{l_{i,j}}-c^n_{l_{i,j-1}})+\alpha^2_{l_{i,\jph}}(c^n_{l_{i,j+1}}-c^n_{l_{i,j}})\notag\\
&\hspace{1.1 cm}-\alpha^3_{i,\jmh}(c^n_{l_{i,j}}-c^{n}_{l_{i,j-1}})+\alpha^4_{i,\jph}(c^n_{l_{i,j+1}}-c^n_{l_{i,j}})\notag\\
&=c^n_{l_{i,j}}-C^n_{\imh,j}(c^n_{l_{i,j}}-c^n_{l_{i-1,j}})+D^n_{\iph,j}(c^n_{l_{i+1,j}}-c^n_{l_{i,j}})\notag\\
&\hspace{1.05 cm}-C^n_{i,\jmh}(c^n_{l_{i,j}}-c^n_{l_{i,j-1}})+D^n_{i,\jph}(c^n_{l_{i,j+1}}-c^n_{l_{i,j}})\label{harten},
\end{align}
where
\begin{align*}
 C^n_{\imh,j}=\alpha^1_{\imh,j}+\alpha^3_{\imh,j},\quad D^n_{\iph,j}=\alpha^2_{\iph,j}+\alpha^4_{\iph,j}\\
 C^n_{i,\jmh}=\alpha^1_{i,\jmh}+\alpha^3_{i,\jmh},\quad D^n_{i,\jph}=\alpha^2_{i,\jph}+\alpha^4_{i,\jph}
\end{align*}
and 
\begin{align*}
 \alpha^1_{\imh,j}&=\lx \frac{F^+_{\iph,j}}{(s^{n+1}_{i,j}+a'_l(\zeta^{n+\half}_{i,j}))(c^n_{l_{i,j}}-c^n_{l_{i-1,j}})} \frac{\delta^x_{i,j}}{2},\\
 \alpha^2_{\iph,j}&=-\lx \frac{F^-_{\iph,j}}{(s^{n+1}_{i,j}+a'_l(\zeta^{n+\half}_{i,j}))}(1-\frac{\delta^x_{i+1,j}}{2(c^n_{l_{i+1,j}}-c^n_{l_{i,j}})}),\\
 \alpha^3_{\imh,j}&=\lx \frac{F^+_{\imh,j}}{(s^{n+1}_{i,j}+a'_l(\zeta^{n+\half}_{i,j}))}(1-\frac{\delta^x_{i-1,j}}{2(c^n_{l_{i,j}}-c^{n}_{l_{i-1,j}})}),\\
 \alpha^4_{\iph,j}&=-\lx \frac{F^-_{\imh,j}}{(s^{n+1}_{i,j}+a'_l(\zeta^{n+\half}_{i,j}))(c^n_{l_{i+1,j}}-c^n_{l_{i,j}})}\frac{\delta^x_{i,j}}{2},\\
 \alpha^1_{i,\jmh}&=\ly\frac{F^+_{i,\jph}}{(s^{n+1}_{i,j}+a'_l(\zeta^{n+\half}_{i,j}))(c^n_{l_{i,j}}-c^n_{l_{i,j-1}})} \frac{\delta^y_{i,j}}{2},\\
 \alpha^2_{i,\jph}&=-\ly \frac{F^-_{i,\jph}}{(s^{n+1}_{i,j}+a'_l(\zeta^{n+\half}_{i,j}))}(1-\frac{\delta^y_{i,j+1}}{2(c^n_{l_{i,j+1}}-c^n_{l_{i,j}})}),\\
 \alpha^3_{i,\jmh}&=\ly \frac{ F^+_{i,\jmh}}{(s^{n+1}_{i,j}+a'_l(\zeta^{n+\half}_{i,j}))}(1-\frac{\delta^y_{i,j-1}}{2(c^n_{l_{i,j}}-c^{n}_{l_{i,j-1}})}),\\
 \alpha^4_{i,\jph}&=-\ly \frac{ F^-_{i,\jmh}}{(s^{n+1}_{i,j}+a'_l(\zeta^{n+\half}_{i,j}))(c^n_{l_{i,j+1}}-c^n_{l_{i,j}})}\frac{\delta^y_{i,j}}{2}. 
\end{align*}
From the property of the limiter it is easy to see that 
\begin{align*}
 0\leq \frac{\delta^x_{i+1,j}}{2(c^n_{l_{i+1,j}}-c^n_{l_{i,j}})}, \quad \frac{\delta^y_{i+1,j}}{2(c^n_{l_{i,j+1}}-c^n_{l_{i,j}})} \leq 1,
\end{align*}
which in turn implies that 
\begin{align*}
 C^n_{\imh,j},C^n_{i,\jmh}, D^n_{\iph,j},D^n_{i,\jph}\geq 0.
\end{align*}
Now we prove the maximum principle through following cases.
\par 
{\bf Case 1:}
 Suppose that \\
 (a) $c^n_{l_{i,j}}$ lies between $ c^n_{l_{i-1,j}},$ and $c^n_{l_{i+1,j}}$ and \\ 
 (b) $c^n_{l_{i,j}}$ lies between $c^n_{l_{i,j+1}},$ and $c^n_{l_{i,j-1}},$ then
 \begin{align}
  c^n_{l_{i,j}}&=\theta^xc^n_{l_{i-1,j}}+(1-\theta^x)c^n_{l_{i+1,j}} \quad \mbox{ for some } \theta^x \in [0,1] \mbox{ and }\label{con-1}\\
  c^n_{l_{i,j}}&=\theta^y c^n_{l_{i,j-1}}+(1-\theta^y)c^n_{l_{i,j+1}}\quad \mbox{ for some } \theta^y \in [0,1].\label{con-2} 
 \end{align}
Now 
\begin{align*}
 c^n_{l_{i,j}}-c^n_{l_{i-1,j}}&=(1-\theta^x)(c^n_{l_{i+1,j}}-c^n_{l_{i-1,j}}),\\
 c^n_{l_{i+1,j}}-c^n_{l_{i,j}}&=\theta^x(c^n_{l_{i+1,j}}-c^n_{l_{i-1,j}}),\\
 c^n_{l_{i,j}}-c^n_{l_{i,j-1}}&=(1-\theta^y)(c^n_{l_{i,j+1}}-c^n_{l_{i,j-1}}),\\
 c^n_{l_{i,j+1}}-c^n_{l_{i,j}}&=\theta^y(c^n_{l_{i,j+1}}-c^n_{l_{i,j-1}}).
\end{align*}
By writing $c^n_{l_{i,j}}$ as $\tfrac{1}{2} (c^n_{l_{i,j}}+c^n_{l_{i,j}})$ and substituting the values from (\ref{con-1}) and (\ref{con-2}) the 
equation (\ref{harten}) becomes
\begin{align}
 c^{n+1}_{l_{i,j}}&=\tfrac{1}{2}(\theta^xc^n_{l_{i-1,j}}+(1-\theta^x)c^n_{l_{i+1,j}})
-C^n_{\imh,j}(1-\theta^x)(c^n_{l_{i+1,j}}-c^n_{l_{i-1,j}})\notag\\
&\hspace{0.5 cm}+D^n_{\iph,j}\theta^x(c^n_{l_{i+1,j}}-c^n_{l_{i-1,j}})+\tfrac{1}{2}(\theta^y c^n_{l_{i,j-1}}+(1-\theta^y)c^n_{l_{i,j+1}})\notag\\
&\hspace{0.5 cm}-C^n_{i,\jmh}(1-\theta^y)(c^n_{l_{i,j+1}}-c^n_{l_{i,j-1}})+D^n_{i,\jph}\theta^y(c^n_{l_{i,j+1}}-c^n_{l_{i,j-1}}),\notag\\
&=\lambda_1c^n_{l_{i-1,j}}+\lambda_2c^n_{l_{i+1,j}}+\lambda_3c^n_{l_{i,j-1}}+\lambda_4c^n_{l_{i,j+1}},\label{convex-com}
\end{align}
where 
\begin{align*}
\lambda_1=(1-\theta^x)C^n_{\imh,j}+\theta^x(\tfrac{1}{2}-D^n_{\iph,j}),\\
\lambda_2=(1-\theta^x)(\tfrac{1}{2}-C^n_{\imh,j})+\theta^xD^n_{\iph,j}\\
\lambda_3=(1-\theta^y)C^n_{i,\jmh}+\theta^y(\tfrac{1}{2}-D^n_{i,\jph}),\\
\lambda_4=(1-\theta^y)(\tfrac{1}{2}-C^n_{i,\jmh})+\theta^yD^n_{i,\jph}.
\end{align*}
Note that $\lambda_1+\lambda_2+\lambda_3+\lambda_4=1$ and 
with the CFL condition (\ref{cfl}) we have $$C^n_{\imh,j},C^n_{i,\jmh},D^n_{\iph,j},D^n_{i,\jph}\leq \tfrac{1}{2}$$
which gives $\lambda_1,\lambda_2,\lambda_3,\lambda_4\geq0.$ Hence from (\ref{convex-com}) the maximum principle follows.
\par
{\bf Case2:}
Suppose that \\
(a) $c^n_{l_{i,j}}$ does not lie between $ c^n_{l_{i-1,j}},$ and $c^n_{l_{i+1,j}}$ and \\
(b) $c^n_{l_{i,j}}$ does not lie between  $c^n_{l_{i,j+1}},$ and $c^n_{l_{i,j-1}},$
then we have $\delta^x_{i,j}=\delta^y_{i,j}=0.$ i.e.,
\begin{align*}
 C^n_{\imh,j}=\alpha^3_{\imh,j},D^n_{\iph,j}=\alpha^2_{\iph,j},C^n_{i,\jmh}=\alpha^3_{i,\jmh} \mbox{ and } D^n_{i,\jph}=\alpha^2_{i,\jph}.
\end{align*}
The equation (\ref{harten}) can be written as 
\begin{align*}
 c^{n+1}_{l_{i,j}}&=(1-C^n_{\imh,j}-D^n_{\iph,j}-C^n_{i,\jmh}-D^n_{i,\jph})c^n_{l_{i,j}}\\
 &\hspace{1 cm}+C^n_{\imh,j}c^n_{l_{i-1,j}}+D^n_{\iph,j}c^n_{l_{i+1,j}}
 +C^n_{i,\jmh}c^n_{l_{i,j-1}}+D^n_{i,\jph}c^n_{l_{i,j+1}}.
\end{align*}
Note that 
\begin{align*}
 C^n_{\imh,j}+D^n_{\iph,j}+C^n_{i,\jmh}+D^n_{i,\jph}\leq 1,\quad \mbox{ under the CFL condition (\ref{cfl}).}
\end{align*}
This proves the maximum principle. Other cases can be handled in a similar way and the maximum principle can be shown.
\subsection{Numerical experiments}\label{2d-expt}
For numerical simulation  we have chosen an example of the quarter five-spot problem in the domain  $[0,1] \times [0,1]$. 
To show the effect of gravity numerical experiments are performed  in the  presence of gravity as well as in the absence of gravity.
 Also to study the polymer flooding effect numerical experiments are performed for various concentration of the polymers.
 The behavior of water saturation is studied when the polymers are injected with different 
concentrations. The flux function $F=(F_1,F_2)$ takes the same form as in equation (\ref{eq:ff}) and (\ref{eq:ff1}) with 
  $$\lambda_w=\frac{s^2}{\mu_w(c_1,c_2)} , \, \lambda_o={(1-s)^2},\,  \rho_w g=2 \mbox{ and } \rho_o g=1,\,a_l(c_l)=1+0.5c_l\,(l=1,2)$$ and velocity $v$ 
  across the grid point is calculated by using 
(\ref{eq:p}).
\subsection{Initial and boundary conditions}
\begin{figure}
\begin{minipage}[b]{0.45\linewidth}
\centering
\includegraphics[scale=0.35]{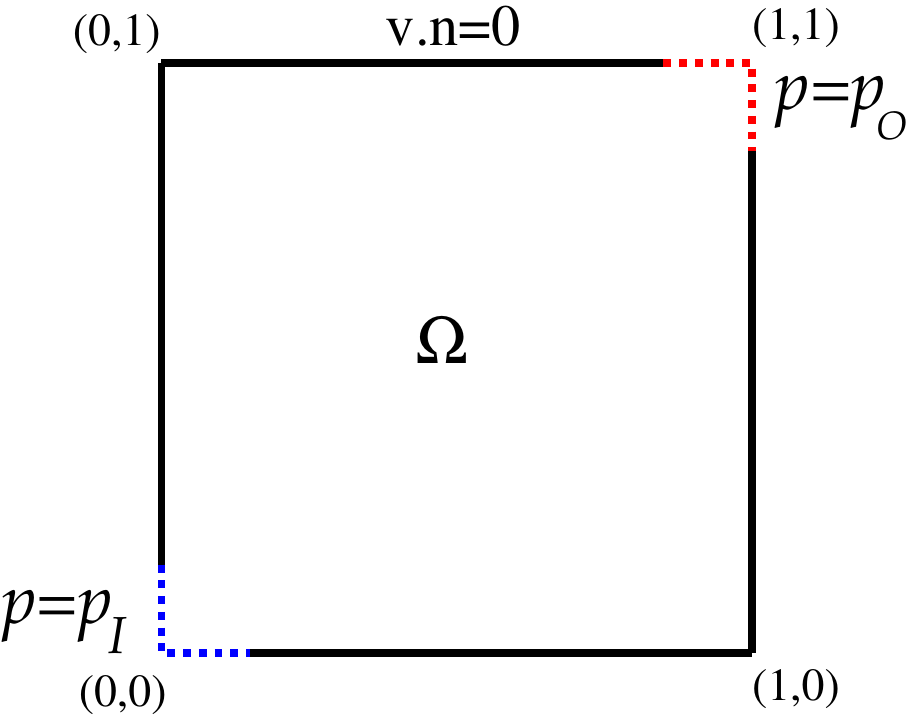}
\caption{Reservoir domain and boundary conditions}
\label{fig:bc}
\end{minipage}
\hspace{0.7cm}
\begin{minipage}[b]{0.45\linewidth}
\centering
\includegraphics[scale=0.24]{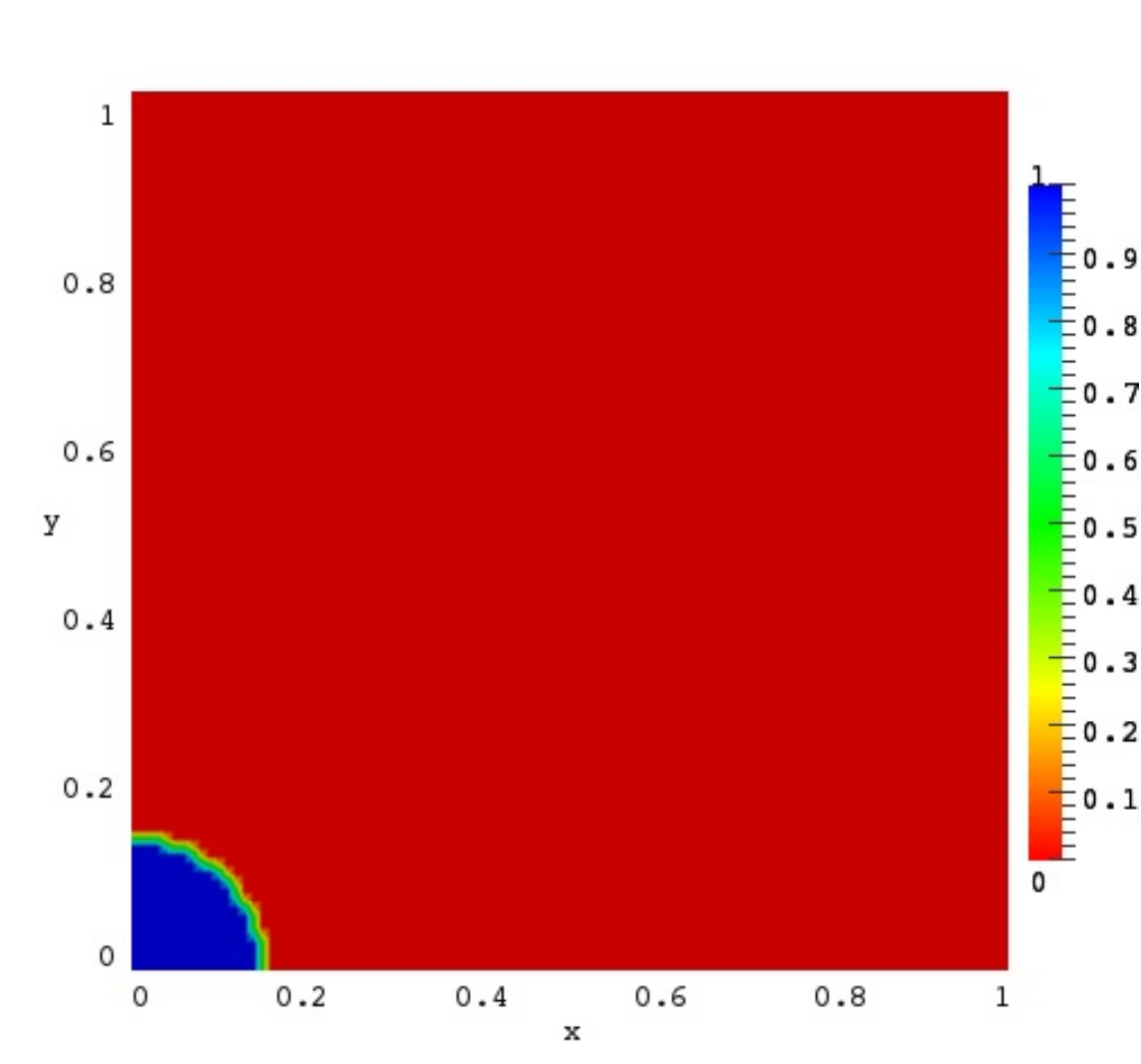}
 \caption{Pumping water through the inlet boundary}
 \label{ini}
\end{minipage}
\end{figure}

 The simulations are performed in a computational domain  $\Omega=[0,1] \times [0,1]$ for $t\in [0,1]$.
 The initial condition is $s(x,0)=0,$ i.e. In the inlet part of the boundary we pump water with 
a pressure $p=p_I$  and we keep the outlet part of the boundary with a pressure $p=p_O(p_I>p_O)$, on the remaining part of the boundary normal velocity is set to zero.
The initial inlet saturation is shown in Fig.\ref{fig:bc} and \ref{ini}
\subsection{Permeability of the porous media}
\begin{figure}
\centering
\subfigure[] {\includegraphics[scale=0.32]{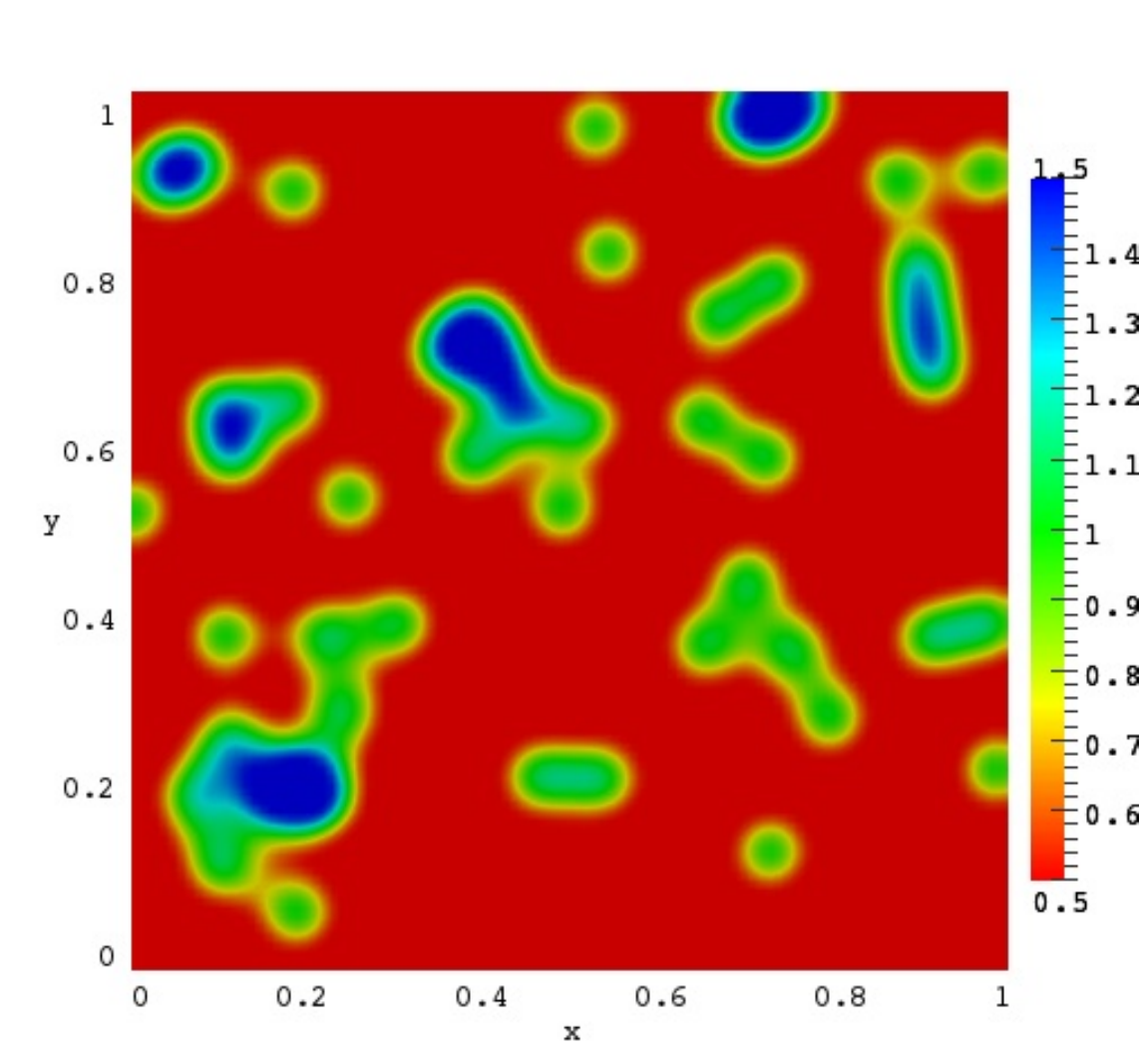}}
\subfigure[]{\includegraphics[scale=0.32]{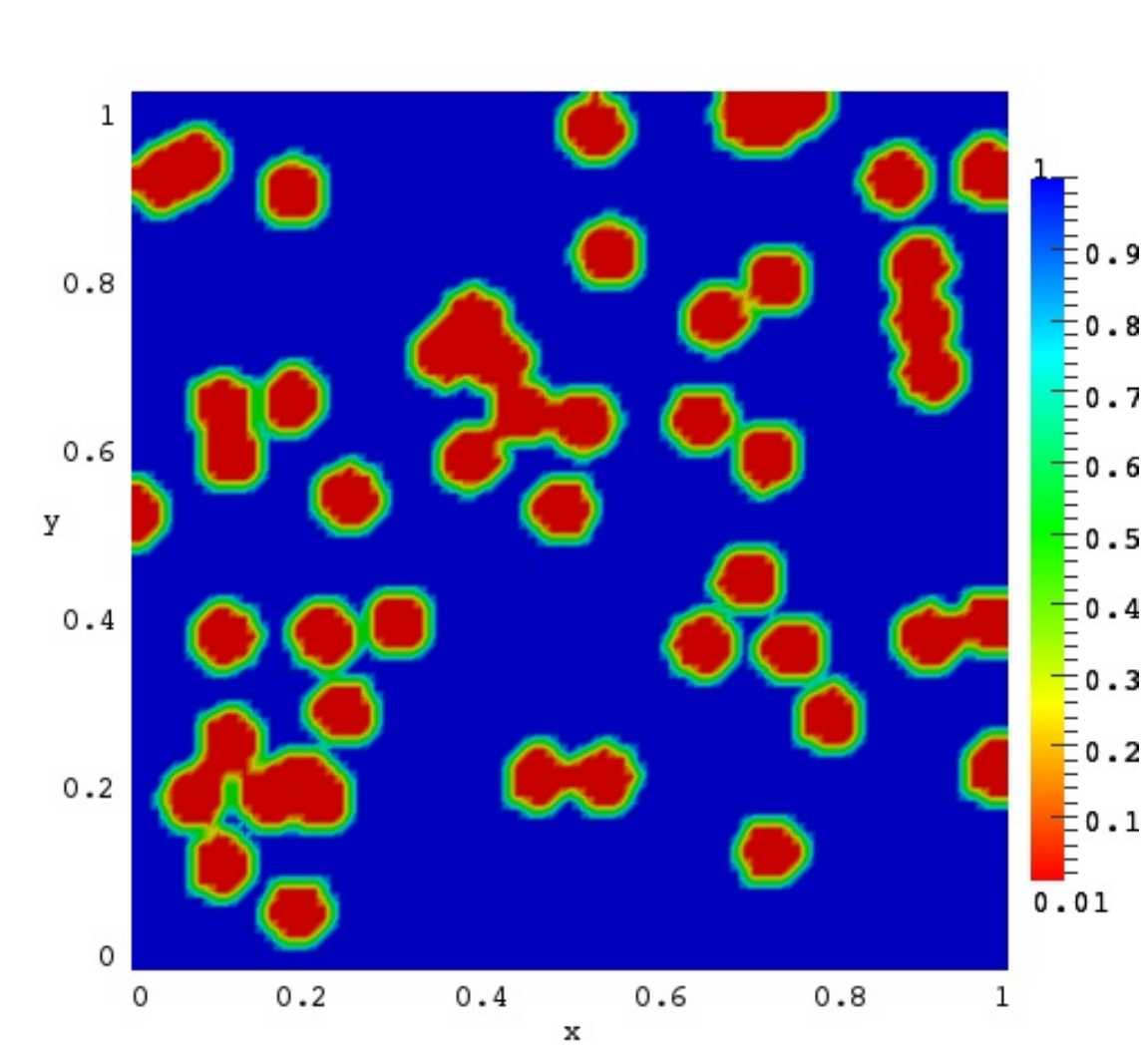}}
\caption{(a) Permeability fields for (\ref{pm1}), (b) Permeability fields for (\ref{pm2}) .}
\label{perm}
\end{figure}

We consider a heterogeneous porous medium with an absolute permeability $K(x).$ In order to illustrate the robustness of the proposed numerical scheme we
consider the two model porous media. The first test case corresponds to a heterogeneous medium with a continuous random  permeability given by 
\begin{equation}
\label{pm1}
K(x)= \min \{ \max \{ \sum\limits_{i=0}^N  \varPhi _i (x),0.5 \}, 1.5 \}
\end{equation}
and $$ \varPhi_i(x)=exp(-(\frac{|x-x_i|}{0.05})^2)$$ where $x_i $ are $N$ randomly chosen locations inside the domain. Here we have taken 
$N=100.$
The second test case corresponds to a heavily heterogeneous medium with hard rocks and the permeability is given by choosing $N$ random locations $x_i$
and
\begin{equation}
\label{pm2}
K(x)= \begin{cases}
0.01 & \textrm{if } x \in \mbox{B}(x_i,0.0015)\,\, \textrm{for some } i\in \{1,2, \,\,,, N\} \\
1 & \textrm{eslewhere } 
\end{cases}
\end{equation}
The permeability fields for these two test cases are shown in Fig.\ref{perm}
\par
\begin{figure}
\centering
\subfigure[] {\includegraphics[scale=0.35]{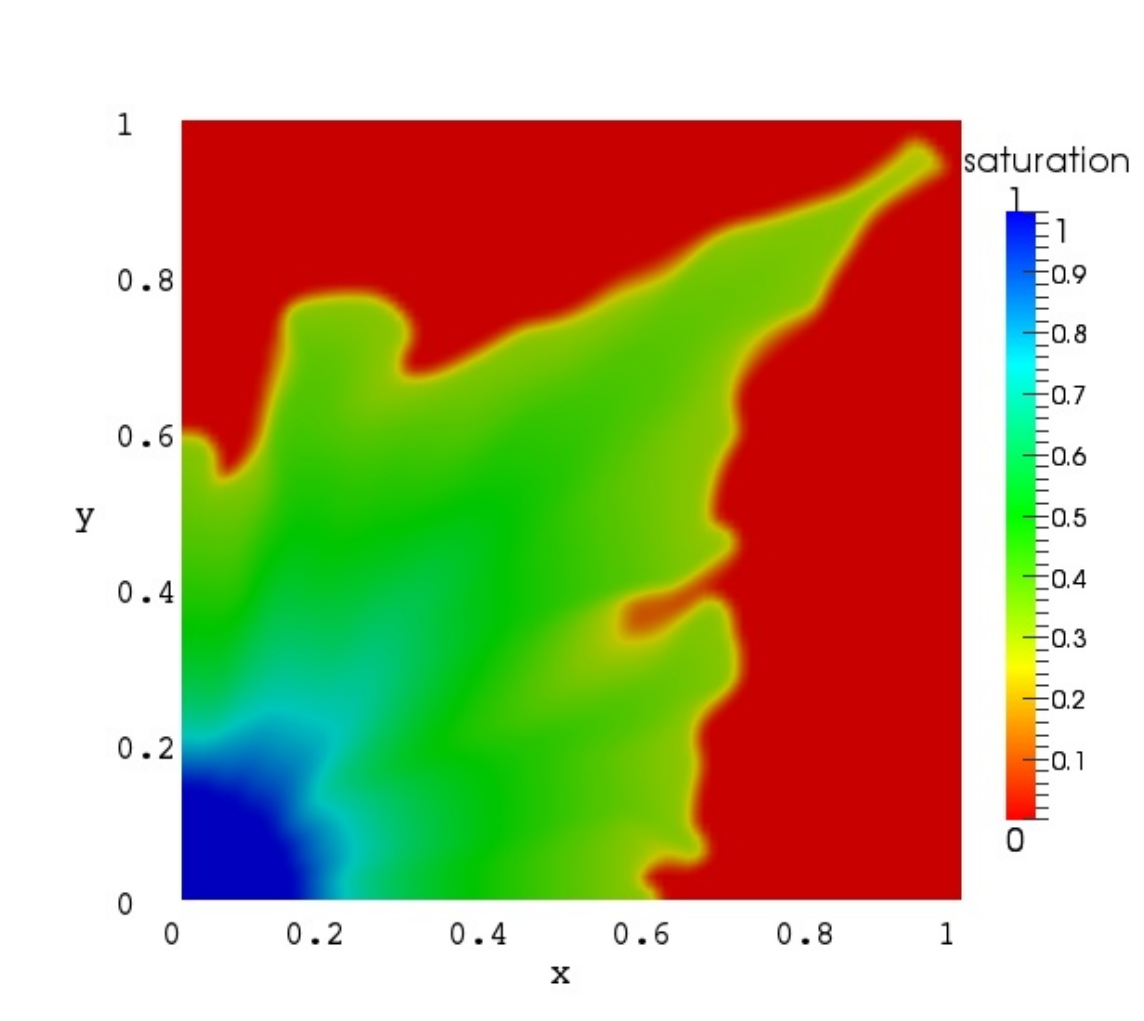}}
\subfigure[]{\includegraphics[scale=0.35]{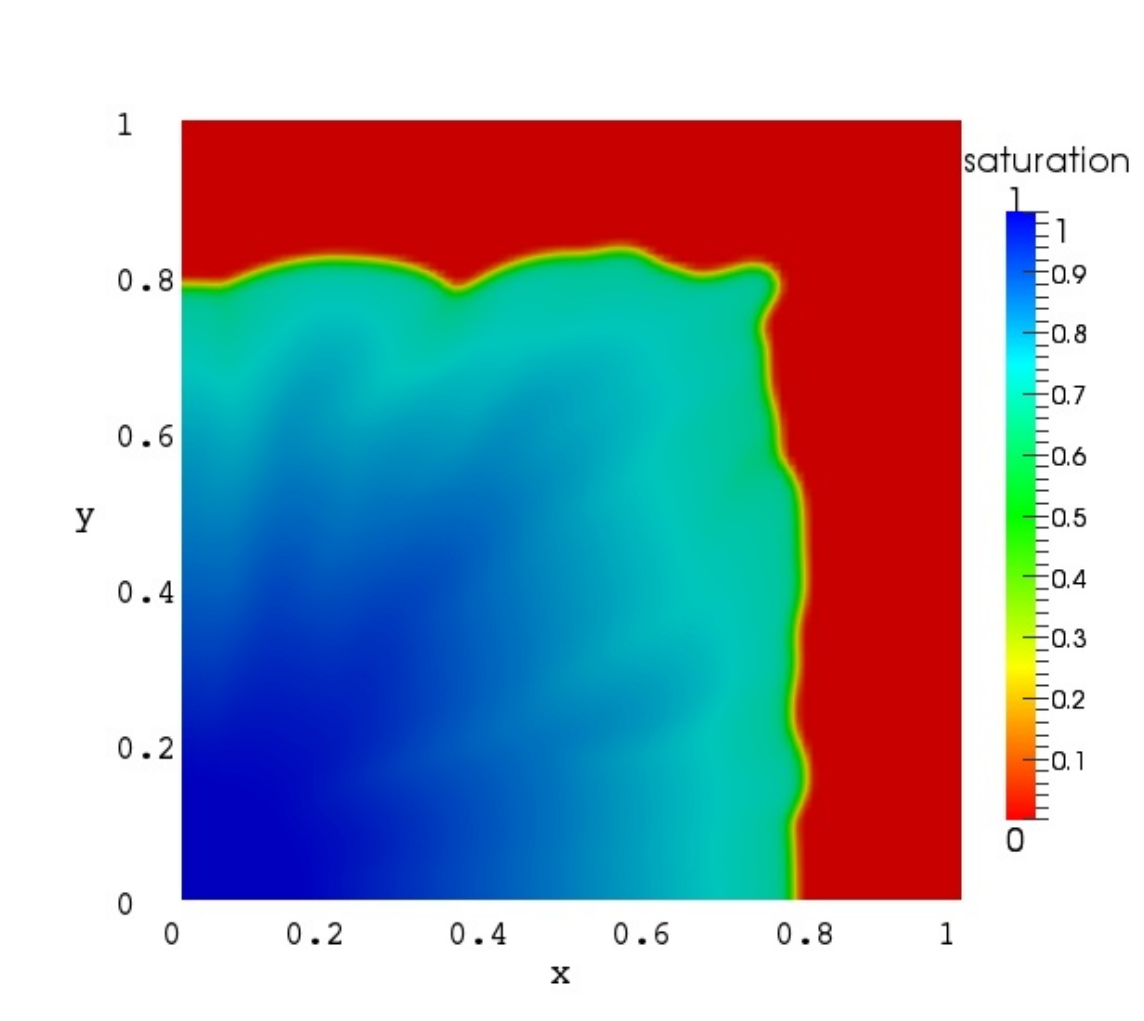}}
\caption{(a) Saturation s in the absence of  a polymer (b) Saturation s in the presence of a polymer.}
\label{fig1}
\end{figure}

{\bf{Experiment 1:}} Simulations in this experiment was performed using the spatial permeability distribution given in 
(\ref{pm1}), shown in Fig.\ref{perm}(a). The viscosity of water is given by $\mu_w(c_1,c_2)=0.5+c_1+c_2.$ We inject water through the inlet boundary with an inlet pressure $p_I=8$
and inlet concentration $c_1=0$ and $c_2=0.$ This is the case of without polymer. As expected it produces fingering effects, consequently when the water front touches the outlet boundary, a large
amount of oil is stuck in the remaining portion of the domain, which reduces the efficiency of oil-recovery, this is shown in Fig.\ref{fig1}(a).
To avoid this instability polymer is dissolved with water and injected through the inlet wall. In the presence of polymer say
 $c_1=7$ and $c_2=0 $ the fingering instability almost disappears and the amount of oil produced at the recovery well (outlet boundary) is increased. This is shown in
Fig.\ref{fig1}(b).
\par
\begin{figure}
\centering
\subfigure[] {\includegraphics[scale=0.35]{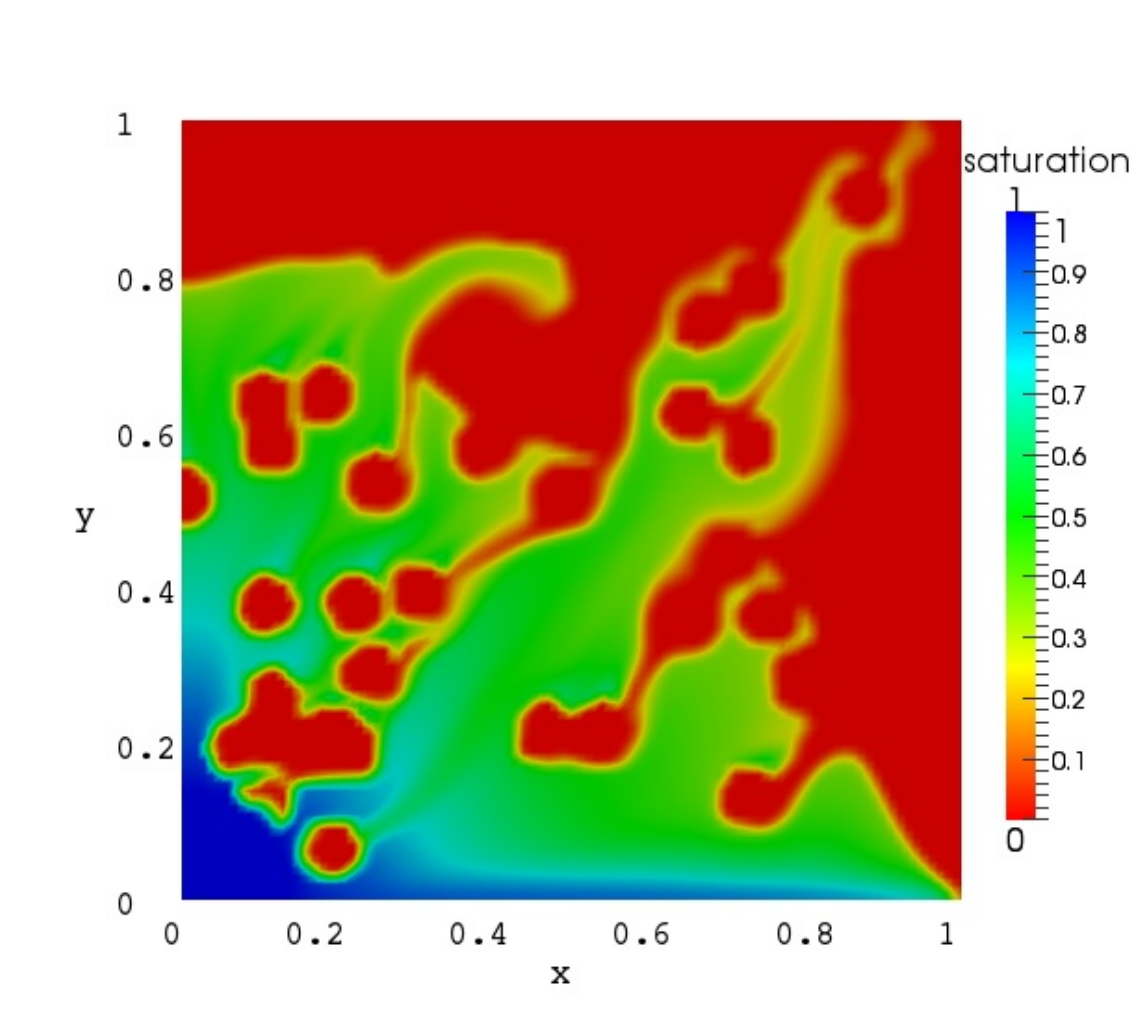}}
\subfigure[]{\includegraphics[scale=0.35]{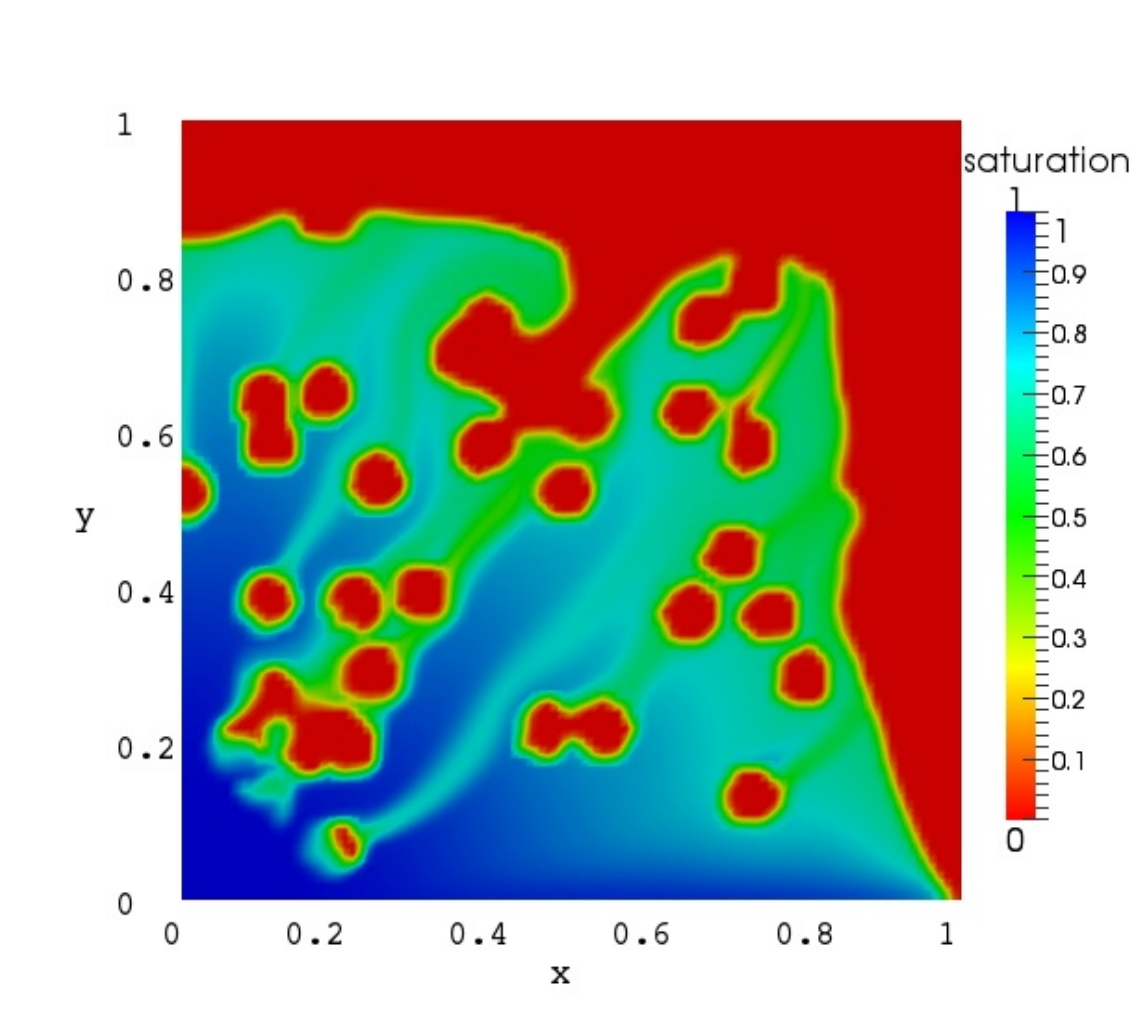}}
\caption{(a) Saturation s in experiment 2, in the abscence of polymer (b) Saturation s in experiment 2, in the presence of polymer.}
\label{fig2}
\end{figure}

{\bf{Experiment 2:}} 
The permeability fields is chosen as in expression (\ref{pm2}) with the presence of gravity, shown in Fig.\ref{perm}(b). Which corresponds to a heavily heterogeneous media with hard rocks.
 Here we have taken viscosity of water as  $\mu_w(c_1,c_2)=0.5+c_1+c_2.$ The result obtained in Fig.\ref{fig2}(a) corresponds to the saturation profile 
 with the inlet concentration $c_1=0$ and $c_2=0.$ The result obtained in Fig.\ref{fig2}(b) corresponds to the saturation profile with inlet concentration
 $c_1=5$ and $c_2=3. $ A consistent behavior of the saturation profile shows that our proposed scheme works well with varying spatial discontinuity in the media.
\par
\begin{figure}
\centering
\subfigure[] {\includegraphics[scale=0.35]{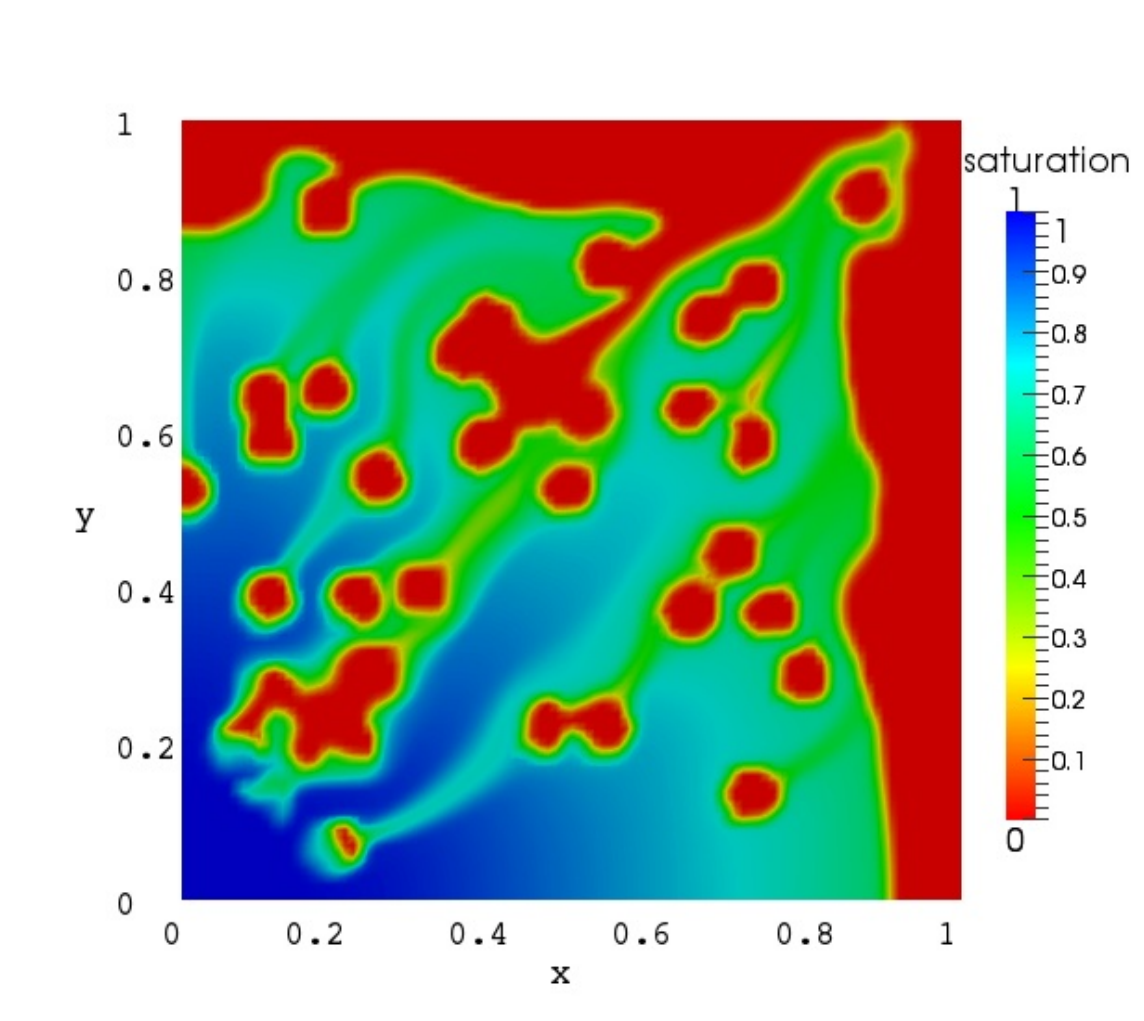}}
\subfigure[]{\includegraphics[scale=0.35]{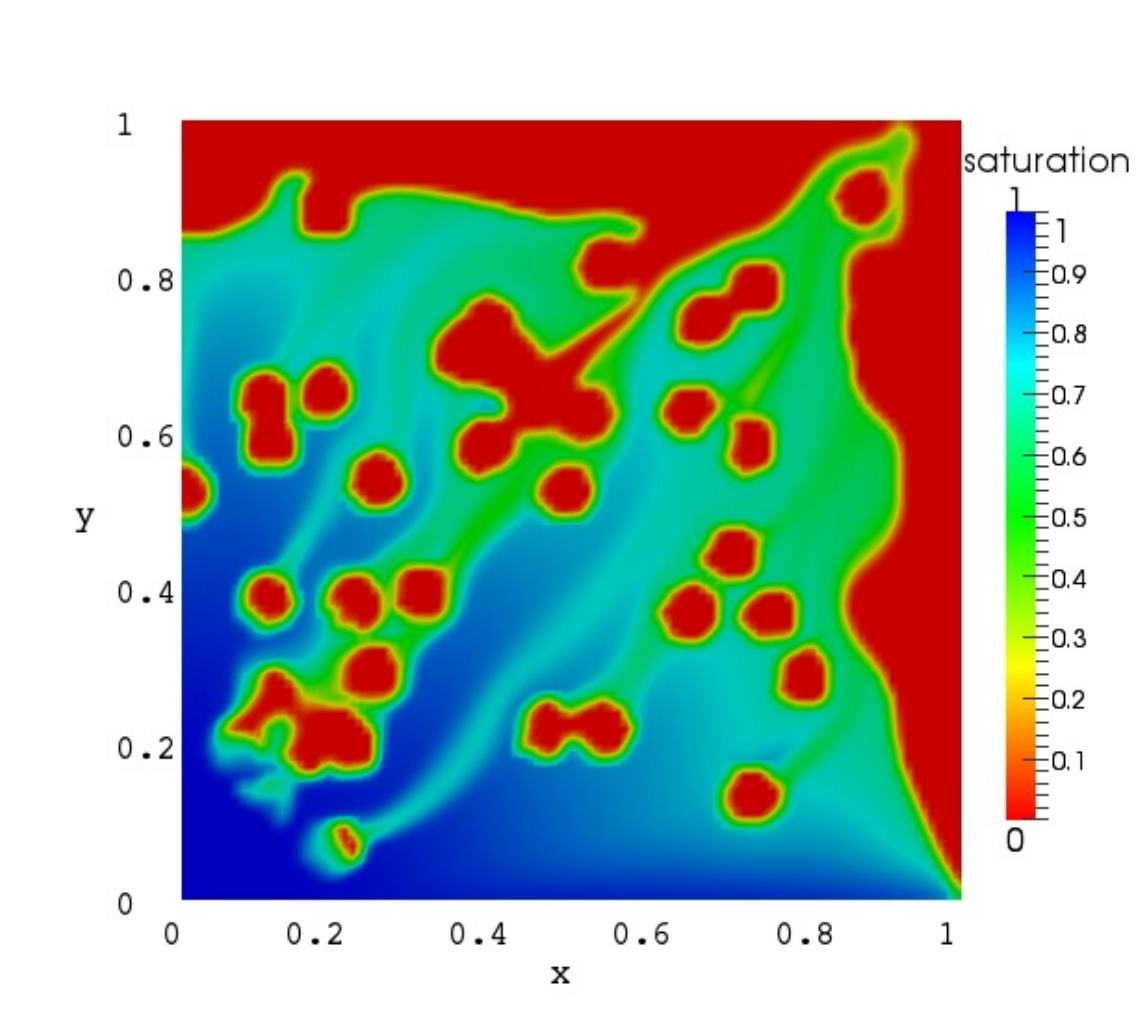}}
\caption{(a) With out the effect of gravity (b) With the effect of gravity}
\label{fig3}
\end{figure}

{\bf{Experiment 3:}}
This experiment is mainly to study the effect of gravity in  saturation profile. This experiment 
is performed using spatial permeability distributions given in (\ref{pm2}), shown in Fig.\ref{perm}(b). Viscosity takes the form 
$\mu_w(c_1,c_2)=0.5+c_1+c_2.$ We chose the inlet concentrations as $c_1=7$ and $c_2=0.$ The expression involving gravity term is considered along the
$y$ direction (see eqn (\ref{eq:ff1})). The resulting figures are shown in Fig.\ref{fig3}(a) with the absence  of gravity and Fig.\ref{fig3}(b) with presence of gravity. Observe that
presence of gravity significantly effects the saturation profile.
\par
\begin{figure}[H]
\centering
\subfigure[] {\includegraphics[scale=0.35]{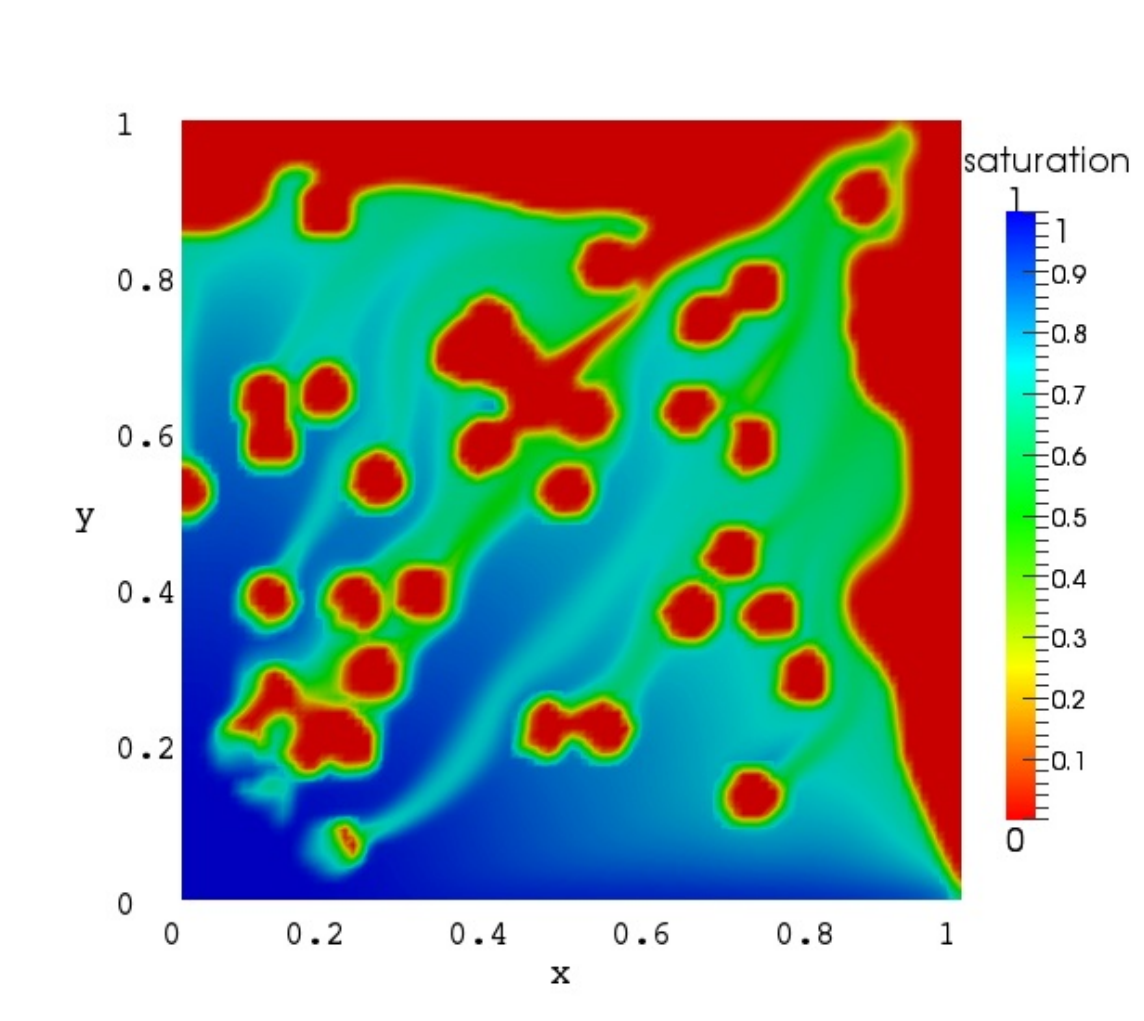}}
\subfigure[]{\includegraphics[scale=0.35]{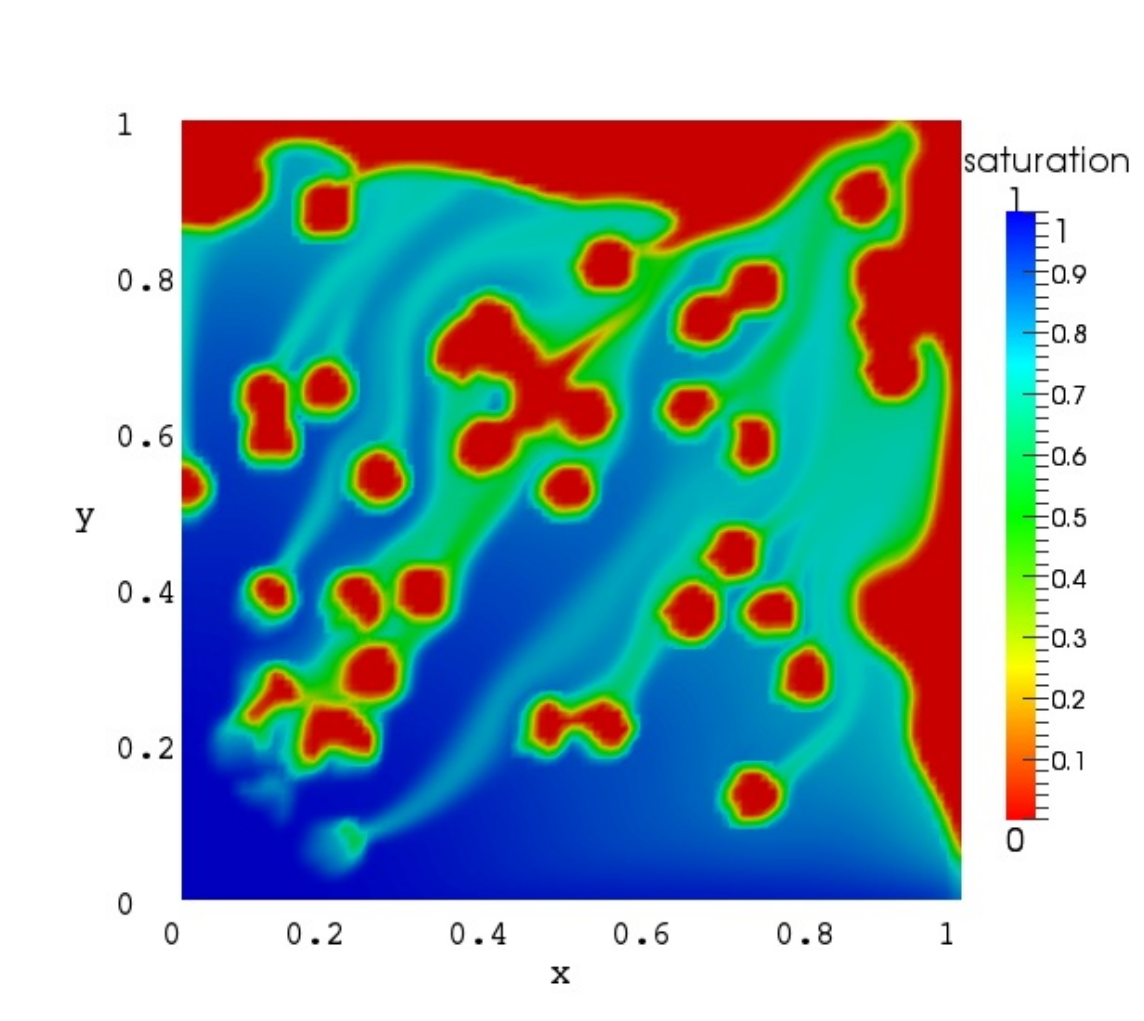}}
\caption{(a) Single component (b) Multicomponent}
\label{fig4}
\end{figure}

{\bf{Experiment 4:}} This experiment is to study the effect of adding more than one  polymer with different concentrations.
In this model we have taken $\mu_w(c_1,c_2)=0.5+\sqrt{c_1}+\sqrt{c_2}$ and permeability field is chosen as in (\ref{pm2}).
 Figure \ref{fig4}(a) corresponds to the case with concentrations $c_1=49,c_2=0.$ In figure \ref{fig4}(b) we have taken the concentrations to be $c_1=25,c_2=24.$ Observe that
 the total amount $(c_1+c_2) $ of injected concentrations in both the case are the same. But in the second case by adding two concentrations the sweeping profile of 
 water saturation is improved considerably. This is reflected in fig \ref{fig4}(b). It is clear from this fact that by adding multiple polymers
 and by taking a suitable viscosity $\mu_w(c_1,c_2)$ it may be possible to maximize the oil-recovery.
\par
\section{Conclusion.}
 A high resolution finite volume scheme is  developed to study the two-phase flow in porous media
by using the idea of discontinuous flux. 
The idea of discontinuous flux helps to reduce the system to an uncoupled scalar equation with discontinuous coefficients.
Discontinuous flux uses the solution of the Riemann problem of the scalar equation where as the Godunov flux needs solution of the
Riemann problem of the coupled system which is difficult to construct especially in the presence of gravity, heterogeneity and multiple components.
The results obtained from the idea of discontinuous flux agrees well with the results obtained from the Godunov flux.
The two-phase flow is studied in the presence as well as  in the absence  of gravity.
It is shown that the presence of gravity affects the saturation profile.
Also the efficiency of the numerical method is demonstrated by performing numerical simulations corresponding to 
two-phase flow in heterogeneous media.

\end{document}